\documentclass[10pt]{amsart}
\usepackage{amsmath, amssymb, amsthm, amscd, graphicx}

\theoremstyle{plain}
\newtheorem{prop}{Proposition}[section]
\newtheorem*{thrmA}{Theorem A}
\newtheorem*{thrmB}{Theorem B}
\newtheorem{coro}[prop]{Corollary}
\newtheorem{lemm}[prop]{Lemma}
\newtheorem{ques}[prop]{Question}
\newtheorem{thrm}[prop]{Theorem}

\theoremstyle{definition}
\newtheorem{defi}[prop]{Definition}
\newtheorem{nota}[prop]{Notation}
\newtheorem{exam}[prop]{Example}
\newtheorem{rema}[prop]{Remark}
\newtheorem*{rema*}{Remark}

\numberwithin{equation}{section}

\renewcommand{\a}{\alpha}
\renewcommand{\AA}{A}
\renewcommand{\aa}{a}%
\newcommand\AlgoD{\mathcal{A}}

\renewcommand{\b}{\beta}
\newcommand{\bb}{b}%
\newcommand{\BB}[1]{B_{#1}^{\scriptscriptstyle+}}
\newcommand{\BBB}[1]{\boldsymbol{B}_{\!\,#1}^{\scriptscriptstyle+}}
\newcommand\BSucc[2]{#2^{[#1]}}

\newcommand\card{\mathrm{card}}
\newcommand{\cc}{c}

\newcommand{\cl}[1]{[#1]}
\newcommand{\CL}[1]{\overline{\vrule height5pt width0pt#1}}
\newcommand{\conc}{\mathbin{\hbox to 2mm{$^{\scriptscriptstyle\frown}$}}}

\renewcommand{\d}{\delta}%
\let\D=\Delta
\newcommand{\dd}{d}
\newcommand{\DD}[2]{D\!_{{#1}}(#2)}
\newcommand{\ddd}[1]{\delta_{#1}}
\newcommand{\DDhat}[2]{\widehat{\Delta}_{#1, #2}}
\newcommand{\DDfe}[2]{\exp{D}_{#1}(#2)}
\newcommand{\DDf}[2]{D_{#1}(#2)}
\renewcommand{\div}{\prec}

\newcommand{\divL}{\prec_{_L}}
\newcommand{\DR}{\mathrm{Div}_{\scriptscriptstyle\! R}}

\newcommand{\e}{\varepsilon}
\newcommand{\ea}{{\scriptstyle\varnothing}}%
\newcommand\ee{e}
\newcommand\emin{e^{\scriptscriptstyle\mathtt{min}}}
\newcommand{\etc}{\emph{etc.}}
\renewcommand{\exp}[1]{{#1}^{*}}

\newcommand\ff[1]{\Phi_{#1}}

\let\ge=\geqslant
\newcommand{\gen}{g}
\newcommand{\gf}{>^{\scriptscriptstyle+}}
\newcommand{\gfe}{\ge^{\scriptscriptstyle+}}
\newcommand\Gen[1]{\gen_{#1}}
\newcommand\GGen{\boldsymbol{\gen}}
\newcommand{\gs}{>}
\newcommand{\gse}{\ge}

\newcommand{\h}{\theta}

\newcommand\ie{\emph{i.e.}}
\newcommand\ii{i}
\newcommand\II{I}
\newcommand\ince{\subseteq}
\newcommand\inv{^{-1}}

\newcommand{\jj}{j}
\newcommand{\JJ}{J}

\newcommand{\kk}{k}

\renewcommand{\l}{\lambda}
\let\le=\leqslant
\newcommand\Ldots{...\,}
\newcommand\lf{<^{\scriptscriptstyle+}}

\renewcommand{\lg}[1]{\vert#1\vert}
\newcommand{\lmin}{\mathrel{\widetilde{<}}}
\newcommand\ls{<}

\newcommand\lSL{<^{\mathtt{\scriptscriptstyle ShortLex}}}

\newcommand{\mm}{m}
\newcommand{\MM}{M}
\newcommand{\MMM}{\boldsymbol{M}}
\newcommand{\Mon}[1]{\langle#1\rangle^{\scriptscriptstyle\!+}}
\newcommand{\mult}{\succ}
\newcommand{\multe}{\succcurlyeq}

\newcommand{\Nat}{\mathbb{N}}
\newcommand{\nn}{n}
\newcommand{\nno}{{n-1}}
\newcommand{\nnt}{{n-2}}
\newcommand{\norm}[1]{\vert\!\vert#1\vert\!\vert}

\newcommand\om{\omega}
\newcommand\Ord{\mathbf{Ord}}
\newcommand\OVER[1]{\,#1\inv}

\newcommand{\perm}{\pi}
\newcommand{\pp}{p}
\newcommand{\ppo}{{p-1}}

\newcommand{\qq}{q}

\newcommand{\rr}{r}
\newcommand{\rro}{{r-1}}
\newcommand{\rrt}{{r-2}}
\def\resp{\emph{resp.}~}

\newcommand{\s}{\sigma}
\newcommand{\ShortLex}{\mathtt{ShortLex}}
\renewcommand{\ss}[1]{\sigma_{\!#1}^{\vrule height5pt
width0pt}}
\renewcommand{\SS}{S}
\newcommand{\sss}[1]{\sigma_{#1}^{-1}}
\newcommand\Succ[2]{#2^{(#1)}}

\newcommand{\Tail}[2]{\mathrm{tail}(#1,#2)}
\newcommand{\TTT}[2]{\underline{D}_{#1}(#2)}
\newcommand{\TTTe}[2]{\exp{\underline{D}}_{#1}(#2)}

\newcommand{\uu}{u}
\newcommand{\uuu}{\boldsymbol{s}}

\newcommand{\vv}{v}
\newcommand{\vvv}{\boldsymbol{t}}

\newcommand{\ww}{w}
\newcommand{\WW}[1]{\underline{B}_{#1}^{\scriptscriptstyle+}}
\newcommand{\www}{\boldsymbol{s}}
\newcommand{\WWW}[1]{\underline{\boldsymbol{B}}_{\!\,#1}^{\scriptscriptstyle+}}

\newcommand{\xs}{s}
\newcommand{\xt}{t}
\newcommand{\xx}{x}
\newcommand{\XX}{X}

\newcommand{\yy}{y}

\newcommand{\zz}{z}
\begin{document}

\hfill{\tiny 2008-02}

\author{Patrick DEHORNOY}
\address{Laboratoire de Math\'ematiques Nicolas Oresme,
Universit\'e de Caen, 14032 Caen, France}
\email{dehornoy@math.unicaen.fr}
\urladdr{//www.math.unicaen.fr/\!\hbox{$\sim$}dehornoy}

\title{Alternating normal forms for braids and locally
Garside monoids}



\begin{abstract}
We describe new types of normal forms for braid
monoids, Artin--Tits monoids, and, more generally, for all
monoids in which divisibility has some convenient lattice
properties (``locally Garside monoids''). We show that, in
the case of braids, one of these normal forms coincides with the
normal form introduced by Burckel and  deduce that the latter
can be computed easily. This approach leads to a new, simple
description for the standard ordering (``Dehornoy order'')
of~$B_\nn$ in terms of that of~$B_\nno$, and to a quadratic
upper bound for the complexity of this ordering.
\end{abstract}
\maketitle

The first aim of this paper is to improve our understanding
of the well-ordering of positive braids and of the Burckel normal form 
of~\cite{BuT, Bur}, which after more than ten years remain
mysterious objects. This aim is achieved, at least
partially, by giving  a new, alternative definition for the
Burckel normal form that makes it natural and easily
computable. This new description is direct, involving right
divisors only, while Burckel's original approach resorts to
iterating some  tricky reduction procedure. It turns out that
the construction we describe below relies on a very general
scheme for which many monoids are eligible, and we hope for
further applications beyond the case of braids.

After the seminal work of F.A.\,Garside~\cite{Gar},
we know that braid monoids and, more generally, spherical Artin--Tits
monoids and Garside monoids that generalize them, can be 
equipped with a normal form, namely the so-called
greedy normal form of~\cite{BrS, Adj, ElM, Thu}, which gives
for each element of the monoid a distinguished representative
word. This normal form is excellent both in theory and in
practice as it provides a bi-automatic structure and it is
easily computable~\cite{Eps, Cha, Dgk}.

What we do in this paper is to construct a new type of
normal form for braid monoids and their generalizations.
Our construction keeps one of the ingredients of the (right)
greedy normal form, namely considering the maximal right
divisor that lies in some subset~$\AA$, but, instead of taking for~$\AA$ the
set of so-called simple elements, \ie, the divisors of the
Garside element~$\D$, we choose $\AA$ to be some
standard parabolic submonoid~$\MM_\II$ of~$\MM$,
\ie, the monoid generated by some subset~$\II$ of the
standard generating set~$\SS$. When $\II$ is a proper
subset of~$\SS$, the submonoid~$\MM_\II$ is a proper subset
of~$\MM$, and the construction stops after one step.
However, by considering two parabolic submonoids~$\MM_\II$, $\MM_\JJ$ which
together generate~$\MM$, we can obtain a well-defined, unique
decomposition consisting of alternating factors in~$\MM_\II$
and~$\MM_\JJ$, as in the case of an amalgamated product. By considering
convenient families of  submonoids, we can iterate the process and
obtain a unique normal form for each element of~$\MM$.
When it exists, typically in all Artin--Tits monoids, such a
normal form is exactly as easy to compute as the greedy
normal form, and, as the greedy form, it solves the word problem in quadratic time. 

The above construction is quite general, as it only
requires the ground monoid~$\MM$ to be what is
now called locally right Garside---or locally left
Gaussian in the obsolete terminology of~\cite{Dgl}.
However, our main interest in the current paper lies in the case of braids and, more specifically, their ordering. For a
convenient choice of the parameters, the alternating normal form  turns out to
coincide with the Burckel normal form of~\cite{Bur}. As a
consequence, we at last obtain both an easy algebraic
description of the latter, and an efficient algorithm for
computing it. Mainly, because of the connection between the
Burckel normal form and the standard ordering of braids
(``Dehornoy order''), we obtain a new characterization of the
latter. The result can be summarized as follows. As usual,
$\BB\nn$ denotes the monoid of positive $\nn$-strand
braids. We use $\ff\nn$ for the involutive flip
automorphism of~$\BB\nn$ that maps~$\ss\ii$
to~$\ss{\nn-\ii}$ for each~$\ii$, \ie, for conjugation by the Garside element~$\D_\nn$, and $<$ for the upper
version of the standard braid ordering.

\begin{thrmA}
$(i)$ Every positive $\nn$-strand braid~$\xx$ admits a unique decomposition
$$\xx= \ff\nn^\ppo(\xx_\pp) \cdot ...  \cdot \ff\nn^2(\xx_3) \cdot \ff\nn(\xx_2) \cdot \xx_1$$ 
with $\xx_\pp, \Ldots, \xx_1$ in~$\BB\nno$ such that, for
each $\rr \ge2$, the only generator~$\ss\ii$ that divides
$\ff\nn^{\pp-\rr}(\xx_\pp) \cdot ... \cdot \ff\nn(\xx_{\rr+1}) \cdot
\xx_{\rr}$ on the right is~$\ss1$. Starting from~$\xx^{(0)} = \xx$, the
element~$\xx_\rr$  is determined  by the condition that $\xx_\rr$ is the
maximal right divisor of~$\xx^{(\rr)}$ that lies in~$\BB\nno$, and
$\xx^{(\rr)}$ is $\ff\nn(\xx^{(\rr-1)}\xx_\rr\inv)$.

$(ii)$ Let $\xx, \yy$ be positive $\nn$-strand braids. Let $(\xx_\pp, \Ldots, \xx_1)$ and $(\yy_\qq,\Ldots,\yy_1)$
be the sequences associated with~$\xx$ and~$\yy$ as in~$(i)$. Then $\xx<\yy$ holds in~$\BB\nn$ if and only if we
have either $\pp<\qq$, or $\pp=\qq$ and, for some~$\rr \le \pp$, we have $\xx_{\rr'}=\yy_{\rr'}$ for
$\rr < \rr' \le \pp$ and $\xx_\rr<\yy_\rr$ in~$\BB\nno$.
\end{thrmA}

In other words, via the above decomposition, the ordering of~$\BB\nn$ is a $\ShortLex$-extension of that
of~$\BB\nno$, this meaning the variant of lexicographical extension in which the length is given priority.
In the above  statement, Point~$(i)$---Proposition~\ref{P:FSplitting} below---is easy, but
Point~$(ii)$---Corollary~\ref{C:RecOrder}---is not. Another outcome of the
current approach is the following complexity upper
bound for the braid ordering---Corollary~\ref{C:Complexity}:

\begin{thrmB}
For each~$\nn$, the standard ordering of~$B_\nn$ has at most a quadratic
complexity: given two $\nn$-strand braid words~$\uu, \vv$ of
length~$\ell$, we can decide whether the braid represented by~$\uu$ is
smaller than the braid represented by~$\vv$ in time~$O(\ell^2)$.
\end{thrmB}

We think that the tools developed in this paper might be useful for addressing other types of questions, typically
those involving conjugacy in~$B_\nn$.

The paper is organized as follows. In Section~\ref{S:Alt}, we describe the alternating decompositions obtained when
considering two submonoids in a locally Garside monoid. In Section~\ref{S:Iter}, we show how to iterate the
construction using a binary tree of nested submonoids. In Section~\ref{S:NormalForm}, we deduce a
normal form result in the case when the base submonoids are generated by atoms. From Section~\ref{S:Flip}, we
concentrate on the specific case of braids and investigate what we call the $\Phi$-splitting and the $\Phi$-normal
form of a braid. Finally,  in Section~\ref{S:Burckel}, we investigate the connection between the $\Phi$-normal form
and the Burckel normal form, and deduce the above mentioned applications to the braid ordering.

\begin{rema*}
All constructions developed in this paper involve right
divisibility and the derived notions. This choice is dictated
by the braid applications of Section~\ref{S:Burckel}. Of
course, we could use left divisibility instead and obtain
symmetric versions, in the framework of  monoids that are
locally Garside on the left.
\end{rema*}

We use~$\Nat$ for the set of all nonnegative integers.

\section{Alternating decompositions}
\label{S:Alt}

We construct unique decompositions for the elements of
monoids in which enough least common left multiples (left
lcm's) exist. If $\MM$ is such a monoid
and $\AA$ is a subset of~$\MM$ that is closed under the left lcm
operation, then, under weak additional assumptions, every
element~$\xx$ admits a distinguished decomposition $\xx= \xx'
\xx_1$, where $\xx_1$ is a maximal right divisor of~$\xx$ that
lies in~$\AA$. The element~$\xx_1$ will be called the
$\AA$-tail of~$\xx$. If we assume that every non-trivial (\ie,
$\not=1$) element of~$\MM$ has a non-trivial $\AA$-tail, we can
consider the $\AA$-tail of~$\xx'$, and, iterating the process,
obtain a distinguished decomposition of~$\xx$ as a product of
elements of~$\AA$, as done for the standard greedy normal
form of Garside monoids. Here, we drop the assumption
that every non-trivial element has a non-trivial $\AA$-tail, but
instead consider two subsets~$\AA_1, \AA_2$ of~$\MM$ with
the property that, for every non-trivial~$\xx$, at least one of
the $\AA_1$- or $\AA_2$-tails of~$\xx$ is non-trivial. Then, we
obtain a distinguished decomposition of~$\xx$ as an alternating
product of elements of~$\AA_1$ and of~$\AA_2$.

\subsection{Locally Garside monoids}

Divisibility features play a key r\^ole throughout 
the paper, and we first fix some notation. 

\begin{nota}
For $\MM$ a monoid and $\xx,\yy\in\MM$, we say that
$\yy$ is a \emph{right divisor} of~$\xx$, or, equivalently,
that $\xx$ is a \emph{left multiple} of~$\yy$, denoted
$\xx\multe\yy$, if $\xx=\zz\yy$ holds
for some~$\zz$; we write $\xx\mult\yy$ if $\xx=\zz\yy$
holds for some~$\zz\not=1$. The set of all right divisors
of~$\xx$ is denoted by~$\DR(\xx)$. 
\end{nota}

The approach considered below turns out to be relevant
for the following monoids.

\begin{defi}
\label{D:LocGarside}
We say that a monoid~$\MM$ is a \emph{locally right Garside} if:

$(C_1)$ The monoid $\MM$ is right cancellative, \ie,
$\xx\zz=\yy\zz$ implies $\xx=\yy$;

$(C_2)$ Any two elements of~$\MM$ that admit a
common left multiple admit a left~lcm; 

$(C_3)$ For every~$\xx$ in~$\MM$, there is no infinite
chain $\xx_1 \div \xx_2 \div ... $ in~$\DR(\xx)$.
\end{defi}

If $\MM$ is a locally right Garside monoid, and
$\xx,\yy$ are elements of~$\MM$ satisfying
$\xx\multe\yy$, the element~$\zz$ satisfying $\xx=\zz\yy$
is unique by right cancellativity, and we denote it
by~$\xx\OVER{\yy}$.

\begin{exam}
\label{X:Artin}
According to~\cite{BrS} and~\cite{Mic}, all Artin--Tits
monoids are locally right (and left) Garside. We recall that
an Artin--Tits monoid is a monoid generated by  a
set~$\SS$ and relations of the form
$\xs\xt\xs...=\xt\xs\xt...$ with $\xs,\xt \in \SS$, both sides of the same length,
and at most one such relation for each pair~$\xs,\xt$. An important
example is Artin's braid monoid~$\BB\nn$ \cite{KaT}, which corresponds to
$\SS=\{\ss1, \Ldots, \ss\nno\}$ with
\begin{equation}
\label{Pres}
\ss i \ss j = \ss j \ss i
\quad\mbox{for $\vert i - j\vert \ge 2$,}
\qquad
 \ss i \ss j \ss i = \ss j \ss i \ss j 
\quad\mbox{for $\vert i - j\vert = 1$}.
\end{equation}
As the name suggests, more general examples of locally
Garside monoids are the Garside monoids
of~\cite{Dfx, Dgk, CMW, ChM, Pid}, which include torus
knot monoids~\cite{Pik}, dual braid monoids~\cite{BKL},
and many more.
\end{exam}

If $\MM$ is locally right Garside, then no non-trivial
 element of~$\MM$ is invertible: if we
had $\xx\yy=1$ with $\xx\not=1$, hence $\yy\not=1$, the
sequence $\xx, 1, \xx, 1,...$ would
contradict~$(C_3)$. So right divisibility
is antisymmetric, and, therefore, it is a
partial ordering on~$\MM$. As a consequence, the left lcm, when it exists, is
unique.

Definition~\ref{D:LocGarside}---which also appears
in~\cite{DiM}---is satisfactory in that it exclusively involves
the right divisibility relation, and it directly leads to
Lemma~\ref{L:Lcm} below. Actually, it does not coincide with
the definitions of~\cite{Dgk} and~\cite{Dfx}, where
$(C_3)$ is replaced with some condition involving left
divisibility. However, both definitions are equivalent. For a 
while, we use $\divL$ for the proper left divisibility relation,
\ie, $\xx\divL\yy$ means $\yy=\xx\zz$ with $\zz\not=1$. We
denote by~$\Ord$ the class of ordinals.

\begin{lemm}
\label{L:LocGarside}
$(i)$ If $\MM$ is right cancellative,
Condition~$(C_3)$ is equivalent~to

$(C'_3)$ There is no infinite descending
chain in~$(\MM,\divL)$.

$(C''_3)$ There exists $\l: \MM\to\Ord$ such that $\yy\not=1$
implies $\l(\xx\yy) > \l(\xx)$.

\noindent $(ii)$ In any case, Conditions~$(C_3)$--$(C''_3)$
follow from

$(C_3^+)$ There exists $\l\!:\!\MM\to\Nat$
such that $\yy\not=1$ implies
$\l(\xx\yy)\ge\l(\xx)+\l(\yy) >\nobreak \l(\xx)$.
\end{lemm}

\begin{proof}
$(i)$ Assume that $\MM$ is right cancellative and
$(C_3)$ fails in~$\MM$. There exists~$\xx$
in~$\MM$ and a sequence $\xx_1, \xx_2,...$ in~$\DR(\xx)$
such that
$\xx_{\nn+1}\mult \xx_\nn$ holds for every~$\nn$. So,
for each~$\nn$, there exists~$\yy_\nn\not=1$ satisfying
$\xx_{\nn+1}=\yy_\nn\xx_\nn$. On the other hand, as
$\xx_\nn$ belongs to~$\DR(\xx)$, there exist~$\zz_\nn$
satisfying $\xx=\zz_\nn\xx_\nn$. We find
$$\xx=\zz_\nn\xx_\nn=\zz_{\nn+1}\xx_{\nn+1}=
\zz_{\nn+1}\yy_\nn\xx_\nn.$$
By cancelling~$\xx_\nn$
on the right, we deduce $\zz_\nn=\zz_{\nn+1}\yy_\nn$,
hence $\zz_{\nn+1}\divL\zz_\nn$ for each~$\nn$, and the
sequence $\zz_0, \zz_1, ...$ witnesses that $(C'_3)$ fails.

Conversely, assume that $(C'_3)$ fails in~$\MM$.
Let $\zz_0, \zz_1, ...$ be a descending chain for~$\divL$. For
each~$\nn$, choose~$\yy_\nn\not=1$ satisfying
$\zz_\nn=\zz_{\nn+1}\yy_\nn$. Let $\xx=\zz_0$,
$\xx_1=1$, and, inductively,
$\xx_{\nn+1}=\yy_\nn\xx_\nn$. By construction, we
have $\xx_{\nn+1} \mult\xx_\nn$ for each~$\nn$.
Now, we also have $\xx=\zz_\nn\xx_\nn$ for each~$\nn$,
so all elements~$\xx_\nn$ belong to~$\DR(\xx)$, and the sequence
$\xx_1, \xx_2, ...$ witnesses that $(C_3)$ fails.

The equivalence of~$(C'_3)$ and~$(C''_3)$ is standard, and
$(C_3^+)$ strengthens~$(C''_3)$.
\end{proof}

Condition~$(C_3^+)$ holds in particular in every monoid
that is presented by homogeneous relations, \ie, relations of
the form~$\uu=\vv$ where $\uu$ and~$\vv$ are words of
the same length: then we can define~$\l(\xx)$ to
be the length of any word representing~$\xx$. This is the
case for the Artin--Tits monoids of Example~\ref{X:Artin}.

Lemma~\ref{L:LocGarside} implies that locally right Garside
monoids coincide with the monoids called locally left Gaussian
in~\cite{Dgk}, in connection with the left Gaussian monoids
of~\cite{Dfx}. The reason for changing terminology is that
the current definition is coherent with~\cite{DiM} and it is
more natural: locally right Garside monoids involve right
divisibility, and the normal forms we discuss below are
connected with what is usually called the right normal form.

Assume that $\MM$ is a locally right Garside monoid.
Condition~$(C_2)$ is equivalent to saying that, for
every~$\xx$ in~$\MM$, any two elements of~$\DR(\xx)$
admit a left lcm, and it follows that any finite subset
of~$\DR(\xx)$ admits a global left lcm. By the Noetherianity
condition~$(C_3)$, the result extends to arbitrary subsets.
We say that a set~$\XX$ is \emph{closed under left lcm} if the
left lcm of any two elements of~$\XX$ exists and lies
in~$\XX$ whenever it exists in~$\MM$, \ie, by~$(C_2)$,
whenever these elements admit a common left multiple
in~$\MM$.

\begin{lemm}
\label{L:Lcm}
Assume that $\MM$ is a locally right Garside monoid, and
$\xx\in\MM$. Then every nonempty subset~$\XX$ of~$\DR(\xx)$ admits
a global left lcm~$\xx_1$; if moreover $\XX$ is closed
under left lcm, then $\xx_1$ belongs to~$\XX$.
\end{lemm}

\begin{proof}
Assume first that $\XX$ is closed under left lcm. By the
axiom of dependent choices, Condition~$(C_3)$ implies
that $(\DR(\xx),\mult)$ is a well-founded poset, so
$\XX$ admits some $\mult$-minimal, \ie, some
$\div$-maximal, element~$\xx_1$: so $\xx'\xx_1\in\XX$
implies~$\xx'=1$.  Then $\xx_1$ is a global left lcm for~$\XX$.
Indeed, assume $\yy_1 \in \XX$. By hypothesis,
$\xx_1$ and~$\yy_1$ lie in~$\DR(\xx)$, so, by~$(C_2)$, they
admit a left lcm~$\zz$, which can be expressed
as~$\zz = \yy'\yy_1=\xx'\xx_1$. The hypothesis that $\XX$ is
closed under left lcm implies $\zz\in\XX$. The choice
of~$\xx_1$ implies $\xx'=1$, hence $\xx_1 \multe \yy_1$.

If the assumption that $\XX$ is closed under left
lcm is dropped, we can apply the above result to the
closure~$\widehat\XX$ of~$\XX$ under left lcm. Then the global left
lcm~$\xx_1$ of~$\widehat\XX$ is a global left lcm for~$\XX$, but we
cannot be sure that
$\xx_1$ lies in~$\XX$---yet it is certainly the left lcm of some
finite subset of~$\XX$.
\end{proof}

Although standard, the previous result is crucial. By applying
Lemma~\ref{L:Lcm} to the subset $\DR(\xx)\cap\DR(\yy)$
of~$\DR(\xx)$, we deduce that any two elements~$\xx,\yy$
of a locally right Garside monoid~$\MM$ admit a right gcd
(greatest common divisor), and, therefore, for every~$\xx$
in~$\MM$, the structure $(\DR(\xx),\multe)$ is a lattice,
with minimum~$1$ and maximum~$\xx$. 

\subsection{The $\AA$-tail of an element}

If $\MM$ is a  monoid and $\xx, \xx_\pp, \Ldots, \xx_1$ belong to~$\MM$, we say
that $(\xx_\pp, \Ldots, \xx_1)$ is a \emph{decomposition} of~$\xx$ if $\xx = \xx_\pp
... \xx_1$ holds. The basic observation is that,  for each subset~$\AA$
of the monoid~$\MM$ that contains~$1$ and is closed under left
lcm, and every~$\xx$ in~$\MM$, Lemma~\ref{L:Lcm} leads to a distinguished
decomposition~$(\xx', \xx_1)$ of~$\xx$ with~$\xx_1 \in \AA$. 

\begin{lemm}
\label{L:Tail}
Assume that $\MM$ is a locally right Garside monoid and
$\AA$ is a  subset of~$\MM$ that contains~$1$ and is closed under left lcm.
Then, for each element~$\xx$ of~$\MM$, there exists a
unique right divisor~$\xx_1$ of~$\xx$ that lies in~$\AA$
and is maximal with respect to right divisibility, namely the
left lcm of $\DR(\xx)\cap\AA$.
\end{lemm}

\begin{proof}
Apply Lemma~\ref{L:Lcm} with $\XX=\DR(\xx)\cap\AA$.
The latter set  is nonempty as it contains at least~$1$, and it is closed under left lcm as it is the
intersection of two sets that are closed under left lcm.
\end{proof}

\begin{defi}
Under the hypotheses of Lemma~\ref{L:Tail}, the
element~$\xx_1$ is called the \emph{$\AA$-tail} of~$\xx$,
and denoted~$\Tail\xx\AA$.
\end{defi}

\begin{exam}
\label{X:Closed}
Let $\MM$ be an Artin--Tits monoid with standard set of
generators~$\SS$. We assume in addition that $\MM$ is of
spherical type, which means that the Coxeter group
obtained by adding the relation~$\xs^2=1$ for each~$\xs$ in~$\SS$ is finite.
Then, Garside's theory shows that any two elements of~$\MM$ 
admit a common left multiple, hence a left lcm. We shall
consider two types of closed subsets of~$\MM$. The first,
standard choice consists in considering the set~$\Sigma$ of
so-called simple elements in~$\MM$, namely the divisors
of the lcm~$\D$ of~$\SS$. By construction,
$\Sigma$ contains~$1$ and is closed under left (and right) divisor, and under
left (and right) lcm. For each~$\xx$
in~$\MM$, the $\Sigma$-tail of~$\xx$ is the right gcd
of~$\xx$ and~$\D$.

A second choice consists in considering $\II \subseteq \SS$, and taking
for~$\AA$ the standard parabolic submonoid~$\MM_\II$ of~$\MM$ generated
by~$\II$. The specific form of the Artin--Tits relations
implies that $\MM_\II$ is closed under left (and right)
divisor, and under left (and right) lcm, hence it is
eligible for our approach. Denote by~$\D_\II$
the lcm of~$\II$. Then, for every
element~$\xx$ of~$\MM$, the
$\MM_\II$-tail~$\xx_1$ of~$\xx$ is the right gcd of~$\xx$
and~$\D_\II^{\lg\xx}$, where $\lg\xx$ denotes the common
length of all words representing~$\xx$. Indeed, let~$\xx'_1$
be the latter gcd, and let $\ell=\lg\xx$. By definition,
$\xx_1$ is a right divisor of~$\xx$, so we have
$\lg{\xx_1}\le\ell$, and, as for each~$\zz$ in~$\MM_\II$
satisfying
$\lg\zz\le\ell$ is, we have $\D_\II^\ell \multe \xx'_1$,
hence $\xx'_1 \multe \xx_1$. Conversely, $\xx'_1$
is an element of~$\DR(\xx)\cap\MM_\II$, hence we have
$\xx_1 \multe \xx'_1$, and, finally, $\xx_1=\xx'_1$.
Note that the previous approach does not require
that $\MM$ be of spherical type, but only that
$\MM_\II$ is. Actually, $\MM_\II$ is a closed submonoid
even if it is not of spherical type---but, then, the
characterization of the $\MM_\II$-tail in terms of
powers of~$\D_\II$ vanishes.
\end{exam}

\subsection{Alternating decompositions}

In the second case of Example~\ref{X:Closed}, the
involved subset is a submonoid of~$\MM$, \ie, in
addition to being closed under left lcm, it is closed
under multiplication. From now on, we shall
concentrate on this situation. Then, the decomposition of
Lemma~\ref{L:Tail} takes a specific form.

\begin{defi}
\label{D:Closed}
Assume that $\MM$ is a locally right Garside monoid. We say that a submonoid~$\MM_1$ of~$\MM$ is \emph{closed} if it is closed under both left lcm and left divisor, \ie,
every left lcm of elements of~$\MM_1$ belongs to~$\MM_1$
and every left divisor of an element of~$\MM_1$ belongs
to~$\MM_1$.
\end{defi}

\begin{exam}
If $\MM$ is an Artin--Tits monoid with standard set of
generators~$\SS$, then every standard parabolic submonoid
of~$\MM$ is closed. This need not be the case in every locally
right Garside monoid, or even in every Garside monoid. For
instance, the monoid $\Mon{\aa,\bb \mid
\aa\bb\aa=\bb^2}$ is Garside, hence locally right
Garside---the associated Garside group is the braid
group~$B_3$. However, the submonoid generated by~$\bb$
is not closed, as it contains~$\bb^2$, which is
$\aa\bb\aa$, but it contains neither~$\aa$ nor~$\aa\bb$,
which are left divisors of~$\bb^2$.
\end{exam}

\begin{nota}
For~$\MM$ a monoid, $\xx\in\MM$ and $\AA\ince\MM$,
we write $\xx\perp\AA$ if no non-trivial element of~$\AA$ is a
right divisor of~$\xx$, \ie, if $\DR(\xx)\cap\AA$ is
either~$\emptyset$ or~$\{1\}$.
\end{nota}

\begin{lemm}
\label{L:Decomp}
Assume that $\MM$ is a locally right Garside monoid, and
$\MM_1$ is a closed submonoid of~$\MM$. Then, for each~$\xx$ in~$\MM$,
there exists a unique decomposition
$(\xx', \xx_1)$ of~$\xx$ satisfying
\begin{equation}
\label{E:Decomp}
\xx' \perp\MM_1\text{\quad and\quad}
\xx_1\in\MM_1,
\end{equation}
namely the one given by
$\xx_1=\Tail\xx{\MM_1}$ and $\xx'=\xx\OVER{\xx_1}$.
\end{lemm}

\begin{proof}
Let $\xx_1=\Tail\xx{\MM_1}$ and $\xx'=\xx\OVER{\xx_1}$. We
claim that, for each decomposition $(\yy', \yy_1)$ of~$\xx$ with
$\yy_1\in\MM_1$, we have
\begin{equation}
\label{E:CharTail}
\yy' \perp \MM_1
\quad\Longleftrightarrow\quad
\yy_1=\xx_1.
\end{equation}
First, assume $\zz \in \DR(\xx')\cap\MM_1$.
Then we have $\xx'=\xx''\zz$ for some $\xx''$, hence
$\xx=\xx''\zz\xx_1$, and $\zz\xx_1\in\DR(\xx)$.
As $\zz$ and $\xx_1$ belong to~$\MM_1$ and the
latter is a submonoid of~$\MM$, we deduce
$\zz\xx_1\in\MM_1$, hence $\zz=1$ by definition
of~$\xx_1$. So $\xx' \perp \MM_1$ holds, and the $\Longleftarrow$
implication in~\eqref{E:CharTail} is true.

Conversely, assume $\xx=\yy'\yy_1$ with
$\yy_1\in\MM_1$. By definition of the $\MM_1$-tail, we have
$\xx_1=\zz\yy_1$ for some~$\zz$. The assumption that $\MM_1$ is closed
under left divisor implies $\zz\in\MM_1$. Then we find
$\yy'\yy_1=\xx= \xx'\xx_1 = \xx'\zz\yy_1$, hence
$\yy'=\xx'\zz$ by cancelling~$\yy_1$,
and finally $\zz\in\DR(\yy')\cap\MM_1$. Then
$\DR(\yy')\cap\MM_1=\{1\}$ implies $\zz=1$, \ie,
$\yy_1=\xx_1$, and, from there, $\yy'=\xx'$. So the
$\Longrightarrow$ implication in~\eqref{E:CharTail} is true.
\end{proof}

By definition, the relation $\xx' \perp \MM_1$
of~\eqref{E:CharTail} is equivalent to $\Tail{\xx'}{\MM_1} =
1$. This shows that iterating the decomposition of
Lemma~\ref{L:Decomp} makes no sense: we extracted the
maximal right divisor of~$\xx$ that lies in~$\MM_1$, so,
after that, there remains nothing to extract any longer. But
assume that $\MM$ is locally right Garside, and that
$\MM_2,\MM_1$ are \emph{two} closed submonoids of~$\MM$.
For each~$\xx$ in~$\MM$, Lemma~\ref{L:Decomp} gives
a distinguished decomposition $(\xx', \xx_1)$ of~$\xx$
with~$\xx_1$ in~$\MM_1$. If $\xx'$ is not~$1$, and if
$\MM_2 \cup \MM_1$ generates~$\MM$, the $\MM_2$-tail
of~$\xx'$ is not~$1$, and we obtain a new decomposition
$(\xx'' , \xx_2 , \xx_1)$ of~$\xx$ with $\xx_2\in\MM_2$
and $\xx_1\in\MM_1$. If $\xx''$ is not~$1$, we repeat the
process with~$\MM_1$, \etc\  finally obtaining a
decomposition of~$\xx$ as an alternating sequence of
elements of~$\MM_2$ and~$\MM_1$.

\begin{defi}
\label{D:Covering}
If $\MM$ is a locally right Garside monoid, we say that $(\MM_2,\MM_1)$
is a \emph{covering} of~$\MM$ if $\MM_2$ and $\MM_1$ are
closed submonoids of~$\MM$ and  $\MM_2\cup\MM_1$
generates~$\MM$ (as a monoid).
\end{defi}

\begin{exam}
Let $\MM$ be an Artin--Tits monoid with standard
set of generators~$\SS$, and let $\SS_2,
\SS_1$ be two subsets of~$\SS$ satisfying
$\SS_2\cup\SS_1=\SS$. For $\kk =2,1$, let
$\MM_\kk$ be the standard parabolic submonoid
of~$\MM$ generated by~$\SS_\kk$. Then $(\MM_2, \MM_1)$
is a covering of~$\MM$. Indeed, we already observed that
$\MM_1$ and $\MM_2$ are closed submonoids of~$\MM$.
Moreover, $\SS$ is included in~$\MM_2\cup\MM_1$, so
the latter generates~$\MM$. 
\end{exam}

Similar results hold for every locally right Garside monoid
that is generated by the union of two sets~$\SS_2, \SS_1$
provided we define $\MM_\kk$ to be the smallest \emph{closed}
submonoid of~$\MM$ generated by~$\SS_\kk$. 

\begin{nota}
\label{N:Parity}
For each (nonnegative) integer~$\rr$, we define $\cl\rr$ to
be~$1$ if $\rr$ is odd, and $2$ if $\rr$ is even.
\end{nota}

\begin{prop}
\label{P:AltDec}
Assume that $\MM$ is a locally right Garside monoid and
$(\MM_2,\MM_1)$ is a covering of~$\MM$. Then, for 
every non-trivial element~$\xx$ of~$\MM$, there exists a
unique decomposition $(\xx_\pp, \Ldots, \xx_1)$ of~$\xx$ satisfying
$\xx_\pp\not=1$ and,  for each~$\rr \ge 1$,
\begin{equation}
\label{E:AltDec1}
\xx_p...\xx_{\rr+1} \perp \MM_{\cl\rr}
\text{\quad and\quad}
\xx_\rr\in\MM_{\cl\rr}.
\end{equation}
The elements~$\xx_\rr$ are determined from
$\xx^{(0)}= \xx$ by
\begin{equation}
\label{E:AltDec2}
\xx_\rr= 
\Tail{\xx^{(\rr-1)}}{\MM_{\cl\rr}}
\text{\quad and\quad}
\xx^{(\rr)}=\xx^{(\rr-1)}\OVER{\xx_\rr}. 
\end{equation}
Moreover, we have $\xx_\rr\not=1$ for $\rr\ge2$.
\end{prop}

\begin{proof}
Let $\xx$ belong to~$\MM$, and let $\xx_\rr$,
$\xx^{(\rr)}$ be as specified by~\eqref{E:AltDec2}.
Using induction on $\rr\ge1$, we first prove the
relations
\begin{gather}
\label{E:Alt3}
\xx=\xx^{(\rr)}\xx_\rr\cdots\xx_1,\\
\label{E:Alt4}
\xx^{(\rr)} \perp\MM_{\cl\rr}.
\end{gather}
For $\rr=1$, Lemma~\ref{L:Decomp} for~$\xx$
and~$\MM_1$ gives $\xx=\xx^{(1)}\xx_1$, which
is~\eqref{E:Alt3}, and
$\xx^{(1)} \perp \MM_1$, which is~\eqref{E:Alt4}.
Assume $\rr\ge2$. Then \eqref{E:AltDec2} implies
$\xx^{(\rr-1)}=\nobreak\xx^{(\rr)}\xx_\rr$, and, 
susbtituting in $\xx=\xx^{(\rr-1)}\xx_{\rr-1}...\xx_1$, which
holds by induction hypothesis, we obtain~\eqref{E:Alt3}.
Moreover, Lemma~\ref{L:Decomp} for~$\xx^{(\rr)}$ and
$\MM_{\cl\rr}$ gives~\eqref{E:Alt4}.

By construction, the sequence $\xx_1$, $\xx_2\xx_1$,
$\xx_3\xx_2\xx_1$, ... is increasing
in~$(\DR(\xx), \div)$. By Condition~$(C_3)$, it is
eventually constant. By right cancellability, this implies
that there exists~$\pp$ such that $\xx_\rr=\xx^{(\rr)}=1$
holds for all~$\rr\ge\pp$. Then \eqref{E:Alt3} implies
$\xx=\xx_\pp...\xx_1$, with $\xx_\pp\not=1$ provided
$\pp$ is chosen to be minimal and $\xx$ is not~$1$.

So the expected sequence $(\xx_\pp, \Ldots, \xx_1)$ exists and
satisfies~\eqref{E:AltDec1} and~\eqref{E:AltDec2}. We show now
$\xx_\rr\not=1$ for $\rr\ge2$. Indeed, assume
$\xx^{(\rr-1)}\not=1$. By hypothesis,
$\MM_2\cup\MM_1$ generates~$\MM$, implying
$\xx^{(\rr-1)} \not\perp (\MM_2\cup\MM_1)$.
By~\eqref{E:Alt4}, we have $\xx^{(\rr-1)} \perp
\MM_{\cl{\rr-1}}$, hence $\xx^{(\rr-1)} \not\perp
\MM_{\cl\rr}$. Therefore the $\MM_{\cl\rr}$-tail
of~$\xx^{(\rr-1)}$, which by definition is~$\xx_\rr$, is
not~$1$---the argument fails for $\rr=1$ because $\xx^{(0)}
\perp \MM_{\cl0}$ need not hold.

We turn to uniqueness. Consider any decomposition
$(\yy_\qq,  \Ldots, \yy_1)$ of~$\xx$ satisfying $\yy_\qq\not=1$ with 
$\yy_\rr\in\MM_{\cl\rr}$ and
$\yy_\qq...\yy_{\rr+1} \perp\MM_{\cl\rr}$ for
each~$\rr$. We inductively prove $\yy_\rr=\xx_\rr$
and $\yy_\qq...\yy_{\rr+1}=\xx^{(\rr)}$ for
$\rr\ge1$. For $\rr=1$, the hypotheses $\xx=
(\yy_\qq...\yy_2)\yy_1$ with $\yy_1\in\MM_1$ and
$\yy_\qq...\yy_2 \perp \MM_1$ imply
$\yy_1=\xx_1$ and $\yy_\qq...\yy_2= \xx^{(1)}$ by
Lemma~\ref{L:Decomp}. Assume $\rr\ge2$. By induction
hypothesis, we have
$\yy_\qq...\yy_{\rr}= \xx^{(\rr-1)}$, and the hypotheses
about the elements~$\yy_\jj$ give
$\xx^{(\rr-1)} = (\yy_\qq...\yy_{\rr+1}) \yy_\rr$ with
$\yy_\rr\in\MM_{\cl\rr}$ and
$\yy_\qq...\yy_{\rr+1} \perp \MM_{\cl\rr}$. Then
Lemma~\ref{L:Decomp} implies
$\yy_\rr=\Tail{\xx^{(\rr-1)}}{\MM_{\cl\rr}} = \xx_\rr$ and
$\yy_\qq...\yy_{\rr+1}= \xx^{(\rr-1)}\OVER{\xx_\rr} =
\xx^{(\rr)}$. Finally,  $\qq>\pp$ would imply
$\xx_\qq=\yy_\qq\not=1$, contradicting the choice of~$\pp$.
\end{proof}

\begin{defi}
\label{D:Decomp}
In the framework of Proposition~\ref{P:AltDec}, the sequence
$(\xx_\pp, \Ldots, \xx_1)$ is called the  \emph{$(\MM_2,\MM_1)$-decomposition} of~$\xx$.
\end{defi}

\begin{figure}[htb]
\begin{picture}(72,23)
\put(0,0){\includegraphics{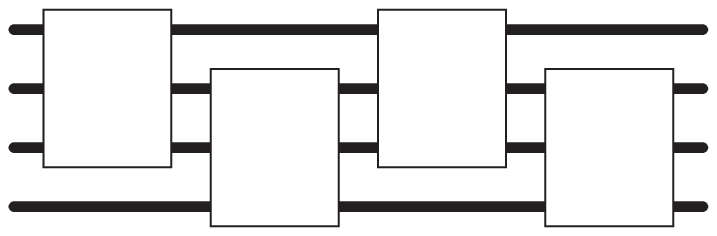}}
\put(60.8,10){$\xx_1$}
\put(44,16){$\xx_2$}
\put(26.8,10){$\xx_3$}
\put(10,16){$\xx_4$}
\put(-4,13){$\dots$}
\end{picture}
\caption{\sf Diagram associated with the $(\MM_2,
\MM_1)$-decomposition of a $4$-braid: starting from the right,
we alternatively select the maximal right divisor that does not
involve the $\nn$th strand and the first strand.}
\label{F:Alt}
\end{figure}

\begin{exam}
\label{X:Delta}
Consider the $4$-strand braid monoid~$\BB4$. Let
$\MM_1$ be the submonoid generated by~$\ss1$
and~$\ss2$, \ie, $\BB3$, and $\MM_2$ be the submonoid
generated by~$\ss2$ and~$\ss3$. Choose
$\xx=\D_4^2=(\ss1\ss2\ss1\ss3\ss2\ss1)^2$. The
computation of the $(\MM_2,\MM_1)$-decomposition
of~$\xx$ is as follows:
\begin{align*}
\smash{\xx^{(0)}} 
&= \xx = \smash{\D_4^2}
&\xx_1
&=\smash{\Tail{\xx^{(0)}}{\MM_1}
=\D_3^2},\\
\smash{\xx^{(1)}}
&=\smash{\xx^{(0)}}\OVER{\xx_1}
=\ss3\ss2\s_1^2\ss2\ss3
&\xx_2 
&=\smash{\Tail{\xx^{(1)}}{\MM_2}}
=\ss2\ss3,\\
\xx^{(2)}
&=\xx^{(1)}\OVER{\xx_2}
=\smash{\ss3\ss2\s_1^2}
&\xx_3 
&=\Tail{\xx^{(2)}}{\MM_1}
=\smash{\ss2\s_1^2},\\
\xx^{(3)}
&=\xx^{(2)}\OVER{\xx_3}
=\ss3
&\xx_4 
&=\Tail{\xx^{(3)}}{\MM_2}
=\ss3,\\
\xx^{(4)}
&=\smash{\xx^{(3)}\OVER{\xx_4}}
=1.
\end{align*}
Thus the $(\MM_2,\MM_1)$-decomposition of~$\D_4^2$ is the
sequence $(\ss3, \ss2\s_1^2, \ss2\ss3, \D_3^2)$---see Figure~\ref{F:Alt} for an illustration in terms of
standard braid diagrams. Note that the decomposition depends on the
order of the submonoids: the
$(\MM_1,\MM_2)$-decomposition of~$\D_4^2$ is
$(\ss1, \ss2\s_3^2, \ss2\ss1, (\ss2\ss3\ss2)^2)$.
\end{exam}

\begin{rema}
In the framework of Proposition~\ref{P:AltDec}, $\xx_1$ is the
left lcm of all right divisors of~$\xx$ that lie in~$\MM_1$. Comparing with the case of the greedy normal form, we
might expect that, similarly, $\xx_2 \xx_1$ is the left lcm of
all right divisors of~$\xx$ of the form~$\yy_2 \yy_1$ with
$\yy_\kk \in \MM_\kk$, \ie, lying in $\MM_2\MM_1$. This is
not the case. Consider Example~\ref{X:Delta} again, and let
$\xx = \s_1^\ee\ss2\ss1$ with~$\ee\ge1$. Then the
$(\MM_2, \MM_1)$-decomposition of $\xx$ is $(\s_1^\ee,
\ss2, \ss1)$, so $\xx_2 \xx_1$ is $\ss2\ss1$ here. Now, we also
have $\xx = \ss2 \ss1 \s_2^\ee$, so $\s_2^\ee$, \ie,
$\s_2^\ee \cdot 1$, is a right divisor of~$\xx$ that belongs
to~$\MM_2\MM_1$ and does not divide $\ss2 \ss1$.
More generally, we see that the braids~$\ss\ii$ that are
right divisors of~$\xx$ cannot be retrieved from the last two
elements of the $(\MM_2, \MM_1)$-decomposition of~$\xx$.
\end{rema}

\begin{rema}
\label{R:Shortest}
Assume that $\MM$ is a locally right Garside monoid, and $(\MM_2, \MM_1)$ is a covering of~$\MM$.
Define an \emph{$(\MM_2, \MM_1)$-sequence} to be any finite
sequence $(\xx_\pp, \Ldots, \xx_1)$ such that $\xx_\rr$ belongs
to~$\MM_{\cl\rr}$ for each~$\rr$.
Then the $(\MM_2, \MM_1)$-decomposition of~$\xx$ is a
certain decomposition of~$\xx$ that is a $(\MM_2,
\MM_1)$-sequence. As we take the maximal right divisor at
each step, we might expect to obtain a short
$(\MM_2, \MM_1)$-sequence, possibly the shortest possible one.
We shall see in Section~\ref{S:Burckel} below that this is
indeed the case for the covering of Example~\ref{X:Delta}.
However, this is not the case in general.
Indeed, keep the braid monoid~$\BB4$, but consider
the covering $(\MM'_2, \MM'_1)$, where $\MM'_1$
(\resp $\MM'_2$) is the submonoid generated by~$\ss1$
and~$\ss3$ (\resp by~$\ss2$ and~$\ss3$). Let $\xx$ be $\ss3\ss2\ss1\ss2\s_3^2\ss2$. The
$(\MM'_2,\MM'_1)$-decomposition of~$\xx$ turns out to be
$(\ss1 , \s_2^2 , \ss3\ss1 , \ss2 , \s_1)$, a sequence of
length~$5$, but another decomposition of~$\xx$ is the
$(\MM'_2, \MM'_1)$-sequence 
$(\ss3\ss2, \ss1, \ss2\s_3^2\ss2, 1)$, which has length~$4$: choosing the
maximal right divisor at each step does \emph{not} guarantee that we obtain the
shortest sequence.
\end{rema}

Finally, it should be clear that, instead of
considering two closed submonoids $\MM_2,\MM_1$
of~$\MM$, we could consider any finite family of such
submonoids $\MM_{\mm}, \Ldots, \MM_1$. Provided the union
of all~$\MM_\jj$'s generates~$\MM$, we can extend
Proposition~\ref{P:AltDec} and obtain for every
element~$\xx$ of~$\MM$ a distinguished decomposition
$(\xx_\pp, \Ldots, \xx_1)$ such that $\xx_\rr$ belongs to
$\MM_{\cl\rr}$ and
$\xx_p...\xx_{\rr+1}\perp\MM_{\cl\rr}$ holds for
every~$\rr$, where $\cl\rr$ now denotes the unique
element of~$\{1, \Ldots, \mm\}$ that equals~$\rr$
mod~$\mm$. The only difference is that
the condition $\xx_\rr\not=1$ for $\rr\ge2$ has to be relaxed to
$\xx_{\rr+\mm-2}...\xx_\rr\not=1$ for $\rr \ge \mm$, since
the conjunction of $\xx\not=1$ and
$\xx\perp\MM_{\cl\rr}$ need not guarantee 
$\xx\not\perp\MM_{\cl{\rr+1}}$, but only
$\xx\not\perp
(\MM_{\cl{\rr+\mm-1}}\cup...\cup\MM_{\cl{\rr+1}})$. 
Adapting is easy---see~\cite{Fro} for an example.

\subsection{Algorithmic aspects}

Computing the alternating decomposition is easy provided one can efficiently
perform right division in the ground monoid. To give a precise statement, we
recall from~\cite{Dfx} the notion of word norm (or pseudolength)
that generalizes the standard notion of word length. In the
sequel, for $\SS$ included in~$\MM$ and~$\ww$ a word on~$S$, we denote by~$\CL\ww$ the element of~$\MM$
represented by~$\ww$.

\begin{defi}
\label{D:Norm}
Assume that $\MM$ is a locally right Garside monoid that satisfies Condition~$(C_3^+)$, and $\SS$
generates~$\MM$. For~$\ww$ a word on~$\SS$, we denote
by~$\norm\ww$ the maximal length of a word~$\ww'$
satisfying~$\CL{\ww'}=\CL\ww$. 
\end{defi}

Condition~$(C_3^+)$ is precisely what is needed to guarantee that $\norm\ww$ exists for every word~$\ww$.
Indeed, if $\lambda: \MM \to \Nat$ witnesses that $(C_3^+)$ is satisfied, then every word~$\ww'$
satisfying~$\CL{\ww'}=\CL\ww$ must satisfy $\lg\ww \le \l(\CL\ww)$. Conversely, if $\norm\ww$ exists for each
word~$\ww$, then the map $\ww \mapsto \norm\ww$ induces a well-defined map of~$\MM$ to~$\Nat$ that
witnesses~$(C_3^+)$. In the case of Artin--Tits monoids and, more generally, of monoids presented by homogeneous
relations,  $\norm\ww$ coincides with the length~$\lg\ww$.

\begin{prop}
\label{P:Complexity1}
Assume that $\MM$ is a locally right Garside monoid,
generated by some finite set~$\SS$, and satisfying
Condition~$(C_3^+)$ plus:
\begin{quote}
$(*)$ There exists an algorithm~$\AlgoD$ that, for~$\ww$ a word on~$\SS$ and $\xs$ in~$\SS$, runs in
time~$O(\norm\ww)$, recognizes whether
$\CL\ww\multe\xs$ holds and, if so, returns a word
representing~$\CL\ww\OVER{\xs}$. 
\end{quote} 
Let $\SS_2,\SS_1\ince\SS$ satisfying $\SS_2\cup\SS_1=\SS$. Let
$\MM_\kk$ be the submonoid of~$\MM$ generated
by~$\SS_\kk$, and suppose that $\MM_2, \MM_1$ are closed. 
Then there exists a algorithm that, for~$\ww$
a word on~$\SS$, runs in time~$O(\norm\ww^2)$ and
computes the $(\MM_2,\MM_1)$-decomposition
of~$\CL\ww$.
\end{prop}

\begin{proof}
Having listed the elements of~$\SS_1$ and~$\SS_2$, and
starting with~$\ww$, we use~$\AlgoD$ to divide by
elements of~$\SS_1$ until division fails, then we divide by
elements of~$\SS_2$ until division fails, \etc\ We stop when
the remainder is~$1$. If we start with a word~$\ww$
satisfying $\norm\ww = \ell$, then the words~$\ww_\rr$
subsequently occurring represent the elements~$\xx^{(\rr)}$
of~\eqref{E:AltDec2}, which are left divisors of~$\xx$, and,
hence, we have $\lg{\ww_\rr} \le \norm{\ww_\rr} \le
\ell$. Moreover, at each step, $\norm{\ww_\rr}$ decreases by at least~$1$,
so termination occurs after at most
$\card(\SS) \times \ell$~division steps. By hypothesis, the
cost of each division step is bounded above by~$O(\ell)$,
whence a quadratic global upper bound.
\end{proof}

\begin{exam}
Let $\MM$ be an Artin--Tits of spherical type, or, more
generally, a Garside monoid, and let $\SS$ be the set of
atoms in~$\MM$. Then there exist division
algorithms running in linear time, \emph{e.g.}, those
involving a rational transducer based on the (right)
automatic structure~\cite{Eps}. Alternatively, for the specific
question of dividing by an atom, the reversing method
of~\cite{Dgp} is specially convenient.
\end{exam}

\section{Iterated alternating decompositions}
\label{S:Iter}

If the submonoids involved in a covering are monogenerated, it makes no sense to iterated the alternating
decomposition. But, in general, for instance in the case of Example~\ref{X:Delta}, the covering submonoids
need not be monogenerated, and they can in turn be covered by smaller submonoids. In such cases, it is natural to
iterate the alternating decomposition using a sequence of nested coverings. This is the idea we develop in this
section. The main observation is that the result of the iterated decomposition can be obtained directly, without any
iteration.

\subsection{Iterated coverings}
\label{S:IterCov}

The possibility of iterating the alternating decomposition
relies on the following trivial observation:

\begin{lemm}
\label{L:Sub}
Every closed submonoid of a locally right Garside monoid is
locally right Garside.
\end{lemm}

\begin{proof}
Assume that $\MM_1$ is a closed submonoid of
a locally right Garside monoid~$\MM$. First,
$\MM_1$ admits right cancellation as every submonoid of a
right cancellative monoid does. Then, if $\xx,
\yy$ belong to~$\MM_1$ and admit a common left
multiple~$\zz$ in~$\MM_1$, then $\zz$ is a common left
multiple of~$\xx$ and~$\yy$ in~$\MM$, so, in~$\MM$,
the left lcm~$\zz'$ of~$\xx$ and~$\yy$ exists. The
hypothesis that $\MM_1$ is closed under left lcm implies
$\zz' \in \MM_1$, and, then,
$\zz'$ must be a left lcm for~$\xx$ and~$\yy$ in the sense
of~$\MM_1$. Finally, the right divisibility
relation of~$\MM_1$ is included in the right divisibility
relation of~$\MM$, so a sequence contradicting
Condition~$(C_3)$ in~$\MM_1$ would also
contradict~$(C_3)$ in~$\MM$. 
\end{proof}

Assume that $\MM$ is a locally right Garside monoid and $(\MM_2,\MM_1)$ is a covering of~$\MM$.
By Lemma~\ref{L:Sub}, $\MM_2$ and $\MM_1$ are locally
right Garside, and we can repeat the process: assuming
that $(\MM_{\kk,2},\MM_{\kk,1})$ is a covering of~$\MM_\kk$
for~$\kk=2,1$, every element of~$\MM_\kk$ admits a
$(\MM_{\kk,2},\MM_{\kk,1})$-decomposition, and,
therefore, every element of~$\MM$ admits a distinguished
decomposition in terms of the four
monoids~$\MM_{22}$, $\MM_{21}$,
$\MM_{12}$, and $\MM_{22}$---we drop commas in
indices.

\begin{exam}
\label{X:IterCov}
As in Example~\ref{X:Delta}, consider the $4$-strand
braid monoid~$\BB4$, and let $\MM_2, \MM_1$ be the
parabolic submonoids respectively generated
by~$\ss3, \ss2$, and by $\ss2, \ss1$. Then let $\MM_{22}$,
$\MM_{21}$, $\MM_{12}$, and $\MM_{11}$ be the
submonoids respectively generated by~$\ss2$, $\ss3$,
$\ss2$, and~$\ss1$. Then  $(\MM_{\kk2},\MM_{\kk1})$ is
a covering of~$\MM_\kk$ for~$\kk=2,1$.
\end{exam}

To make the construction formal, we introduce the notion of an iterated covering.

\begin{defi}
Assume that $\MM$ is a locally right Garside monoid. We say that $\MM$ is a \emph{$0$-covering} of itself, and
then, for $\nn \ge 1$, we define an \emph{$\nn$-covering} of~$\MM$ to be a pair $(\MMM_2, \MMM_1)$
such that there exists a covering $(\MM_2, \MM_1)$ of~$\MM$ such that $\MMM_\kk$ is an $(\nno)$-covering
of~$\MM_\kk$ for $\kk = 1,2$.
\end{defi}

So a $1$-covering of~$\MM$ is just an ordinary covering, and, for instance, a $2$-covering of~$\MM$ consists of a
covering $(\MM_2, \MM_1)$ of~$\MM$, plus coverings of~$\MM_2$ and~$\MM_1$, as in Example~\ref{X:IterCov}.

An iterated covering of a monoid~$\MM$ has the structure of a binary tree, and we can specify the various
submonoids by using finite sequences of twos and ones---or of ones and zeroes, or of letters `L' and `R'---to
indicate at each forking which direction is to be taken. In the sequel, such a finite sequence of length~$\nn$ is called
a \emph{binary $\nn$-address}. In this way, an $\nn$-covering of a monoid~$\MM$ is a sequence of submonoids
$\MM_\a$ indexed by binary addresses of length at most~$\nn$, such that, for each~$\a$ of length smaller
than~$\nn$, the pair~$(\MM_{\a2}, \MM_{\a1})$ is a covering of~$\MM_\a$, and $\MM_\ea$ is~$\MM$---using
$\ea$ for the empty address. In the sequel, if $\MMM$ is an iterated covering, we shall always use~$\MM_\a$ for
the $\a$-entry in~$\MMM$.

\begin{figure}[htb]
\begin{picture}(30,22)
\put(0,-1.5){\includegraphics{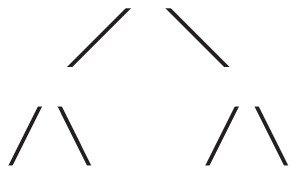}}
\put(0,0){$\ss2$}
\put(10,0){$\ss3$}
\put(20,0){$\ss2$}
\put(30,0){$\ss1$}
\put(2,10){$\ss2,\!\ss3$}
\put(22,10){$\ss1,\!\ss2$}
\put(10,20){$\ss1,\!\ss2,\!\ss3$}
\end{picture}
\caption{\sf Skeleton of the $2$-covering of~$\BB4$  of Example~\ref{X:Iter}: a depth~$2$ binary
tree displaying the inclusions between the generating sets of the successive submonoids; this example corresponds
to $\SS_{22}=\SS_{12}=\{\ss2\}$, $\SS_{21}=\{\ss3\}$, and $\SS_{11}=\{\ss1\}$;  we find for instance $\MM_2 =
\Mon{\ss2,\ss3}$, and $\MM_{12} = \Mon{\ss2}$.}
\label{F:Skeleton}
\end{figure}

If the ground monoid~$\MM$ has some distinguished generating set~$\SS$, we can specify an
$\nn$-covering by choosing a subset~$\SS_\a$ of~$\SS$ for
each~$\a$ in~$\{2,1\}^\nn$, and, for~$\b$ in~$\{2,1\}^\nn$ with $\mm\le\nn$, defining~$\MM_\b$ to be
the submonoid generated by all~$\SS_{\a}$'s such that $\b$ is a
prefix of~$\a$. We obtain an $\nn$-covering provided each
submonoid~$\MM_\b$ is closed. For such coverings, we can 
display the inclusions in a binary tree---see Figure~\ref{F:Skeleton}.

\subsection{Iterated $\MMM$-decomposition}
\label{S:IterDec}

As was shown in Section~\ref{S:Alt}, each covering $(\MM_2, \MM_1)$ of a monoid~$\MM$ leads to a
distinguished decomposition for the elements of~$\MM$ in terms of elements of~$\MM_2$ and~$\MM_1$. An
iterated covering similarly leads to what can be called an iterated decomposition.

\begin{defi}
\label{D:IterDec}
Assume that $\MM$ is a locally right Garside monoid, and $\MMM$ is an $\nn$-covering of~$\MM$. For
$\xx$ in~$\MM$, we define the \emph{$\MMM$-decomposition}~$\DD\MMM\xx$ of~$\xx$ by $\DD\MMM\xx=\xx$ for
$\nn=0$, and, for $\nn \ge 1$ and $\MMM = (\MMM_2, \MMM_1)$, by
\begin{equation}
\label{E:IterDec}
\DD\MMM\xx=
(\DD{\MMM_{\cl\pp}}{\xx_\pp} \,,\, ... \,,\,
\DD{\MMM_1}{\xx_1}),
\end{equation}
where $(\xx_\pp, \Ldots, \xx_1)$ is the $(\MM_2,
\MM_1)$-decomposition of~$\xx$. 
\end{defi}

\begin{exam}
\label{X:Iter}
Consider the braid~$\D_4^2$ of $\BB4$ and the covering~$\MMM$ of Example~\ref{X:IterCov}. We saw in
Example~\ref{X:Delta} that the $(\MM_2, \MM_1)$-decomposition of~$\D_4^2$ is
$$(\ss3, \ss2\s_1^2,  \ss2\ss3,  \D_3^2).$$
Now, the $(\MM_{12}, \MM_{11})$-decomposition of~$\D_3^2$ turns out to be
$(\ss2, \s_1^2,\ss2,\s_1^2)$. Similarly, the $(\MM_{22}, \MM_{21})$-decomposition of~$\ss2\ss3$
is $(\ss2, \ss3)$. Continuing in this way, we obtain 
\begin{equation}
\label{E:XIterDec1}
\DD\MMM{\D_4^2} = ((\ss3), (\ss2,\s_1^2),(\ss2,\ss3),
(\ss2, \s_1^2,\ss2,\s_1^2)), 
\end{equation}
corresponding to the factorization
$\D_4^2 = (\ss3)\,\cdot\, (\ss2 \cdot \s_1^2)\,\cdot\, (\ss2
\cdot \ss3)\,\cdot\, (\ss2\cdot\s_1^2\cdot \ss2 \cdot
\s_1^2)$.
\end{exam}

\begin{figure}[htb]
\begin{picture}(95,27)
\put(-1,3.3){\includegraphics{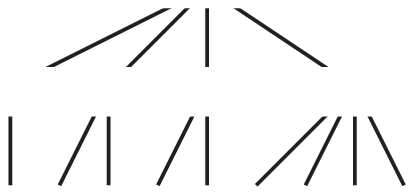}}
\put(49.2,3){\includegraphics{Tree.eps}}
\put(0,0){$\ss3$}
\put(5,0){$\ss2$}
\put(10,0){$\s_1^2$}
\put(15,0){$\ss2$}
\put(20,0){$\ss3$}
\put(25,0){$\ss2$}
\put(30,0){$\s_1^2$}
\put(35,0){$\ss2$}
\put(40,0){$\s_1^2$}
\put(0,12){$\ss3$}
\put(7,12){$\ss2\s_1^2$}
\put(17,12){$\ss2\ss3$}
\put(33,12){$\D_3^2$}
\put(19,23){$\D_4^2$}
\put(49,0){$41$}
\put(54,0){$32$}
\put(59,0){$31$}
\put(64.5,0){$22$}
\put(69.5,0){$21$}
\put(74.5,0){$14$}
\put(79.5,0){$13$}
\put(84.5,0){$12$}
\put(90,0){$11$}
\put(50,12){$4$}
\put(60,12){$3$}
\put(70,12){$2$}
\put(85,12){$1$}
\put(70,23){$\ea$}
\end{picture}
\caption{\sf The tree associated with the
$2$-sequence of~\eqref{E:XIterDec1}: on the
left,  the braid entries, on the right, the addresses; the entry list specifies
the name of the leaves, while the address
list specifies the shape of the tree; for each address~$\h$, the $\h$-subsequence~$\www_\h$ corresponds to what
lies below~$\h$ in~$\www$; here, the
$31$-subsequence is~$\s_1^2$, while the $2$-subsequence
is~$(\ss2, \ss3)$. The $23$-subsequence does not exist.}
\label{F:Tree}
\end{figure}

For $\nn \ge 2$, the $\MMM$-decomposition of an element is a sequence of sequences. More precisely, it is an
\emph{$\nn$-sequence}, defined to be a single element for $\nn = 0$, and to be a sequence of $(\nno)$-sequences
for $\nn \ge 1$. Such iterated sequences can naturally be viewed as trees, on the model of Figure~\ref{F:Tree} (left).

Entries in an ordinary sequence of length~$\pp$ are usually specified using numbers from~$1$
to~$\pp$---or rather $\pp$ to~$1$ in the context of this paper where we start from the right. Entries in an iterated
sequence are then specified using finite sequence of numbers, as done in Section~\ref{S:IterCov} with binary
addresses. In the sequel, a length~$\nn$ sequence of positive numbers is called an \emph{$\nn$-address}: for
instance, $32$ is a typical
$2$-address---in examples, we drop brackets and separating commas. If $\www$ is an $\nn$-sequence, and $\h$ is
an $\mm$-address with $\mm \le \nn$, we denote by~$\www_\h$ the $\h$-subsequence of~$\www$, \ie, the
$(\nn-\mm)$-sequence made by those entries in~$\www$ whose address begins with~$\h$---when it exists, \ie,
when the considered sequences are long enough---see Figure~\ref{F:Tree} (right). 

Note that addresses are just a way of specifying brackets in an iterated sequence: an
$\nn$-sequence is determined by its unbracketing---that is, the (ordinary) sequence obtained by removing all
inner brackets---and its address list. For instance, in the $2$-sequence of~\eqref{E:XIterDec1}, the unbracketing and
the address list are
\begin{equation}
\label{E:XIterDec3}
(\ss3, \ss2,\s_1^2, \ss2,\ss3, \ss2, \s_1^2,\ss2,\s_1^2)
\text{\quad and \quad}
(41,32,31,22,21,14,13, 12,11).
\end{equation}

Assume that $\www$ is the $\MMM$-decomposition of an element~$\xx$. For each~$\h$ that is the address of a
node of~$\www$ (viewed as a tree), write $\xx_\h$ for the product of the subsequence~$\www_\h$. Then, by
definition, if
$\h$ is the address of an inner node and $\h\pp, \Ldots,
\h1$ are the addresses of the nodes that lie immediately below~$\h$ in~$\www$, the sequence
$(\xx_{\h\pp}, \Ldots, \xx_{\h1})$ is the $(\MM_{\cl\h2}, \MM_{\cl\h1})$-decomposition of~$\xx_\h$, where 
$\cl\h$ denotes the binary address obtained by replacing each~$\rr$ ocurring in~$\h$ with~$\cl{\rr}$---which is
coherent with Notation~\ref{N:Parity}. Applying Proposition~\ref{P:AltDec} immediately gives the following
characterization.

\begin{prop}
\label{P:IterDec}
Assume that $\MM$ is a locally right Garside monoid, $\MMM$ is an $\nn$-covering of~$\MM$, and $\www =
\DD\MMM\xx$. For each address~$\h$ in~$\www$, let
$\xx_\h$ denote the product of~$\www_\h$. Assume that $\h$ is the address of an inner node and
$\h\pp, \Ldots, \h1$ are the addresses of the nodes that lie immediately below~$\h$ in~$\www$. Then, the
elements~$\xx_{\h\rr}$ are determined from
$\xx_\h^{(0)}= \xx_\h$ by
\begin{equation}
\label{E:IterDec1}
\xx_{\h\rr}=\Tail{\xx_\h^{(\rr-1)}}{\MM_{\cl{\h\rr}}}
\text{\quad and\quad}
\xx_\h^{(\rr)} = \xx_\h^{(\rr-1)}\OVER{\xx_{\h\rr}}.
\end{equation}
\end{prop}

\begin{exam}
\label{X:Complicated}
In the context of Example~\ref{X:Iter} and Figure~\ref{F:Tree}, \eqref{E:IterDec1} gives
$$\xx_1=
\ss2\s_1^2\ss2\s_1^2 = 
\Tail\xx{\BB3}, \quad
\xx_2= \ss2\ss3 = 
\Tail{\ss3\ss2\s_1^2\ss2\ss3}{\Mon{\ss2,\ss3}}, \etc$$
which involve the whole of~$\xx$, but also, at the next level,
we have
$$\xx_{11}= \s_1^2 = 
\Tail{\ss2\s_1^2\ss2\s_1^2}{\BB2}, \quad
\xx_{12}= \ss2 = 
\Tail{\ss2\s_1^2\ss2}{\Mon{\ss2}}, \etc
$$
which only involve the element
$\xx_1$, namely~$\ss2\s_1^2\ss2\s_1^2$, and not the
whole of~$\xx$.
\end{exam}

\subsection{A transitivity lemma}
\label{S:Iter2}

Proposition~\ref{P:IterDec} looks intricate, and it is not
satisfactory in that it does not give a global characterization
of the $\MMM$-decomposition and a way to obtain it directly. This is
what we shall do now. The point is that, according to the following result, there is no
need to consider local remainders when computing
iterated tails.

\begin{lemm}
\label{L:BiTail}
Assume that $\MM$ is a locally right Garside monoid, that $\MM_1$ is a closed submonoid of~$\MM$, and
that $\MM_{11}$ is a closed submonoid of~$\MM_1$. Then,
for every~$\zz$ in~$\MM$ and every left divisor~$\yy$
of~$\Tail\zz{\MM_1}$, we have
\begin{equation}
\label{E:Bitail}
\Tail{(\zz\OVER{\Tail\zz{\MM_1}})\yy}{\MM_{11}}.
= \Tail{\yy}{\MM_{11}}.
\end{equation}
\end{lemm}

\begin{proof}
Put $\zz_1= \Tail\zz{\MM_1}$ and $\zz'=\zz\OVER{\zz_1}$. By
definition,
$\Tail{\yy}{\MM_{11}}$ is a right divisor of
$\Tail{\zz'\yy}{\MM_{11}}$, hence the point is to prove
that every right divisor of~$\zz'\yy$ lying
in~$\MM_{11}$ is a right divisor of~$\yy$. So assume
$\zz'\yy=\xx'\xx$ with
$\xx\in\MM_{11}$. By hypothesis, we have
$\zz_1=\yy\zz'_1$ for some~$\zz'_1$, necessarily lying
in~$\MM_1$. Then, we have $\zz= \zz'\zz_1 =
\zz'\yy\zz'_1 = \xx'\xx\zz'_1$. Now $\xx\in\MM_{11}$
implies $\xx\in\MM_1$, hence $\xx\zz'_1\in\MM_1$, and
$\xx\zz'_1$ has to be a right divisor of~$\Tail\zz{\MM_1}$,
\ie, of~$\zz_1$, which is also~$\yy\zz'_1$. It follows that
$\xx$ is a right divisor of~$\yy$, as expected.
\end{proof}

In particular, when we choose~$\yy$ to be~$\zz_1$ itself,
\eqref{E:Bitail} gives
\begin{equation}
\label{E:BiTail2}
\Tail\zz{\MM_{11}}=\Tail{\Tail\zz{\MM_1}}{\MM_{11}},
\end{equation}
which is vaguely reminiscent of the equality
$\Tail{\zz\yy}\Sigma = 
\Tail{\Tail\zz\Sigma\yy}\Sigma$ that is crucial
in the construction of the right greedy normal form in a Garside
monoid.

\subsection{Global characterization of the iterated decomposition}
\label{S:Global}

We shall now give a direct description of the
$\MMM$-decomposition not involving the
intermediate values~$\xx_\h$. Consider Examples~\ref{X:Iter} and~\ref{X:Complicated} again. The
problem is as follows: in the case of the $1$-covering
of~$\BB3$, only two submonoids are involved, and the final
decomposition consists of alternating blocks belonging to
each of them; in the case of the $2$-covering
of~$\BB4$, the decomposition consists of blocks
of~$\ss1$'s, $\ss2$'s, and~$\ss3$'s, but the order in which
these blocks appear is not so simple. Indeed, on the left of a
block of~$\ss2$'s, there may be either a block of~$\ss1$'s or
a block of~$\ss3$'s, depending on the
current address, \ie, on the position in
(the skeleton of) the covering, typically on which of the
two occurrences of~$\ss2$ in the tree of
Figure~\ref{F:Skeleton} the considered block of~$\ss2$'s is
to be associated: on the left of a block of~$\ss2$'s associated
with the rightmost~$\ss2$ in Figure~\ref{F:Skeleton}, $\ss1$ is expected, while $\ss3$
is expected in the other case. This is what
Proposition~\ref{P:GlobalIter} below says, namely that the
$\MMM$-decomposition can be obtained directly provided
we keep track of some position specified by a binary address. 

To make the description precise, we introduce the notion of
successors of an address. It comes in two versions, one for
general addresses, one for binary addresses. 

\begin{defi}
\label{D:Successor}
For $\h$ an $\nn$-address and $0\le\mm\le\nn$, the
\emph{$\mm$-successor}~$\Succ\mm\h$ of~$\h$ is 
the $\nn$-address obtained by keeping the first
$\mm$ digits of~$\h$, adding~$1$ to the next one, and
completing with~$1$'s, \ie, for $\h=\dd_1...\dd_\nn$,
the $\mm$-successor is $\dd'_1...\dd'_\nn$ with
$\dd'_\rr=\dd_\rr$ for $\rr\le\mm$, and, if $\mm<\nn$
holds,  $\dd'_{\mm+1}=\dd_{\mm+1}+1$ and $\dd'_\rr=1$
for $\rr>\mm+1$. For $\a$ a binary $\nn$-address, the
\emph{binary $\mm$-successor}~$\BSucc\mm\a$ of~$\a$ is
defined to be~$\cl{\Succ\mm\a}$. 
\end{defi}

\begin{exam}
Let $\h=3612$. The successors of~$\h$ are
$$\Succ0\h=4111, \quad
\Succ1\h=3711, \quad
\Succ2\h=3621, \quad
\Succ3\h=3613, \quad
\Succ4\h=3612.$$
Similarly, the binary successors of~$\a=1212$ are
$$\BSucc0\a=2111, \quad
\BSucc1\a=1111, \quad
\BSucc2\a=1221, \quad
\BSucc3\a=1211, \quad
\BSucc4\a=1212.$$
\end{exam}

Note that $\Succ\nn\h=\h$ holds for every $\nn$-address~$\h$. We recall that specifying an iterated sequence
amounts to specifying both its unbracketing and its address list.

\begin{prop}
\label{P:GlobalIter}
Assume that $\MM$ is a locally right Garside monoid, and $\MMM$ is an $\nn$-covering of~$\MM$.
Then, for~$\xx$ in~$\MM$, the unbracketing
$(\xx_\pp, \Ldots, \xx_1)$ and the address list $(\h_\pp, \Ldots,
\h_1)$ of~$\DD\MMM\xx$ are inductively determined
from
$\xx^{(0)}=\xx$ and $\h_1=1^\nn$ by
\begin{equation}
\label{E:GlobalIter}
\xx_\rr= 
\Tail{\xx^{(\rr-1)}}{\MM_{\cl{\h_\rr}}}
\text{\ , \quad}
\xx^{(\rr)}=\xx^{(\rr-1)}\OVER{\xx_\rr}
\text{\ , \quad and\quad}
\h_{\rr+1}=\Succ\mm{\h_\rr},
\end{equation}
where $\mm$ is the length of the longest prefix~$\h$
of~$\h_\rr$ that satisfies $\xx^{(\rr)} \not\perp\MM_{\cl\h}$.
\end{prop}

\begin{proof}
As can be expected, we use an induction on~$\nn$. The argument relies on the
transivity relation of Lemma~\ref{L:BiTail}.

For $\nn=0$, everything is
trivial, and, for $\nn=1$, the result is a restatement of
Proposition~\ref{P:AltDec}: in this case, the
$1$-address~$\h_\rr$ is~$\rr$, the longest
prefix of~$\h_\rr$ satisfying $\xx^{(\rr)} \not\perp
\MM_{\cl\h}$ is~$\ea$, and the induction  rule
reduces to $\h_{\rr+1}=\rr+1$.

Assume $\nn\ge2$. Let $(\yy_\qq, \Ldots, \yy_1)$ be the
$(\MM_2, \MM_1)$-decomposition of~$\xx$. By
definition, we have
\begin{equation}
\label{E:GlobalDec1}
\DD\MMM\xx=
(\DD{\MMM_{\cl\qq}}{\yy_\qq} \Ldots, 
\DD{\MMM_1}{\yy_1}).
\end{equation}
For $\qq\ge\jj\ge1$, let $(\yy_{\jj,\pp_\jj}, \Ldots,
\yy_{\jj,1})$ and
$(\h_{\jj,\pp_\jj}, \Ldots, \h_{\jj,1})$ be the unbracketing and
the address list in~$\DD{\MMM_{\cl\jj}}{\yy_\jj}$.
Then, by~\eqref{E:GlobalDec1}, we have
\begin{equation}
\label{E:Concat1}
(\xx_\pp, \Ldots, \xx_1) = (\yy_{\qq,\pp_\qq}, \Ldots,
\yy_{\qq,1}) \conc ... \conc (\yy_{1,\pp_1}, \Ldots,
\yy_{1,1}),
\end{equation}
where $\conc$ denotes concatenation, and, similarly,
\begin{equation}
\label{E:Concat2}
(\h_\pp, \Ldots, \h_1) = (\qq\h_{\qq,\pp_\qq}, \Ldots,
\qq\h_{\qq,1}) \conc ... \conc (1\h_{1,\pp_1}, \Ldots,
1\h_{1,1}).
\end{equation}
By induction hypothesis, the sequences of~$\yy_\jj$'s
and~$\h_{\jj,\kk}$'s satisfy the counterpart
of~\eqref{E:GlobalIter}, and we wish to
deduce~\eqref{E:GlobalIter}, \ie, dropping the
elements~$\xx^{(\rr)}$, to prove
$$\xx_\rr = \Tail{\xx_\pp ... \xx_\rr}{\MM_{\h_\rr}}
\text{\quad and\quad}
\h_{\rr+1}=\Succ\mm{\h_\rr}$$ where $\mm$ is the
length of the maximal prefix~$\h$ of~$\h_\rr$ satisfying
$\xx_\pp ... \xx_{\rr+1} \not\perp \MM_{\cl\h}$.
We use induction on $\rr\ge1$.

Assume that $\xx_\rr$ corresponds
to some entry~$\yy_{\jj,\kk}$ in~\eqref{E:Concat1}.  By
construction, we have $\h_\rr=\jj\h_{\jj,\kk}$.
Let $\yy =\yy_{\jj,\pp_\jj}...\yy_{\jj,\kk}$. The induction
hypothesis gives
\begin{equation}
\label{E:Tail1}
\xx_\rr = \yy_{\jj,\kk} =
\Tail{\yy}{\MM_{\cl{\jj\h_{\jj,\kk}}}} =
\Tail{\yy}{\MM_{\cl{\h_\rr}}}.
\end{equation}
On the other hand, by construction, $\yy$ is a left divisor of
$\yy_{\jj,\pp_\jj}...\yy_{\jj,1}$, \ie, of~$\yy_\jj$, and
$\yy_\jj$ is the $\MM_{\cl\jj}$-tail of
$\yy_\qq...\yy_\jj$, \ie, putting $\zz=\yy_\qq ...
\yy_\jj$, we have
\begin{equation}
\label{E:Tail2}
\yy_\jj = \Tail{\zz}{\MM_{\cl\jj}}.
\end{equation}
Applying Lemma~\ref{L:BiTail} to the
monoids $\MM_{\cl{\h_\rr}} \subseteq \MM_{\cl\jj}
\subseteq \MM$, we deduce from~\eqref{E:Tail1}
and~\eqref{E:Tail2} the relation
$\xx_\rr = \Tail{(\zz\OVER{\yy_\jj})\yy}{\MM_{\cl{\h_\rr}}}$,
which is $\xx_\rr = \Tail{\xx_\pp ...
\xx_\rr}{\MM_{\cl{\h_\rr}}}$, as, by construction,
we have $(\zz\OVER{\yy_\jj})\yy = \xx_\pp...\xx_\rr$.

Consider now~$\h_{\rr+1}$. Two
cases are possible, according to whether $\xx_\rr$
corresponds to an initial or a non-initial entry in
some sequence of~$\yy$'s, \ie, with the above notation,
according to whether $\kk = \pp_\jj$ holds or not. 
Assume first $\kk < \pp_\jj$. Then $\h_{\jj,\kk+1}$ exists, and
the induction hypothesis implies that
$\h_{\jj,\kk+1}$ is the $\mm$-successor of~$\h_{\jj,\kk}$,
where $\mm$ is the length of the maximal prefix~$\h$
of~$\h_{\jj,\kk}$ for which $\yy_{\jj,\pp_\jj}...
\yy_{\jj,\kk+1} \not\perp \MM_{\cl{\jj\h}}$ holds. The latter
relation is equivalent to $\xx_\pp ... \xx_{\rr+1} \not\perp
\MM_{\cl{\jj\h}}$: indeed,
$\xx\not\perp\AA$ is equivalent to $\Tail\xx\AA\not=1$,
and, as above, Lemma~\ref{L:BiTail} implies
$\Tail{\xx_\pp ... \xx_{\rr+1}}{\MM_{\cl{\jj\h}}}
= \Tail{\yy_{\jj,\pp_\jj}...
\yy_{\jj,\kk+1}}{\MM_{\cl{\jj\h}}}$. Therefore,
$\h_{\rr+1}$, which is $\jj\h_{\jj,\kk+1}$, is the
$\mm+1$-successor of~$\jj\h_{\jj,\kk}$, \ie, of~$\h_\rr$,
where $\mm$ is the length of the maximal prefix~$\h$
of~$\h_{\jj,\kk}$ for which $\xx_\pp ... \xx_{\rr+1} \not\perp
\MM_{\cl{\jj\h}}$ holds, hence $\mm+1$ is the length of the
maximal prefix~$\h'$ of~$\h_\rr$ (namely~$\jj\h$) for which 
$\xx_\pp ... \xx_{\rr+1} \not\perp \MM_{\cl{\h'}}$ holds.

Finally, assume $\kk = \pp_\jj$, \ie, $\h_{\jj,\kk}$ is the
leftmost address in the $\MMM_{\cl\jj}$-decomposition
of~$\yy_\jj$. In this case, by hypothesis, we have
$\h_{\rr+1}=(\jj+1)1^\nno$. Now, the hypothesis implies
$\yy_\qq...\yy_{\jj+1} \perp \MM_{\cl\jj}$, \ie,
$\xx_\pp...\xx_{\rr+1} \perp \MM_{\cl\jj}$. So, in
this case, the only prefix~$\h$ of~$\h_\rr$, \ie,
of~$\jj\h_{\jj,\pp_\jj}$, for which $\xx_\pp...\xx_{\rr+1}
\not\perp \MM_{\cl\h}$ may hold is the empty
address~$\ea$, which is the expected relation with
$\mm=0$. 
\end{proof}

\begin{exam}
\label{X:IterBis}
Consider the case of§ $\BB4$ and $\D_4^2$ again.
Proposition~\ref{P:GlobalIter} directly gives
the $\MMM$-decomposition of~$\D_4^2$ as follows. We
start with $\xx = \D_4^2$ and $\h_1=11$. Then we compute
$\MM_{11}$-tail, \ie, here the $\Mon{\ss1}$-tail,
of~$\xx^{(0)}$, which turns out to be~$\s_1^2$, and call
the quotient~$\xx^{(1)}$. Then the address~$\h_2$ is
obtained by looking at the maximal prefix~$\h$ of~$\h_1$,
\ie, of~$11$, for which $\MM_{\cl\h}\not\perp\xx^{(1)}$ holds.
In the current case, we have
$\xx^{(1)}\perp\MM_{11}$ and
$\xx^{(1)}\not\perp\MM_1$, hence $\h=1$, so $\h_2$ is
obtained from~$11$ by incrementing the second digit,
leading to~$\h_2=12$, which corresponds to
$\MM_{\cl{\h_2}}=\Mon{\ss2}$. We take the
$\Mon{\ss2}$-tail of~$\xx^{(1)}$, call the
remainder~$\xx^{(2)}$, and iterate. The successive values are
displayed in Table~\ref{T:Iter}.
\end{exam}

\begin{table}[htb]
\begin{tabular}{c|lcccc}
$\rr$\quad
&$\xx^{(\rr)}$
&$\h_\rr$
&$\cl{\h_\rr}$
&$\MM_{\cl{\h_\rr}}$
&$\xx_\rr$
\\
\hline
$0$
&$\ss1\ss2\ss1\ss3\ss2\ss1\ss1\ss2\ss1\ss3\ss2\ss1$
\rule{0pt}{12pt}
\\
$1$
&$\ss2\ss1\ss3\ss2\ss1\ss1\ss2\ss3\ss1\ss2$
&$11$&$11$
&$\Mon{\ss1}$
&$\s_1^2$
\\
$2$
&$\ss2\ss1\ss3\ss2\ss1\ss1\ss2\ss3\ss1$
&$12$&$12$
&$\Mon{\ss2}$
&$\ss2$
\\
$3$
&$\ss2\ss3\ss2\ss1\ss1\ss2\ss3$
&$13$&$11$
&$\Mon{\ss1}$
&$\s_1^2$
\\
$4$
&$\ss3\ss2\ss1\ss1\ss2\ss3$
&$14$&$12$
&$\Mon{\ss2}$
&$\ss2$
\\
$5$
&$\ss3\ss2\ss1\ss1\ss2$
&$21$&$21$
&$\Mon{\ss3}$
&$\ss3$
\\
$6$
&$\ss3\ss2\ss1\ss1$
&$22$&$22$
&$\Mon{\ss2}$
&$\ss2$
\\
$7$
&$\ss3\ss2$
&$31$&$11$
&$\Mon{\ss1}$
&$\s_1^2$
\\
$8$
&$\ss3$
&$32$&$12$
&$\Mon{\ss2}$
&$\ss2$
\\
$9$
&$1$
&$41$&$21$
&$\Mon{\ss3}$
&$\ss3$
\end{tabular}
\bigskip
\caption{\sf Direct determination of the iterated
decomposition of~$\D_4^2$: at step~$\rr$, we extract the
maximal right divisor~$\xx_\rr$ of the current
remainder~$\xx^{(\rr-1)}$ that lies in the
monoid~$\MM_{\cl{\h_\rr}}$, we update the remainder
into~$\xx^{(\rr)}$, and we define the next address~$\h_{\rr+1}$
to be the maximal successor~$\h$ of~$\h_\rr$ for which
$\xx^{(\rr)}$ is not orthogonal to~$\MM_{\cl\h}$; we stop
when only $1$ is left.}
\label{T:Iter}
\end{table}

\section{The alternating normal form}
\label{S:NormalForm}

We shall now deduce normal form results in (good) locally Garside monoids. The initial
remark is that, if $\MM$ is a locally Garside monoid generated by an element~$\gen$, then $\MM$
must be torsion-free by Condition~$(C_3)$, hence it is a free monoid, and every element of~$\MM$
admits a unique expression as~$\gen^\ee$ with~$\ee \in \Nat$. Now, if $\MM$ is an arbitrary locally right Garside
monoid and if $\MMM$ is an (iterated) covering of~$\MM$, then each element of~$\xx$ has been given a
distinguished decomposition in terms of the factor monoids~$\MM_\a$ of~$\MMM$. If, moreover, each of the
monoids~$\MM_\a$ happens to be generated by a single element~$\gen_\a$, the $\MMM$-decomposition gives a
unique distinguished expression in terms of the elements~$\gen_\a$. This situation
occurs for instance in the case of the $2$-covering of Example~\ref{X:IterCov}.

\subsection{Atomic coverings}

From now on, we consider locally right Garside monoids that satisfy Condition~$(C_3^+)$. It is easily seen that such
monoids are generated by atoms, \ie, elements~$\gen$ such that $\gen = \xx\yy$ implies $\xx = 1$ or $\yy =
1$---see for instance~\cite{Dfx}. In view of the above remarks, it is natural to concentrate on
coverings that involve submonoids generated by atoms.

\begin{defi}
\label{D:Atomic}
Assume that $\MM$ is a locally right Garside monoid, and $\GGen$ is an $\nn$-sequence of atoms of~$\MM$. We
say that an $\nn$-covering of~$\MM$ is \emph{atomic}  based on the sequence~$\GGen$ if, for each
$\nn$-address~$\a$, the monoid~$\MM_\a$ is the submonoid of~$\MM$ generated by the atom~$\Gen\a$.
\end{defi}

For instance, the $2$-covering of Example~\ref{X:IterCov} is atomic, based on $((\ss2,\ss3),(\ss2,\ss1))$. Note that a
base sequence must contain all atoms of~$\MM$, as, by definition, it generates~$\MM$. An arbitrary sequence of
atoms need not always define a covering, as a submonoid generated by a family of atoms is not necessarily closed in
the sense of Definition~\ref{D:Closed}. This however is true in
braid monoids---and in all Artin--Tits monoids. 

Before going on and defining the $\MMM$-normal form, we discuss one more general point, namely whether
$\MMM$-decompositions may have gap, this meaning that a trivial factor~$1$ may appear between two non-trivial
factors. 

\begin{exam}
\label{X:Gap}
Let $\MM$ be the $5$-strand braid monoid~$\BB5$, and
$\MMM$ be the $2$-covering based on~$((\ss4,\ss3),(\ss2,\ss1))$. One easily checks that the
$\MMM$-decomposition of~$\xx$ is $((\ss4, 1), (\ss1))$,
which has a trivial entry lying between two non-trivial entries.
\end{exam}

It is easy to state conditions that exclude such gaps.

\begin{lemm}
\label{L:Dense}
Say that an $\nn$-covering~$\MMM$ is  \emph{dense} if, for each binary
address~$\b$ of length~$\mm$ with $0\le\mm<\nn$,
\begin{equation}
\label{E:Dense}
\text{$\MM_\b$ is generated by 
$\MM_{\b1}$ and
$\MM_{\b21^{\nn-\mm-1}}$, and by $\MM_{\b2}$ and
$\MM_{\b1^{\nn-\mm}}$}. 
\end{equation}
Then, decomposi\-tions associated with a dense covering have no gap.
\end{lemm}

\begin{proof}
Owing to Proposition~\ref{P:GlobalIter}, the
point is to prove that, if, for some
binary $\nn$-address~$\a$ and some~$\mm$, 
writing $\b$ (\resp $\b'$) for the length~$\mm$ (\resp
$\mm+1$) prefix of~$\a$, we have both
$\xx \not\perp \MM_\b$ and
$\xx \perp \MM_{\b'}$, then necessarily the
$\MM_{\BSucc\mm\a}$-tail of~$\xx$ is not trivial. Write
$\b' = \b\rr$. For $\rr=1$, a sufficient condition for the
previous implication is that $\MM_{\b}$ is generated
by~$\MM_{\b1}$ and~$\MM_{\b21^{\nn-\mm-1}}$: then, a
non-trivial right divisor of~$\xx$ lying in~$\MM_\b$ cannot
be right divisible by any factor in~$\MM_{\b1}$ and,
therefore, it must be right divisible by some factor
in~$\MM_{\b21^{\nn-\mm-1}}$, and, by definition, we have
$\b21^{\nn-\mm-1}=\BSucc\mm\a$. For $\rr=2$, the
argument is similar, replacing~$\b1$ with~$\bŽ$, and
$\b21^{\nn-\mm-1}$ with $\b1^{\nn-\mm}$. So, the
conditions in~\eqref{E:Dense} are sufficient.
\end{proof}

In the case of an atomic covering, the density condition of Lemma~\ref{L:Dense} requires
that the base sequence be highly redundant. Such conditions are important in practice because they strongly limit
the patterns that can be used in the construction of dense atomic coverings.

\begin{prop}
\label{P:Successors}
Assume that $\MMM$ is a dense atomic
$\nn$-covering of~$\MM$ based on~$\GGen$. Then, for each
$\nn$-address~$\a$, the set
$\{\Gen{\BSucc\mm\a} \mid 0\le\mm\le\nobreak \nn\}$ is the atom
set of~$\MM$, and the latter contains at most
$\nn+1$ elements.
\end{prop}

\begin{proof}
Use induction on~$\nn\ge0$. The case $\nn=0$ is
obvious. Assume $\nn\ge1$. Write $\a=\dd\b$ with
$\dd=1$ or~$2$. Assume first $\dd=1$. By~\eqref{E:Dense},
$\MM$ is generated by $\Gen{21^\nno}$, which is the
$0$-successor of~$\a$, and~$\MM_1$. By induction
hypothesis, the latter is generated by the family of all
$\Gen{1\BSucc\mm\b}$'s, so $\MM$ is generated by the
successors of~$\a$. The argument is symmetric for $\dd=2$,
using the second part of~\eqref{E:Dense}. By
construction, every $\nn$-address admits $\nn+1$
successors, hence there are at most
$\nn+1$ atoms in~$\MM$.
\end{proof}

\begin{figure}[htb]
\begin{picture}(85,34)
\put(0,-0.5){\includegraphics{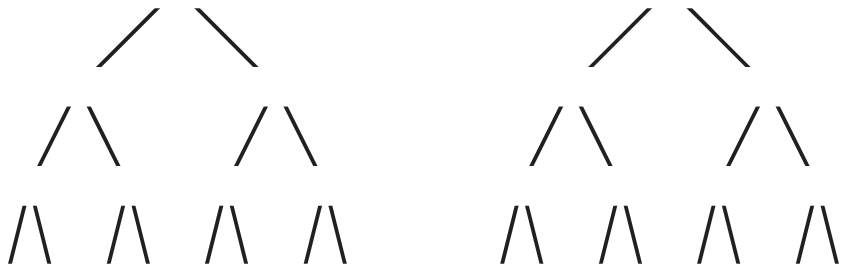}}
\put(0,0){$3$}
\put(5,0){$2$}
\put(10,0){$3$}
\put(15,0){$4$}
\put(20,0){$2$}
\put(25,0){$3$}
\put(30,0){$2$}
\put(35,0){$1$}
\put(2,10){$23$}
\put(12,10){$34$}
\put(22,10){$23$}
\put(32,10){$12$}
\put(6,20){$234$}
\put(26,20){$123$}
\put(15,30){$1234$}
\put(50,0){$2$}
\put(55,0){$3$}
\put(60,0){$2$}
\put(65,0){$4$}
\put(70,0){$2$}
\put(75,0){$3$}
\put(80,0){$2$}
\put(85,0){$1$}
\put(52,10){$23$}
\put(62,10){$24$}
\put(72,10){$23$}
\put(82,10){$12$}
\put(56,20){$234$}
\put(76,20){$123$}
\put(65,30){$1234$}
\end{picture}
\caption{\sf The two possible patterns for a dense
$3$-covering involving four atoms.}
\label{F:Solutions}
\end{figure}

We shall see in Section~\ref{S:Flip} that dense atomic
$\nn$-coverings involving $\nn\nobreak+\nobreak1$ atoms
exist for each~$\nn$. For $\nn=2$, the
only possible pattern is (up to renaming) that of Figure~\ref{F:Skeleton}. For
$\nn\ge3$, several non-isomorphic patterns exist---see Figure~\ref{F:Solutions}.

\subsection{The  $\MMM$-normal form}
\label{S:Normal}

We are now ready to convert the results of
Sections~\ref{S:Iter} into the
construction of a normal form. We recall that,
for~$\SS$ generating~$\MM$ and~$\ww$ a word on~$\SS$,
we denote by~$\CL\ww$ the element of~$\MM$ represented
by~$\ww$. We write~$\ww(\kk)$ for the $\kk$th letter
in~$\ww$ \emph{from the right}.

\begin{defi}
\label{D:NF}
Assume that $\MM$ is a locally right Garside monoid with atom set~$\SS$, and that $\MMM$ is a dense atomic
$\nn$-covering of~$\MM$ based on~$\GGen$. A
length~$\ell$ word~$\ww$ on~$\SS$ is said
to be \emph{$\MMM$-normal} if
\begin{quote}
There exist $\nn$-addresses $\a_\ell, \Ldots, \a_0$ with
$\a_0 = 1^\nn$ such that, for each~$\kk$, 
$\ww(\kk)=\Gen{\a_\kk}$ holds, where $\a_\kk$ is the
maximal successor of~$\a_{\kk-1}$---\ie, is
$\BSucc\mm{\a_{\kk-1}}$ with maximal~$\mm$---for which
$\Gen{\a_\kk}$ is a right divisor of
$\CL{\ww(\ell) ... \ww(\kk)}$.
\end{quote}
\end{defi}

The above definition may look convoluted, but handling a few
examples should make it easily understandable. 
Table~\ref{T:XNF} shows that our favourite example, namely 
$\ss3\ss2\ss1\ss1\ss2\ss3\ss2\ss1\ss1\ss2\ss1\ss1$,
is $\MMM$-normal with respect to the $2$-covering
of~Example~\ref{X:IterCov}.

The expected existence and uniqueness of the $\MMM$-normal form is the following easy result.

\begin{prop}
\label{P:NF}
Assume that $\MM$ is a locally right Garside monoid
with atom set~$\SS$, and $\MMM$ is a dense atomic
$\nn$-covering of~$\MM$ based on~$\GGen$.  Then each
element~$\xx$ of~$\MM$ admits a unique
$\MMM$-normal representative, namely
$\Gen{{\a_\ell}}...\Gen{{\a_1}}$, where
$\a_\ell, \Ldots, \a_1$ are inductively determined
from $\xx^{(0)}=\xx$ and $\a_0=1^\nn$ by
\begin{equation}
\label{E:NF2}
\a_\kk=\BSucc\mm{\a_{\kk-1}}
\text{\quad and \quad}
\xx^{(\kk)}= \xx^{(\kk-1)}\OVER{\Gen{{\a_\kk}}},
\end{equation}
where $\mm$ is maximal such that $\Gen{\BSucc\mm{\a_{\kk-1}}}$ is a right divisor of~$\xx^{(\kk)}$.
Moreover, $\Gen{{\a_\ell}}...\Gen{{\a_1}}$ is the word
obtained from the $\MMM$-decomposition of~$\xx$ by
concatenating the entries and possibly deleting the final~$1$. 
\end{prop}

\begin{proof}
The existence follows from the assumption that
$\MMM$ is dense, which guarantees that, as long as the
remainder~$\xx^{(\kk)}$ is not trivial, there must exist a
successor~$\BSucc\mm{\a_{\kk-1}}$ of the address~$\a_{\kk-1}$ such that
$\Gen{\BSucc\mm{\a_{\kk-1}}}$ is a right divisor of~$\xx^{(\kk)}$.
Uniqueness follows from the choice of that successor.

The inductive construction of~\eqref{E:NF2} is
essentially the construction of the $\MMM$-decomposition
as given in Proposition~\ref{P:GlobalIter}. The only
difference is that, here, we do not extract the whole tail of
the current remainder, but only one letter at each step. For
instance, if, at some point, the generator to be looked for
is~$\gen$ and the current remainder~$\xx^{(\kk-1)}$ is
divisible by~$\gen^2$, then
$\xx^{(\kk)}$ is~$\xx^{(\kk-1)}\OVER{\gen}$, and, at the next step,
$\a_\kk$ is the $\nn$-successor of~$\a_{\kk-1}$, \ie, it
is~$\a_{\kk-1}$ again, and the next letter of the normal form
is~$\gen$ again. In such a case, we have $\mm=\nn$.
By contrast, in Proposition~\ref{P:GlobalIter}, the
parameter~$\mm$ is never~$\nn$.
\end{proof}

Under the hypotheses of Proposition~\ref{P:NF}, the
word~$\ww$ is called the \emph{$\MMM$-normal form} of~$\xx$. The construction described in
Proposition~\ref{P:NF} is an algorithm, displayed in Table~\ref{T:Algo}.  A
typical example  is given in Table~\ref{T:XNF}.

\begin{table}[htb]
\begin{tabular}{l}
\hline
\rule{0pt}{12pt}{\bf Input:} A word~$\ww$ on~$\SS$;\\
\rule{0pt}{12pt}{\bf Procedure:}\\
\hspace{7mm}$\ww':=\mathtt{emptyword}$;\\
\hspace{7mm}$\a:= 1^\nn$;\\
\hspace{7mm}$\mathtt{while}\
\ww\not=\mathtt{emptyword} 
\ \mathtt{do}$\\
\hspace{14mm}$\mm:=\nn$;\\
\hspace{14mm}$\mathtt{while}\ \mathtt{quotient}
(\ww,\Gen{{\BSucc\mm{\a}}})\ =\mathtt{error}\ 
\mathtt{do}$\\
\hspace{21mm}$\mm:= \mm-1$;\\
\hspace{14mm}$\mathtt{od}$;\\
\hspace{14mm}$\a:=\BSucc\mm{\a}$;\\
\hspace{14mm}$\ww:= \mathtt{quotient}(\ww,\Gen{\a})$;\\
\hspace{14mm}$\ww':=\mathtt{concat}(\Gen{\a}, \ww')$;\\
\hspace{7mm}$\mathtt{od}$.\\
\rule[-6pt]{0pt}{12pt}{\bf Output:} The unique
$\MMM$-normal word~$\ww'$ that is equivalent
to~$\ww$.\\
\hline
\end{tabular}
\bigskip
\caption{\sf Algorithm for the $\MMM$-normal
form; we assume that $\SS$ is the atom set of~$\MM$, and
$\MMM$ is a dense atomic $\nn$-covering of~$\MM$
based on~$\GGen$; moreover, we assume
that $\mathtt{quotient}(\ww,\gen)$ is a subroutine that, for
$\ww$ a word on~$\SS$ and $\gen$ in~$\SS$,
returns $\mathtt{error}$ if $\gen$ is not a right divisor
of~$\CL\ww$, and returns a word
representing~$\CL\ww\OVER{\gen}$ otherwise.}
\label{T:Algo}
\end{table}

\begin{table}[htb]
\begin{tabular}{c|ccccccc}
$\kk$\quad
&$\ww_\kk$\hfill$\ww'_\kk$
\rule[-10pt]{0pt}{12pt}%
&\hbox{\quad$\a_{\kk-1}$\quad}
&$\mm$
&$\BSucc\mm{\a_{\kk-1}}$
&$\Gen{\BSucc\mm{\a_{\kk-1}}}$
&$\CL{\ww_\kk} \multe
\Gen{\BSucc\mm{\a_{\kk-1}}}?$
\\
\hline
\rule{0pt}{12pt}%
$0$
&$\ss1\ss2\ss1\ss3\ss2\ss1\ss1\ss2\ss1\ss3\ss2\ss1$
\qquad\hfill
-
&$11$
&$2$
&$11$
&$\ss1$
&yes\\
$1$
&$\ss1\ss2\ss1\ss3\ss2\ss1\ss1\ss2\ss1\ss3\ss2$
\qquad\hfill
$\ss1$
&$11$
&$2$
&$11$
&$\ss1$
&yes\\
$2$
&$\ss2\ss1\ss3\ss2\ss1\ss1\ss2\ss3\ss1\ss2$
\qquad\hfill
$\ss1\ss1$
&$11$
&$2$
&$11$
&$\ss1$
&no\\
&&
&$1$
&$12$
&$\ss2$
&yes\\
$3$
&$\ss2\ss1\ss3\ss2\ss1\ss1\ss2\ss3\ss1$
\qquad\hfill
$\ss2\ss1\ss1$
&$12$
&$2$
&$12$
&$\ss2$
&no\\
&&
&$1$
&$11$
&$\ss1$
&yes\\
$4$
&$\ss2\ss1\ss3\ss2\ss1\ss1\ss2\ss3$
\qquad\hfill
$\ss1\ss2\ss1\ss1$
&$11$
&$2$
&$11$
&$\ss1$
&yes\\
$5$
&$\ss2\ss3\ss2\ss1\ss1\ss2\ss3$
\qquad\hfill
$\ss1\ss1\ss2\ss1\ss1$
&$11$
&$2$
&$11$
&$\ss1$
&no\\
&&
&$1$
&$12$
&$\ss2$
&yes\\
$6$
&$\ss3\ss2\ss1\ss1\ss2\ss3$
\qquad\hfill
$\ss2\ss1\ss1\ss2\ss1\ss1$
&$12$
&$2$
&$12$
&$\ss2$
&no\\
&&
&$1$
&$11$
&$\ss1$
&no\\
&&
&$0$
&$21$
&$\ss3$
&yes\\
$7$
&$\ss3\ss2\ss1\ss1\ss2$
\qquad\hfill
$\ss3\ss2\ss1\ss1\ss2\ss1\ss1$
&$21$
&$2$
&$21$
&$\ss3$
&no\\
&&
&$1$
&$22$
&$\ss2$
&yes\\
$8$
&$\ss3\ss2\ss1\ss1$
\qquad\hfill
$\ss2\ss3\ss2\ss1\ss1\ss2\ss1\ss1$
&$22$
&$2$
&$22$
&$\ss2$
&no\\
&&
&$1$
&$21$
&$\ss3$
&no\\
&&
&$0$
&$11$
&$\ss1$
&yes\\
$9$
&$\ss3\ss2\ss1$
\qquad\hfill
$\ss1\ss2\ss3\ss2\ss1\ss1\ss2\ss1\ss1$
&$11$
&$2$
&$11$
&$\ss1$
&yes\\
$10$
&$\ss3\ss2$
\qquad\hfill
$\ss1\ss1\ss2\ss3\ss2\ss1\ss1\ss2\ss1\ss1$
&$11$
&$2$
&$11$
&$\ss1$
&no\\
&&
&$1$
&$12$
&$\ss2$
&yes\\
$11$
&$\ss3$
\qquad\hfill
$\ss2\ss1\ss1\ss2\ss3\ss2\ss1\ss1\ss2\ss1\ss1$
&$12$
&$2$
&$12$
&$\ss2$
&no\\
&&
&$1$
&$11$
&$\ss1$
&no\\
&&
&$0$
&$21$
&$\ss3$
&yes\\
$12$
&-
\qquad\hfill
$\ss3\ss2\ss1\ss1\ss2\ss3\ss2\ss1\ss1\ss2\ss1\ss1$
&21
&-
&-
&-\\
\end{tabular}
\bigskip
\caption{\sf Computation of the
$\MMM$-normal form of~$\D_4^2$, for
$\MMM$ the
$2$-covering of Example~\ref{X:Iter},
starting from the word
$(\ss1\ss2\ss1\ss3\ss2\ss1)^2$: at
each step, we try to divide the current word~$\ww_\kk$ by
some generator~$\ss\rr$ and, when succesful, we add
this~$\ss\rr$ on the left of~$\ww'_\kk$, until no letter is left
in~$\ww_\kk$; the point is to know in which order the
generators are tried, and this is specified by the
address~$\a_\kk$: we try the successors
of~$\a_{\kk-1}$ starting with the last one, \ie,
with~$\a_{\kk-1}$, and then consider shorter and
shorter prefixes of~$\a_{\kk-1}$; density guarantees that we
cannot get stuck until $\ww_\kk$ is empty.}
\label{T:XNF}
\end{table}

As for complexity, computing the $\MMM$-normal
form is as easy as computing the $\MMM$-decomposition.
In our current atomic context, the existence of the norm
(Definition~\ref{D:Norm}) is guaranteed~\cite{Dfx}.

\begin{prop}
\label{P:Complexity2}
Assume that $\MM$ is a locally right Garside monoid with atom set~$\SS$, that $\MMM$ is a dense atomic
$\nn$-covering of~$\MM$ based on~$\GGen$, and that Condition~$(*)$ of
Proposition~\ref{P:Complexity1} is satisfied. Then, for  each
word~$\ww$ on~$\SS$, the algorithm of
Table~\ref{T:Algo} runs in time~$O(\norm\ww^2)$.
\end{prop}

\begin{proof}
The only change with respect to
Proposition~\ref{P:Complexity1} is that we have to keep
track of binary addresses of fixed length~$\nn$ so as to
know in which order the divisions have to be tried. Getting a
new letter of the normal word under construction requires at
most $\nn+1$~divisions, but the rest is similar.
\end{proof}

\subsection{The exponent sequence}

We conclude this section with an easy remark about  $\MMM$-decompositions in the context of atomic coverings,
namely that an element of the monoid is non-ambiguously determined by the iterated sequence of exponents in its
$\MMM$-decomposition, \ie, we can forget about names of atoms and only keep track of exponents
without losing information.

\begin{defi}
For $\MM, \MMM$ as in Definition~\ref{D:NF}, and for $\www$ an iterated sequence whose entries are of the
form~$\Gen\a^{\ee_\a}$, we define the \emph{exponent sequence}~$\exp\www$ of~$\www$ to be the iterated
sequence obtained by replacing~$\Gen\a^{\ee_\a}$ with~$\ee_\a$ everywhere in~$\www$. 
\end{defi}

For instance, in the context of Example~\ref{X:Iter}, the $\MMM$-decomposition of~$\D_4^2$ is the $2$-sequence 
$((\ss3), (\ss2,\s_1^2), (\ss2,\ss3), (\ss2, \s_1^2,\ss2,\s_1^2))$, so the exponent
sequence is the $2$-sequence of natural numbers
$$((1), (1,2), (1,1), (1,2,1,2)).$$
As in the case of every iterated sequence,  specifying the exponent sequence of~$\DD\MMM\xx$ amounts to
giving two ordinary sequences, namely its unbracketing---in the above
example $(1, 1,2, 1,1, 1,2,1,2)$---and its address list---$(41,32,31,22,21$, $14,13, 12,11)$ above. Easy examples show
that, taken separately, neither of the above sequences is sufficient to recover~$\xx$. But, when we take them
simultaneously, we can recover~$\xx$.

\begin{prop}
\label{P:Exponent}
If $\MMM$ is an atomic $\nn$-covering of~$\MM$, then, for
every~$\xx$ in~$\MM$, the exponent sequence of~$\DD\MMM\xx$
determines~$\xx$.
\end{prop}

\begin{figure}[htb]
\begin{picture}(42,23)
\put(0.5,3){\includegraphics{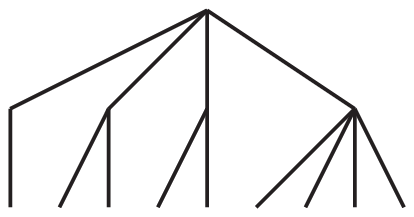}}
\put(0,0){$1$}
\put(5,0){$1$}
\put(10,0){$2$}
\put(15,0){$1$}
\put(20,0){$1$}
\put(25,0){$1$}
\put(30,0){$2$}
\put(35,0){$1$}
\put(40,0){$2$}
\end{picture}
\caption{\sf Tree representation of  the exponent
sequence of~$\DD\MMM{\D_4^2}$, \ie, of  $((1), (1,2), (1,1),
(1,2,1,2))$; Proposition~\ref{P:Exponent} states that the
geometry of the tree determines the missing names: for
instance, the leftmost~$2$ has address~$31$ in the tree, so
it corresponds to the generator~$\Gen{\cl{31}}$,
which is~$\ss1$; hence, this entry~$2$ must correspond to a
factor~$\s_1^2$ in $\DD\MMM{\D_4^2}$.}
\label{F:TreeNat}
\end{figure}

\begin{proof}
Let~$\GGen$ be the base sequence of~$\MMM$, and let
$(\ee_\pp, \Ldots, \ee_1)$ and $(\h_\pp, \Ldots, \h_1)$ be
the unbracketing and the address list in the exponent sequence of~$\DD\MMM\xx$. Then
we recover~$\DD\MMM\xx$ itself, and
therefore~$\xx$, by replacing for each~$\rr$ the
entry~$\ee_\rr$ corresponding to an address~$\h_\rr$
with~$\Gen{\cl{\h_\rr}}^{\ee_\rr}$. The formal proof is an easy induction on the degree of the
covering~$\MMM$---see Figure~\ref{F:TreeNat} for an
example.
\end{proof}

\section{The $\Phi$-normal form of braids}
\label{S:Flip}

From now on, we concentrate on the specific case of braids. In order to apply the previous results, we fix for
each~$\nn$ a covering of~$\BB\nn$ by two copies of~$\BB\nno$, namely
$\BB\nno$ and its image under the flip automorphism~$\ff\nn$. We study the decomposition associated with this
covering, as well as an iterated version and the derived normal form, called the
$\Phi$-normal form. This naturally leads to introducing a certain linear
ordering of~$\BB\nn$, which will be subsequently proved to be connected with the standard braid ordering.

\subsection{The $\Phi$-splitting of a braid}

In the sequel, we always  consider $\BB\nno$ as a submonoid of~$\BB\nn$: an
$(\nn-1)$-strand braid is a particular $\nn$-strand braid. 
We denote by~$\ff\nn$ the flip automorphism of~$\BB\nn$ that exchanges~$\ss\ii$ and~$\ss{\nn-\ii}$ for
each~$\ii$. It is well-known---see for instance~\cite[Chapter~1]{Dhr}---that $\ff\nn$ is the conjugation by the
Garside element~$\D_\nn$. We also use~$\ff\nn$ for $\nn$-strand braid words, thus denoting by~$\ff\nn(\ww)$
the image of a braid word~$\ww$ under~$\ff\nn$ letter by letter.  

The initial, obvious observation is that, for each~$\nn\ge 3$,
the monoids~$\BB\nno$ and~$\ff\nn(\BB\nno)$ are closed submonoids
of~$\BB\nn$, and that the pair $(\ff\nn(\BB\nno), \BB\nno)$ is a
covering of~$\BB\nn$ in the sense of
Definition~\ref{D:Covering}. Thus Proposition~\ref{P:AltDec}
gives for every $\nn$-strand braid a distinguished
decomposition as an alternating product of elements
of~$\BB\nno$ and~$\ff\nn(\BB\nno)$, according to the scheme of Figure~\ref{F:Alt}. We now
restate the general result so as to emphasize the r\^ole of the flip
automorphism.

\begin{prop}
\label{P:FSplitting}
Every braid~$\xx$ in~$\BB\nn$ admits a unique decomposition
\begin{equation}
\label{E:Flip}
\xx=\ff\nn^\ppo(\xx_\pp) \cdot ... \cdot \ff\nn(\xx_2) \cdot \xx_1
\end{equation}
with  $\xx_1, \Ldots, \xx_\pp$ in~$\BB\nno$
such that, for each~$\rr\ge2$, the only $\ss\ii$
that is a right divisor of
$\ff\nn^{\pp-\rr}(\xx_\pp) \cdot ... \cdot \ff\nn(\xx_{\rr+1})
\cdot \xx_\rr$ is~$\ss1$. The braids~$\xx_\rr$
are determined from
$\xx^{(0)}= \xx$ by
\begin{equation}
\label{E:FlipDec2}
\xx_\rr= 
\Tail{\xx^{(\rr-1)}}{\BB\nno}, 
\quad
\xx^{(\rr)}=\ff\nn(\xx^{(\rr-1)}\OVER{\xx_\rr}). 
\end{equation}
\end{prop}

\begin{proof}
As $\ff\nn$ is an automorphism of~$\BB\nn$, the relation
$\yy_1=\Tail\yy{\ff\nn(\BB\nno)}$ is equivalent to
$\ff\nn(\yy_1) =\Tail{\ff\nn(\yy)}{\BB\nno}$. Moreover
$\ff\nn$ is an automorphism for the quotient operation
as well. Then \eqref{E:Flip} and the divisibility constraints
just express that the sequence 
$(\ff\nn^\ppo(\xx_\pp), \Ldots, \ff\nn(\xx_2), \xx_1)$ is the
$(\ff\nn(\BB\nno), \BB\nno)$-decomposition of~$\xx$.
\end{proof}

\begin{defi}
\label{D:Splitting}
The sequence $(\xx_\pp, \Ldots, \xx_1)$ involved
in~\eqref{E:Flip} is called the \emph{$\nn$-splitting}
of~$\xx$; the parameter~$\pp$ is called the \emph{$\nn$-breadth of~$\xx$}.
\end{defi}

The only difference between the $(\ff\nn(\BB\nno),
\BB\nno)$-decomposition and the $\nn$-splitting is that the
flip~$\ff\nn$ is applied to each other entry. The benefit is that
all entries in the $\nn$-splitting of a braid of~$\BB\nn$ are braids of~$\BB\nno$, and not elements of~$\BB\nno$
and~$\ff\nn(\BB\nno)$, alternately.  Note that the
$\nn$-splitting of~$\xx$ is obtained by repeating a
single operation, namely finding the $\BB\nno$-tail of~$\xx$---hence the right gcd of~$\xx$
and~$\D_\nno^\infty$ as was seen in
Example~\ref{X:Closed}---and flipping the quotient.

\begin{exam}
\label{X:Splitting}
Let $\xx$ be the $4$-strand braid $\D_4^2$. The $\BB3$-tail of~$\xx$ is $\D_3^2$, with associated
quotient $\ss3\ss2\s_1^2\ss2\ss3$, hence, after a flip,
$\xx^{(1)} = \ss1\ss2\s_3^2\ss2\ss1$. The $\BB3$-tail of~$\xx^{(1)}$ is $\ss2\ss1$,
with quotient $\ss1\ss2\s_3^2$, hence, after a flip,
$\xx^{(2)} = \ss3\ss2\s_1^2$. The $\BB3$-tail of~$\xx^{(2)}$ is $\ss2\s_1^2$, with
quotient~$\ss3$, hence, after a flip, $\xx^{(3)} = \ss1$, which belongs to~$\BB3$. Thus $\D_4^2$ has
$4$-breadth~$4$, and its
$4$-splitting is $(\ss1, \ss2\s_1^2, \ss2\ss1,
\D_3^2)$---compare with the
$(\ff4(\BB3), \BB3)$-decomposition of~$\D_4^2$ as computed
in Example~\ref{X:Delta}.
\end{exam}

Note that, as in the case of the $(\ff\nn(\BB\nno), \BB\nno)$-decomposition, the non-final entries in an
$\nn$-splitting are never~$1$, but the final (rightmost) entry may:
the $3$-splitting of~$\ss2$ is $(\ss1, 1)$, as $\ss2$ is not divisible by~$\ss1$.

\subsection{The flip covering of~$\BB\nn$}

The $\nn$-splitting operation associates with every braid of~$\BB\nn$ a sequence of braids of~$\BB\nno$. We
can now iterate the construction, so as to associate with every braid of~$\BB\nn$ an iterated sequence of braids
of~$\BB2$. According to the general framework of Section~\ref{S:Iter}, this entails introducing an
iterated $(\nnt)$-covering of the monoid~$\BB\nn$.

\begin{defi}
For $\nn\ge2$, we denote~$\BBB\nn$ the $(\nn\!-\!2)$-covering of~$\BB\nn$ defined by
\begin{equation}
\label{E:DefFlipCov}
\BBB2=\BB2, \quad
\BBB\nn=(\ff\nn(\BBB\nno), \BBB\nno).
\end{equation}
\end{defi}

Applying the recursive definition, we find
\begin{gather*}
\BBB3 = (\ff3(\BB2), \BB2) = (\Mon{\ss2}, \Mon{\ss1}),\\
\BBB4 = (\ff4(\BBB3), \BBB3) = ((\Mon{\ss2}, \Mon{\ss3}), (\Mon{\ss2}, \Mon{\ss1})), 
\end{gather*}
which is the $2$-covering of Example~\ref{X:IterCov}. More generally, writing $\BB{\nn,\a}$ for the $\a$-entry
in~$\BBB\nn$, we deduce from~\eqref{E:DefFlipCov} the rules
\begin{equation}
\label{E:FlipCov4}
\BB{2,\ea}=\BB2, \quad
\BB{\nn,1\a}=\BB{\nn-1,\a},\quad \text{and} \quad
\BB{\nn,2\a}=\ff{\nn}(\BB{\nn-1,\a}).
\end{equation}

The above values show that $\BBB3$ and~$\BBB4$ are dense atomic coverings. This result extends to all values
of~$\nn$, with the following description of the base sequence.

\begin{prop}
\label{P:BraidCovering}
For $\nn \ge 2$, define the $(\nnt)$-sequence~$\GGen_\nn$  by
\begin{equation}
\label{E:FlipCov}
\GGen_2=\ss1, \quad
\GGen_\nn= (\ff\nn(\GGen_\nno), \GGen_\nno).
\end{equation}
Then, for each  binary address~$\a$ of length $\nn\!-\!2$, we have
$\Gen{\a} =\ss\ii$ with
\begin{equation}
\label{E:Name}
\ii = -\mm_1 + \mm_2 - ... +
(-1)^\rr \mm_\rr +
\begin{cases}
1&\text{if $\rr$ is even},\\
\nn&\text{if $\rr$ is odd},
\end{cases}
\end{equation}
if $\a = \dd_1...\dd_{\nn-2}$ and $\mm_1 < ... <
\mm_\rr$ are the~$\mm$'s for which $\dd_\mm$ is even. 
Moreover, $\BBB\nn$ is a dense atomic covering based
on~$\GGen_\nn$. 
\end{prop}

\begin{proof}
Firstly, we prove~\eqref{E:Name} using induction on~$\nn \ge 2$. For $n= 2$, \eqref{E:Name}
reduces to $\Gen\ea=\ss1$, which is true. Assume $\nn\ge3$, and let
$\a'=\dd_2...\dd_{\nn-2}$. Putting $\Gen{\a'}=\ss{\ii'}$,
we aim at proving $\ii=\ii'$ if $\dd_1$ is odd, and
$\ii=\nn-\ii'$ if $\dd_1$ is even.  Write $\SS$ for $-\mm_1 +
\mm_2 - ... + (-1)^\rr
\mm_\rr$, and $\rr'$, $\mm'_1$, $\mm'_2$, \Ldots, $\SS'$,
$\nn'$ for the similar parameters associated with~$\a'$.
Assume first that $\dd_1$ is odd. Then we have $\rr=\rr'$,
and $\mm_\jj=\mm'_\jj+1$ for each~$\jj$, hence
$\SS=\SS'$ if $\rr$ is odd, and
$\SS=\SS'-1$ if $\rr$ is even. The induction
hypothesis gives $\ii'=\SS'+1$ if $\rr$ is odd,
$\SS'+\nn'$ if $\rr$ is even. We deduce
$\ii= \SS+1 = \SS'+1=\ii'$ if $\rr$ is even, and 
$\ii= \SS+\nn = \SS'-1+\nn'+1=\ii' $ if $\rr$ is odd.

Assume now that $\dd_1$ is even. Then we  have
$\rr=\rr'+1$, $\mm_1=1$, and $\mm_{\jj+1} =\mm'_\jj+1$
for each~$\jj\ge1$, hence $\SS=-\SS'$ if $\rr$ is
odd, and $\SS=-\SS'-1$ if $\rr$ is even. The induction
hypothesis gives $\ii'=\SS'+\nn'$ if $\rr$ is
odd, $\SS'+1$ if $\rr$ is even. We deduce
$\ii = \SS+1 = - \SS'+1= \nn - \ii'$ if $\rr$ is odd, and $\ii = 
\SS+\nn = - \SS'-1+\nn=\nn-\ii'$ if $\rr$ is
even.

Nex, the braids~$\ss\ii$ are the atoms of~$\BB\nn$, and every parabolic submonoid of~$\BB\nn$ is closed, so
every surjective sequence of atoms defines a covering. An obvious induction on~$\nn$ shows that, for~$\nn \ge 2$,
each of $\ss1, \Ldots, \ss\nno$ occurs in the sequence~$\GGen_\nn$. Moreover, comparing~\eqref{E:DefFlipCov}
and~\eqref{E:FlipCov} makes it straightforward that $\BBB\nn$ is precisely the covering based on~$\GGen_\nn$.

As for density, the point is
to show that $\BB\nn$ is generated by~$\BB\nno$
and~$\BB{\nn,21^{\nn-3}}$. Now \eqref{E:Name} gives
$\Gen{21^{\nn-3}} = \ss\nno$, precisely
the atom of~$\BB\nn$ missing in~$\BB\nno$.
\end{proof}

It is easy to see that, for each~$\nn$, the unbracketing of ~$\GGen_\nn$ is the length~$2^{\nn-2}$ suffix
of some left infinite sequence $\GGen_\infty$ where indices are
$$\quad...,6,3,4,3,2,4,3,4,5,3,2,3,4,2,3,2,1.$$
An example of application for the rule of~\eqref{E:Name} is as follows: in the length~$7$ address~$1221212$, there
are even digits at positions $2,3,5,7$ (from the left), so \eqref{E:Name}
gives $\ii=(-2+3-5+7)+1=4$, hence $\Gen{1221212}=\ss4$. 

As $\BBB\nn$ is a dense atomic covering of~$\BB\nn$, it is eligible for the results of Section~\ref{S:Iter}. We fix
some specific, simplified notation.

\begin{nota}
For~$\xx$ in~$\BB\nn$, the $\BBB\nn$-decomposition of~$\xx$ is denoted by~$\DDf\nn\xx$, and its exponent
sequence is denoted by~$\DDfe\nn\xx$.
\end{nota}

The recursive definition of~$\BBB\nn$  implies the following connection between the splitting and the
$\BBB\nn$-decomposition.

\begin{lemm}
\label{L:DecompBraid}
For $\nn \ge 3$ and  $\xx$ in~$\BB\nn$, we have
\begin{equation}
\label{E:BiFlip21}
\DDf\nn\xx=(\ff\nn^\ppo(\DDf\nno{\xx_\pp}), \Ldots,
\ff\nn(\DDf\nno{\xx_2}), \DDf\nno{\xx_1}).
\end{equation}
where $(\xx_\pp, \Ldots, \xx_1)$ is the $\nn$-splitting of~$\xx$.
\end{lemm}

\begin{proof}
By definition, the $(\ff\nn(\BB\nno), \BB\nno)$-decomposition of~$\xx$ is the sequence
$$(\ff\nn^\ppo(\xx_\pp), \Ldots, \ff\nn(\xx_2), \xx_1),$$
and, therefore, by definition again, we have
\begin{equation*}
\DDf\nn\xx = 
(\DD{\ff\nn^\ppo(\BBB\nno)}{\ff\nn^\ppo(\xx_\pp)},
\Ldots, 
\DD{\ff\nn(\BBB\nno)}{\ff\nn(\xx_2)}, 
\DD{\BBB\nno}{\xx_1}).
\end{equation*}
Now, as $\ff\nn$ is an automorphism of~$\BB\nn$, we have
$\DD{\ff\nn(\BBB\nno)}{\ff\nn(\yy)} = \ff\nn(\DD{\BBB\nno}{\yy})$
for each~$\yy$ in~$\BB\nno$, \ie, $\DD{\ff\nn(\BBB\nno)}{\ff\nn(\yy)} = \ff\nn(\DDf\nno\yy)$, and 
\eqref{E:BiFlip21} follows.
\end{proof}

\begin{figure}[htb]
\begin{picture}(43,35)
\put(-1,12){\includegraphics{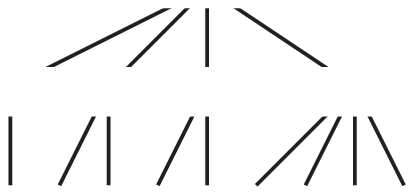}}
\put(0,9){$\ss1$}
\put(5,9){$\ss1$}
\put(10,9){$\s_1^2$}
\put(15,9){$\ss1$}
\put(20,9){$\ss1$}
\put(25,9){$\ss1$}
\put(30,9){$\s_1^2$}
\put(35,9){$\ss1$}
\put(40,9){$\s_1^2$}
\put(0,21){$\ss2$}
\put(7,21){$\ss2\s_1^2$}
\put(17,21){$\ss1\ss2$}
\put(33,21){$\D_3^2$}
\put(19,31.5){$\D_4^2$}
\put(6,27){$\scriptstyle\ff{\!4}$}
\put(18.2,27){$\scriptstyle\ff{\!4}$}
\put(4.5,15){$\scriptstyle\ff{\!3}$}
\put(14.2,15){$\scriptstyle\ff{\!3}$}
\put(25.5,15){$\scriptstyle\ff{\!3}$}
\put(33.2,15){$\scriptstyle\ff{\!3}$}
\put(-40,4){\ie, after reintroducing the flips,}
\put(0,0){$\ss3$}
\put(5,0){$\ss2$}
\put(10,0){$\s_1^2$}
\put(15,0){$\ss2$}
\put(20,0){$\ss3$}
\put(25,0){$\ss2$}
\put(30,0){$\s_1^2$}
\put(35,0){$\ss2$}
\put(40,0){$\s_1^2$}
\put(0,0){$\ss3$}
\put(5,0){$\ss2$}
\end{picture}
\caption{\sf The $\BBB4$-decomposition of~$\D_4^2$ viewed as an iterated
splitting: we split the initial braid of~$\BB4$ into a sequence
of braids in~$\BB3$, then we split each of them into a sequence of
braids in~$\BB2$, \ie, of powers of~$\ss1$; the sequence $\DDf4{\D_4^2}$ is obtained by
iteratively flipping each other entry.}
\label{F:Splitting}
\end{figure}

\begin{exam}
\label{X:FlipNormal}
(See Figure~\ref{F:Splitting})
We saw in Example~\ref{X:Splitting} that the $4$-splitting
of~$\D_4^2$ is $(\ss1, \ss2\s_1^2, \ss2\ss1, \D_3^2)$.
Now, the $3$-splitting of~$\D_3^2$ turns out to be
$(\ss1, \s_1^2, \ss1, \s_1^2)$, that of $\ss2\ss1$ is $(\ss1,
\ss1)$, etc. Gathering the results, and applying the needed
flips, we find
\begin{equation}
\label{E:XFlip}
\DDf4{\D_4^2}=((\ss3), (\ss2,\s_1^2), (\ss2,\ss3), 
(\ss2, \s_1^2,\ss2,\s_1^2)),
\end{equation}
as already seen in Example~\ref{X:Iter}. The associated exponent sequence is
\begin{equation}
\DDfe4{\D_4^2}=((1), (1,2), (1,1), (1,2,1,2)),
\end{equation}
\end{exam}

\subsection{The $\Phi$-normal form}

The iterated covering~$\BBB\nn$ is atomic and, therefore, it gives raise to a unique normal form on~$\BB\nn$. 
According to Proposition~\ref{P:NF}, the $\BBB\nn$-normal form of a braid~$\xx$ of~$\BB\nn$ is the word
obtained by concatenating the (unique) expressions of the successive entries in its $\BBB\nn$-decomposition as
powers of atom. For instance, from the $\BBB4$-decomposition of~$\D_4^2$ given in~\eqref{E:XFlip}, we deduce
the $\BBB4$-normal form $\ss3 \ss2 \s_1^2 \ss2 \ss3 \ss2 \s_1^2 \ss2 \s_1^2$.

If $\xx$ belongs to~$\BB\nno$, then the $\nn$-splitting of~$\xx$ is
the length one sequence~$(\xx)$. Therefore, we have $\DD\nn\xx = (\DD\nno\xx)$, and the normal
form of~$\xx$ as an element of~$\BB\nno$ coincides with its normal form as an element of~$\BB\nn$. Owing to
this remark, we shall forget about subscripts, and put the following without ambiguity.

\begin{defi}
For $\xx$ in~$\BB\nn$, the $\BBB\nn$-normal form of~$\xx$ is called the \emph{$\Phi$-normal form} of~$\xx$.
\end{defi}

Lemma~\ref{L:DecompBraid} implies that the $\Phi$-normal form has the following simple
connection with the splitting operation---which could be taken as an alternative definition:

\begin{prop}
For $\nn \ge 3$ and  $\xx$ in~$\BB\nn$, the $\Phi$-normal form of~$\xx$ is the word
\begin{equation}
\label{E:NFBraid}
\ff\nn^\ppo(\ww_\pp) \cdot ...\cdot \ff\nn(\ww_2) \cdot \ww_1,
\end{equation}
where $(\xx_\pp, \Ldots, \xx_1)$ is the $\nn$-splitting of~$\xx$, and, for each~$\rr$, the word~$\ww_\rr$ is the
$\Phi$-normal form of~$\xx_\rr$.
\end{prop}

The results of Section~\ref{S:Normal} imply that, in addition to the above recursive definitions, the $\Phi$-normal form
also admits direct characterizations. We shall now state such characterizations. Several
equivalent statements are possible---and can be used in practical implementations. The principle is always:
\begin{quote}
An $\nn$-strand braid word~$\ww$ is $\Phi$-normal if, for
each~$\kk$, the $\kk$th letter of~$\ww$ starting from the
right is the smallest~$\ss\ii$ that is a right divisor of the
braid represented by the prefix of~$\ww$ finishing at that
letter, smallest referring to some local ordering of
the~$\ss\ii$'s that is updated at each step and corresponds to
a position in the skeleton of the covering~$\BBB\nn$.
\end{quote}
The formal definition includes a description of the
local ordering of the~$\ss\ii$'s. The latter can be encoded in
several equivalent ways, involving addresses, or numbers, or
permutations. If the local ordering were the
fixed order $\ss1\!< \!...\!< \!\ss\nno$, then being normal
would simply mean being lexicographically minimal. 

We recall that, for $\a$ a binary address,
$\BSucc\mm\aa$ denotes the binary $\mm$-successor
of~$\a$ (Definition~\ref{D:Successor}), and that, for $\ww$ a braid word, $\CL\ww$
denotes the braid represented by~$\ww$.

\begin{prop}
\label{P:NormalBraid}
A length~$\ell$ positive $\nn$-strand braid
word~$\ww$ is $\Phi$-normal if and only if any
one of the following equivalent conditions holds:

$(i)$ There exist binary addresses  $\a_\ell, \Ldots,
\a_0$ with $\a_0 = 1^{\nn-2}$ such that,
for each~$\kk$, $\ww(\kk) =
\Gen{\a_\kk}$ holds, and $\a_\kk$ is
the maximal binary successor of~$\a_{\kk-1}$ such that
$\Gen{\a_\kk}$ is a right divisor of~$\CL{\ww(\ell)  ... \ww(\kk)}$.

$(ii)$ There exist numbers $\mm_\ell, \Ldots, \mm_1$ in
$\{0, \Ldots, \nn\}$ such that, putting $\a_0= 1^{\nn-2}$
and inductively defining  $\a_\kk=
\BSucc{\mm_\kk}{\a_{\kk-1}}$, then, for each~$\kk$, we
have $\ww(\kk)=\Gen{\a_\kk}$ and $\CL{\ww(\ell)  ... \ww(\kk)}
\not\multe\Gen{\a}$ for every $\mm$-successor~$\a$ of~$\a_{\kk-1}$
with $\mm>\mm_\kk$.

$(iii)$ There exist permutations $\perm_\ell, \Ldots,
\perm_0$ of~$\{1, \Ldots, \nn-1\}$ such that $\perm_0$ is
the identity, and, for each~$\kk$, we have
$\ww(\kk) = \ss{\perm_\kk(1)}$ and $\perm_\kk$ is obtained
from~$\perm_{\kk-1}$ as follows: let $\pp$ be minimal
satisfying~$\CL{\ww(\ell)  ... \ww(\kk)}
\multe \ss{\perm_{\kk-1}(\pp)}$; then we have
$\perm_\kk(1)= \perm_{\kk-1}(\pp)$, 
$\perm_\kk(\qq)=\perm_{\kk-1}(\qq)$ for $\qq>\pp$,
and $(\perm_\kk(2), \Ldots, \perm_\kk(\pp))$ is
the increasing (\resp decreasing) enumeration of
$\{\perm_{\kk-1}(1), \Ldots, \perm_{\kk-1}(\pp-1)\}$ if the
latter are larger (\resp smaller) than~$\perm_{\kk}(1)$ in the usual ordering of integers.
\end{prop}

\begin{proof}
Point~$(i)$ is Definition~\ref{D:NF} and~$(ii)$ is a direct
reformulation. As for~$(iii)$, $\perm_\kk$ is the
enumeration of the names of the successors of~$\a_\kk$,
starting from the bottom, \ie, for each~$\mm$, we have
$\Gen{\BSucc\mm{\a_\kk}} =
\ss\ii$ with $\ii = \perm_\kk(\nn-\mm-1)$. At each step, we
select the maximal successor satisfying the divisibility
requirement, hence, here, the first entry in the
permutation~$\perm_{\kk-1}$; the updating rules come
from the specific definition of the covering~$\BBB\nn$. 
\end{proof}

As for complexity, a direct application of
Proposition~\ref{P:Complexity2} gives:

\begin{prop}
\label{P:BraidNF}
Running on a positive $\nn$-strand braid word of
length~$\ell$, the algorithm of Table~\ref{T:Algo} returns the
$\Phi$-normal word that is equivalent to~$\ww$
in~$O(\ell^2 \nn\log\nn)$ steps; in the meanwhile, it also
determines the address list of~$\DDf\nn{\CL\ww}$.
\end{prop}

\begin{proof}
As for~$(ii)$, we recall from~\cite[Chapter 9]{Eps} that
there exists a division algorithm running in
time~$O(\ell\nn\log\nn)$.
\end{proof}

We refer to Table~\ref{T:Algo} for the algorithm 
determining the $\Phi$-normal form, and to
Table~\ref{T:XNF} for the details of the computation
for~$\D_4^2$. Note that, apart from the fact
that letters come gathered in blocks in the former, the only
difference between the unbracketing of the $\BBB\nn$-decomposition and the
$\Phi$-normal form viewed as a sequence of letters is that the
$\BBB\nn$-decomposition always finishes with a power of~$\ss1$,
possibly~$\s_1^0$, \ie,~$1$: for instance, the $\Phi$-normal
form of~$\ss2$ is~$\ss2$, \ie, the length one
sequence~$(\ss2)$, while its $\BBB3$-decomposition is the
length two sequence~$(\ss2, 1)$.

\subsection{A linear ordering on~$\BB\nn$}
\label{S:FlipOrder}

As the monoid~$\BB2$ is isomorphic to~$\Nat$, it is
equipped with a natural linear ordering. Now, as the
$\nn$-splitting associates with every braid of~$\BB\nn$ a
distinguished finite sequence of braids, of~$\BB\nno$, we can recursively define a linear
ordering of~$\BB\nn$.

\begin{defi}
\label{D:FlipOrder}
For $\nn\ge2$, we define the relation~$\lf_\nn$
on~$\BB\nn$ as follows:

$(i)$ For $\xx,\yy$ in~$\BB2$, we say that $\xx\lf_2\yy$ holds
for $\xx=\s_1^\pp$ and $\yy=\s_1^\qq$ with $\pp<\qq$;

$(ii)$ For $\xx,\yy$ in~$\BB\nn$ with $\nn\ge3$, we say that
$\xx\lf_\nn\yy$ holds if, letting $(\xx_\pp, \Ldots, \xx_1)$ and $(\yy_\qq, \Ldots,
\yy_1)$ be the $\nn$-splittings of~$\xx$
and~$\yy$, we have either $\pp<\qq$, or 
$\pp=\qq$ and for some~$\rr \le \pp$ we have
$\xx_{\rr'}=\yy_{\rr'}$ for $\pp \ge \rr' > \rr$ and
$\xx_\rr \lf_\nno \yy_\rr$.
\end{defi}

Thus, $\lf_\nn$ is a sort of lexicographic extension of the
natural order on~$\BB2$, \ie, on~$\Nat$, via  splittings.
The extension is not exactly lexicographic: before comparing
componentwise, we first compare the lengths of the
sequences, \ie, the $\nn$-breadths of the considered braids,
a comparison method called $\ShortLex$ in~\cite{Eps}.

\begin{prop}
\label{P:ShortLex}
$(i)$ For~$\nn\ge2$, the relation~$\lf_\nn$ is a linear
ordering of~$\BB\nn$, which is a well-ordering. For
each braid~$\xx$, the immediate $\lf_\nn$-successor of~$\xx$
is~$\xx\ss1$.

$(ii)$ For $\nn\ge3$, the order~$\lf_\nn$ extends the
order~$\lf_\nno$, and $\BB\nno$ is the initial segment
of~$\BB\nn$ determined by~$\ss\nno$, \ie, we
have $\BB\nno = \{\xx\in\BB\nn \mid \xx\lf_\nn \ss\nno\}$. 
\end{prop}

\begin{proof}
$(i)$ The relation~$\lf_2$ is a linear ordering of~$\BB2$. Then, $\lf_\nn$ being a
linear ordering of~$\BB\nn$ follows from $\lf_\nno$
being a linear ordering of~$\BB\nno$ and
the $\nn$-splitting being unique. That  $\lf_\nn$ is a
well-order results from a similar induction, owing to the
standard result that the
$\ShortLex$-extension of a well-ordering is a well-ordering. Finally,
if the $\nn$-splitting of~$\xx$ is $(\xx_\pp, \Ldots, \xx_1)$,
the $\nn$-splitting of~$\xx\ss1$ is $(\xx_\pp, \Ldots,
\xx_1\ss1)$, making it clear that
$\xx\ss1$ is the immediate successor of~$\xx$.

$(ii)$ For $\xx, \yy$ in~$\BB\nno$, the
$\nn$-splittings of~$\xx$ and~$\yy$ are
the length one sequences~$(\xx)$ and~$(\yy)$, so, by
definition,
$\xx\lf_\nn\yy$ is equivalent to $\xx\lf_\nno\yy$. On the
other hand, the $\nn$-splitting of~$\ss\nno$ is $(\ss1,
1)$, so $\xx\lf_\nn\ss\nno$ holds for each~$\xx$
in~$\BB\nno$. Conversely, assume
$\xx\in\BB\nn$ and $\xx\lf_\nn\ss\nno$. By construction,
if $(\xx_2,\xx_1)$ is a $\nn$-splitting,
$\xx_2$ is not~$1$, hence, by~$(i)$, we have
$\xx_2\gfe_\nn\ss1$. So, if
$\xx\lf_\nn\ss\nno$ holds, the only possibility is that the
$\nn$-breadth of~$\xx$ is~$1$, \ie, that $\xx$ belongs
to~$\BB\nno$.
\end{proof}

Owing to Proposition~\ref{P:ShortLex}$(ii)$, we shall skip the
index~$\nn$ and write~$\lf$ for~$\lf_\nn$.

\begin{exam}
\label{X:Compar}
The $3$-splittings of $\ss1$ and $\ss2$
respectively are $(\ss1)$ and $(\ss1, 1)$, \ie, their
respective $3$-breadths are~$1$ and~$2$. Hence we have
$\ss1 \lf \ss2$.

Similarly, the $3$-splittings of $\D_3$ and
$\s_1^2\s_2^2$ are $(\ss1, \ss1, \ss1)$ and $(\s_1^2, \s_1^2, 1)$. The $3$-breadth is~$3$ in both
cases, and we compare lexicographically. The first entries
are $\ss1$ and $\s_1^2$. The former is smaller, hence $\D_3 \lf \s_1^2\s_2^2$ holds.
\end{exam}

The order~$\lf$ has been introduced above by means of the splitting. It can be introduced equivalently by appealing
to the exponent sequence of the $\BBB\nn$-decomposition and to the following ordering of iterated sequences of
integers.

\begin{defi}
If $\uuu, \vvv$ are $\nn$-sequences of natural numbers, we say that $\uuu$ is $\ShortLex$-smaller than~$\vvv$,
denoted $\uuu \lSL \vvv$, if we have $\nn = 0$ and $\uuu$ is smaller than~$\vvv$ with respect to the standard
order on~$\Nat$, or $\nn \ge 1$ and either $\uuu$---viewed as a sequence of $(\nno)$-sequences---is shorter
than~$\vvv$, or they have equal length and $\uuu$ is lexicographically smaller than~$\vvv$, \ie, writing $\uuu
= (\uuu_\pp, \Ldots, \uuu_1)$ and $\vvv = (\vvv_\pp, \Ldots, \vvv_1)$, there exists $\rr \le \pp$ such that we
have $\uuu_{\rr'} = \vvv_{\rr'}$ for $\pp \ge \rr' > \rr$ and $\uuu_\rr \lSL \vvv_\rr$.
\end{defi}

\begin{lemm}
\label{L:OrderCo1}
For~$\xx,\yy$ in~$\BB\nn$, we have
\begin{equation}
\label{E:FlipOrder}
\xx\lf\yy
\quad\Longleftrightarrow\quad
\DDfe\nn\xx \lSL \DDfe\nn\yy.
\end{equation}
\end{lemm}

\begin{proof}
As the relations involved in both sides of ~\eqref{E:FlipOrder} are linear orderings, it is enough to prove one
implication. We shall prove using  induction on~$\nn\ge2$ that $\xx\lf\yy$ implies $\DDfe\nn\xx \lSL
\DDfe\nn\yy$. The result is obvious for $\nn=2$. Assume $\nn\ge3$ and $\xx \lf \yy$ in~$\BB\nn$. Let
$(\xx_\pp, \Ldots, \xx_1)$ and $(\yy_\qq, \Ldots, \yy_1)$ be the
$\nn$-splittings  of~$\xx$ and~$\yy$.
By~\eqref{E:BiFlip21}, we have
\begin{gather}
\label{E:Expon1}
\DDfe\nn\xx=(\DDfe\nno{\xx_\pp}, \Ldots,
\DDfe\nno{\xx_1}), \quad
\DDfe\nn\yy=(\DDfe\nno{\yy_\qq}, \Ldots,
\DDfe\nno{\yy_1})
\end{gather}
---as the names of the generators are forgotten, the flips do not appear in exponent
sequences. According to the definition of~$\lf$, two cases are possible. If $\pp < \qq$ holds, then the left sequence
in~\eqref{E:Expon1} is shorter than the right sequence, so $\DDfe\nn\xx \lSL\nobreak
\DDfe\nn\yy$ holds. Otherwise, for some~$\rr \le \pp$, we must have $\xx_{\rr'} = \yy_{\rr'}$ for $\pp
\ge \rr' > \rr$ and $\xx_\rr \lf \yy_\rr$. We deduce $\DDfe\nno{\xx_{\rr'}} = \DDfe\nno{\yy_{\rr'}}$ for $\pp
\ge \rr' > \rr$ and, using the induction hypothesis, $\DDfe\nno{\xx_\rr} \lSL \DDfe\nno{\yy_\rr}$. Here again,
we find $\DDfe\nn\xx \lSL \DDfe\nn\yy$.
\end{proof}

For instance, we saw in Example~\ref{X:Compar} that
$\D_3 \lf \s_1^2\s_2^2$ holds. Another way to
see it is to compare $\DDfe3{D_3}$ and $\DDfe3{\s_1^2\s_2^2}$ with respect to~$\lSL$. The
respective values are $(1,1,1)$ and $(2,2,0)$: the former is $\lSL$-smaller.

\subsection{The braids~$\DDhat\nn\pp$}

Few properties of the order~$\lf$ are visible
directly. Typically, whether $\xx\lf\yy$ implies
$\zz\xx\lf\zz\yy$ is unclear because we do not know much
about the $\nn$-splittings of~$\zz\xx$ and~$\zz\yy$ as
compared with those of~$\xx$ and~$\yy$. We shall come
back on the question in Section~\ref{S:Burckel}.
For the moment, we conclude this section with a technical result about~$\lf$, namely we determine the least upper
bound of the braids of~$\BB\nn$ whose $\nn$-breadth is at most~$\pp$. 

\begin{figure}[htb]
\begin{picture}(97,16)
\put(0,-2){\includegraphics{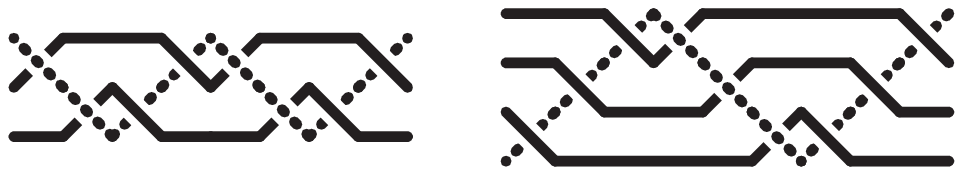}}
\end{picture}
\caption{\sf The braids~$\DDhat34$ (left) and $\DDhat43$ (right):
starting from the right, the upper strand of~$\DDhat\nn\dd$ forms $\dd$~half-twists around all other strands.}
\label{F:DeltaHat}
\end{figure}

\begin{nota}
(See Figure~\ref{F:DeltaHat}) 
For $\nn\ge2$ and $\dd\ge 1$, we set
\begin{equation}
\label{E:dddd}
\ddd\nn=\ss\nno...\ss1
\text{\  and \ }
\DDhat\nn\dd=\ff\nn^{\dd+1}(\ddd\nn) \cdot  ...
\cdot \ff\nn^2(\ddd\nn) \cdot \ff\nn(\ddd\nn).
\end{equation}
\end{nota}

In other words, $\DDhat\nn\dd$ is the
length~$\dd(\nn-1)$ zigzag $...\ss\nno ... \ss1\ss1 ...
\ss\nno$ with $\dd-1$ alternations, finishing
with~$\ss\nno$. For instance, 
$\DDhat42$ is the braid~$\ss3\ss2\s_1^2\ss2\ss3$.

\begin{lemm}
\label{L:dddd}
$(i)$ For $\nn\ge2$ and $\dd\ge 1$, we have
\begin{equation}
\label{E:Delta}
\D_\nn^\dd = \DDhat\nn\dd\,\D_\nno^\dd.
\end{equation}

$(ii)$ For~$\nn\ge2$,  $\dd\ge1$, and $\xx \in \BB\nno$,
the $\nn$-splitting  of~$\DDhat\nn\dd \, \xx$ is
\begin{equation}
\label{E:dddd3}
(\ss1 \ ,\ 
\underbrace{\mbox{$\d_\nno\ss1 \ , \ ... \ , \
\d_\nno\ss1$}}_{\dd-1\text{\ times}}\ ,\  \d_\nno \ ,\
\xx).
\end{equation}
This holds in particular for~$\DDhat\nn\dd$
with~$\xx=1$, and for~$\D_\nn^\dd$ with
$\xx=\D_\nno^\dd$.
\end{lemm}

\begin{proof}
$(i)$ Among the many equivalent inductive definitions
of~$\D_\nn$, we choose the recursive definition $\D_1= 1$ and $\D_\nn= \ss1...\ss\nno\D_\nno$, 
\ie, $\D_\nn=\DDhat\nn1\D_\nno$, for $\nn \ge 2$. Then \eqref{E:Delta} holds
for~$\dd=1$. For~$\dd\ge2$, we use induction:
\begin{multline*}
\D_\nn^\dd
=\D_\nn\D_\nn^{\dd-1}
=\D_\nn\DDhat\nn{\dd-1}\D_\nno^{\dd-1}
=\ff\nn(\DDhat\nn{\dd-1})\D_\nn\D_\nno^{\dd-1}\\
=\ff\nn(\DDhat\nn{\dd-1})\DDhat\nn1\D_\nno\D_\nno^{\dd-1}
=\DDhat\nn\dd\D_\nno^\dd.
\end{multline*}

$(ii)$ When we evaluate the sequence of~\eqref{E:dddd3}
by flipping each other entry, we obtain~$\DDhat\nn\dd\,\xx$. On the other
hand, each entry in~\eqref{E:dddd3} except possibly the last one is
right divisible by~$\ss1$, and by no  other~$\ss\ii$. Hence, by
Proposition~\ref{P:FSplitting}, the considered sequence is the
$\nn$-splitting of the braid it represents.
\end{proof}

In particular, the $3$-splitting of~$\D_3^\dd$ is 
$(\ss1 ,\s_1^2 , \ ... \ , \s_1^2 ,\ss1 ,\s_1^\dd)$, $\dd-1$~times~$\s_1^2$, 
which is $(\ss1, \ss1, \ss1)$ for $\dd=1$, corresponding to
$\D_3=\ss1\ss2\ss1$, and $(\ss1, \s_1^2, \ss1, \s_1^2)$
for $\dd=2$, corresponding to
$\D_3^2=\ss2\s_1^2\ss2\s_1^2$. 

We shall see that $\DDhat\nn\ppo$ is the least upper bound for the
braids of~$\BB\nn$ whose $\nn$-breadth at most~$\pp$. To prove this, we shall show that the $\nn$-splitting
of~$\DDhat\nn\ppo$ is minimal among all $\nn$-splittings of length~$\pp+1$. Therefore, we first investigate the
constraints satisfied by
$\nn$-splittings.

\begin{lemm}
\label{L:Small}
For $\nn\ge2$, the braids in~$\BB\nn$ that satisfy
$\xx\lf\ddd\nn$ are of those of the form
$\ss\nno...\ss\mm\yy$ with $\nn\ge\mm\ge2$ and
$\yy\in\BB{\mm-1}$.
\end{lemm}

\begin{proof}
We use induction on $\nn\ge2$. For $\nn=2$, we
have $\ddd\nn=\ss1$, and the result is true, as $\xx\lf\ss1$
implies $\xx=1$, and $1$ is the only element of~$\BB1$.
Assume $\nn\ge3$, and $\xx\lf\ddd\nn$. The $\nn$-splitting 
of~$\ddd\nn$ is $(\ss1,\ddd\nno)$. By definition, two
cases are possible: either the $\nn$-breadth of~$\xx$
is~$1$, which means that $\xx$ lies in~$\BB\nno$, or the
$\nn$-breadth of~$\xx$ is~$2$ and, letting $(\xx_2,
\xx_1)$ be its $\nn$-splitting, we have
either $\xx_2\lf\ss1$, which is impossible, or
$\xx_2=\ss1$ and $\xx_1\lf\ddd\nno$. In the latter
case, by induction hypothesis, there exist~$\mm$ with
$\nn-1\ge\mm\ge2$ and $\yy$ in~$\BB{\mm-1}$ such that
$\xx_1= \ss{\nn-2}...\ss\mm\yy$ holds, and, then, we
find $\xx=\ss\nno\ss{\nn-2}...\ss\mm\yy$.
\end{proof}

\begin{prop}
\label{P:Constraints}
Assume that $(\xx_\pp, \Ldots, \xx_1)$ is the $\nn$-splitting of some braid in~$\BB\nn$. Then the following
constraints are satisfied:
\begin{equation}
\label{E:Constraints}
\xx_\pp \gfe \ss1, \quad
\xx_\rr \gfe \ddd\nno \ss1
\text{\ for $\pp > \rr \ge 3$}, \quad
\xx_2\gfe \ddd\nno
\text{\ if $\pp \ge 3$ holds}.
\end{equation}
\end{prop}

\begin{proof}
First, we have $\xx_\pp \not= 1$ by hypothesis, hence $\xx_\pp \gfe \ss1$ by
Proposition~\ref{P:ShortLex}$(i)$.

Then, $\xx_\rr$ is right divisible by~$\ss1$ for $\rr\ge2$.
Indeed, by Proposition~\ref{P:AltDec}, we have $\xx_\rr \not=
1$, hence $\xx_\rr \multe \ss\ii$ for some~$\ii$. Now
$\xx_\rr \multe \ss\ii$ implies
$\ff\nn^{\pp-\rr}(\xx_\pp) \cdot  ... \cdot \xx_\rr \multe \ss\ii$,
and $\ii\ge2$ would contradict the $\nn$-splitting
condition of Proposition~\ref{P:FSplitting} at position~$\rr$. 

Assume $\pp > \rr\ge3$, and $\xx_\rr \lf \ddd\nno
\ss1$. Write $\xx_\rr = \yy_\rr \ss1$. By
Proposition~\ref{P:ShortLex}$(i)$,
$\xx_\rr$ is the immediate successor of~$\yy_\rr$, so  $\xx_\rr
\lf \ddd\nno \ss1$ implies $\yy_\rr \lf \ddd\nno$. By
Lemma~\ref{L:Small}, we have $\yy_\rr = \ss{\nn-2} ...
\ss\mm\yy$ with~$\yy$ in~$\BB{\mm-1}$ and
$\nn-1\ge\mm\ge\nobreak2$. The condition $\xx_{\rr+1}
\not= 1$ implies $\ff\nn(\xx_{\rr+1}) \multe
\ss\nno$, hence $\ff\nn(\xx_{\rr+1}) \cdot \xx_\rr \multe
\ss\nno \, ...\, \ss\mm \yy \ss1$. Assume first $\nn > \mm \ge
3$. Then $\ss{\mm}$ commutes with~$\yy_\rr$ and
with~$\ss1$, and we obtain $\ff\nn(\xx_{\rr+1}) \cdot \xx_\rr
\multe \ss\mm$, which contradicts the $\nn$-splitting
condition at position~$\rr$. Assume now $\mm = 2$, hence
$\yy=1$. Then we have
$\ff\nn(\xx_{\rr+1}) \cdot \xx_\rr \multe \ss\nno \, ... \, \ss1$,
hence 
$$\xx_{\rr+1} \cdot \ff\nn(\xx_\rr) \cdot \xx_{\rr-1} \multe
\ss1\, ... \, \ss\nno \xx_{\rr-1}.$$
Now, for~$\ii
\le \nn-2$, we have $\ss1\, ... \, \ss\nno \ss\ii = \ss{\ii+1}
\ss1\, ... \, \ss\nno$, so there exists~$\xx'$ for which
$\ss1\, ... \, \ss\nno \xx_{\rr-1} = \xx' \ss1\, ... \, \ss\nno$  holds. We deduce $\xx_{\rr+1} \cdot
\ff\nn(\xx_\rr)
\cdot \xx_{\rr-1} \multe \ss\nno$, contradicting the
$\nn$-splitting condition at position~$\rr-1$.

Assume finally $\xx_2 \lf \ddd\nno$. By
Lemma~\ref{L:Small}, we can write $\xx_2 = \ss{\nn-2}\, ... \, \ss\mm\xx$ with~$\xx$ in~$\BB{\mm-1}$ and
$\nn-1\ge\mm\ge2$. As above, we deduce $\ff\nn(\xx_3)
\cdot \xx_2 \multe \ss\nno\, ... \, \ss\mm \xx$. If $\nn > \mm
\ge 3$ holds, $\ss\mm$ commutes with~$\xx$, and we obtain 
$\ff\nn(\xx_3) \cdot \xx_2 \multe\nobreak \ss\mm$, which contradicts
the $\nn$-splitting condition at position~$2$. For $\mm
= 2$, hence $\xx=1$, we obtain $\ff\nn(\xx_3) \cdot \xx_2
\multe \ss2$ directly, and the same contradiction.
\end{proof}

\begin{prop}
\label{P:UpperBound}
For $\pp\ge1$, the braid~$\DDhat\nn\ppo$ is the $\lf$-least upper bound of the elements of~$\BB\nn$ whose
$\nn$-breadth is at most~$\pp$.
\end{prop}

\begin{proof}
By Lemma~\ref{L:dddd}$(ii)$, $\DDhat\nn\ppo$ has
$\nn$-breadth~$\pp+1$, hence $\xx\lf\DDhat\nn\ppo$ holds
for~$\xx$ with $\nn$-breadth at most~$\pp$.
Conversely, assume that the $\nn$-breadth of~$\xx$ is at
least~$\pp+1$. If it is~$\pp+2$ or more, then
$\xx\gf\DDhat\nn\ppo$ holds by definition of~$\lf$.
Otherwise, let $(\xx_{\pp+1}, \Ldots, \xx_1)$ be the
$\nn$-splitting of~$\xx$. Proposition~\ref{P:Constraints} says
that the sequence $(\xx_{\pp+1}, \Ldots, \xx_1)$ is at least
$(\ss1, \ddd\nno\ss1, \Ldots, \ddd\nno\ss1, \ddd\nno,
1)$, which is the $\nn$-splitting of~$\DDhat\nn\ppo$.
Hence we have $\xx\gfe\DDhat\nn\ppo$.
\end{proof}

\section{Connection with the braid order}
\label{S:Burckel}

Defining a unique normal representative is of little interest, unless the normal form has some specific
additional properties that make it useful. At the moment, the most interesting
property of the $\Phi$-normal form of braids
seems to be its connection with the so-called Dehornoy order.

\subsection{The braid order}
\label{S:Order}

We shall establish a simple connection between the
$\lf$-ordering of~$\BB\nn$, \ie, the ordering deduced from
the $\nn$-splitting, and the standard linear ordering of braids of~\cite{Dhr}. We recall the definition of the latter.
Considering~$\BB\nno$ as a submonoid of~$\BB\nn$, we
denote by~$\BB\infty$ the union of all~$\BB\nn$'s, and by~$B_\infty$
the group of fractions of~$\BB\infty$, \ie, the braid group on
unboundedly many strands.

\begin{defi}
\label{D:Order}
For $\xx,\yy$ in~$B_\infty$, we say that $\xx\ls\yy$
holds if the braid~$\xx\inv\yy$ admits at least one
word representative in which the generator~$\ss\ii$ with
maximal index occurs positively only, \ie, $\ss\ii$ occurs
but $\sss\ii$ does not.
\end{defi}

\begin{thrm}
\label{T:Order}
$(i)$ \cite{Dfb}
The relation~$\ls$ is a linear ordering of~$B_\infty$ that is
compatible with multiplication on the left.

$(ii)$ \cite{Lve}
The restriction of~$\ls$ to~$\BB\infty$ is a well-ordering.

$(iii)$ \cite{Bur}
For each~$\nn\ge2$, the restriction of~$\ls$ to~$\BB\nn$,
which is the interval $(1, \ss\nn)$ of~$(\BB\infty, \ls)$, is a
well-ordering of type~$\om^{\om^{\nn-2}}$.
\end{thrm}

In the framework of~\cite{Dhr}, the ordering of Definition~\ref{D:Order} is called the upper version of the braid
order. In some sources, in particular the early ones, the lower variant is considered, namely the relation~$\lmin$
referring to the letter~$\ss\ii$ with \emph{minimal} index,
instead of maximal as above. Both relations are similar as $\xx\ls\yy$ is equivalent
to $\ff\nn(\xx)\lmin\ff\nn(\yy)$ for all~$\xx, \yy$
in~$B_\nn$. However, as first noted by S.\,Burckel
in~\cite{BuT}, the statements involving
the well-order property are more natural with~$\ls$.

\subsection{Adding brackets in a braid word}

In order to connect the braid orders~$\lf$ and~$\ls$,
we shall compare the $\Phi$-normal form of Section~\ref{S:Flip}
with some other normal form introduced by S.\,Burckel in his
remarkable work, and we first need to introduce some notions
from~\cite{Bur}. The original description of~\cite{Bur} is
formulated in terms of trees. However, the latter are equivalent to the
iterated sequences of Section~\ref{S:Iter}, and we can easily
describe the fragment of Burckel's construction needed here
in terms of iterated sequences. Here , we give a new
description that is more directly connected to our
approach. In terms of trees, this amounts to starting from the
top and the right, while Burckel's approach starts from the
bottom and the left. The equivalence of both descriptions is
established in Proposition~\ref{P:CoincidenceNF} below.

Our basic observation here is that a free monoid is locally right Garside: this is a trivial result, as the
right divisibility relation of a free monoid is simply the relation of being a suffix. Then, applying the decomposition
process of Sections~\ref{S:Alt} and~\ref{S:Iter} to a word~$\ww$ in a free monoid amounts to grouping the letters
of~$\ww$ into blocks, \ie, in adding brackets in~$\ww$. We shall consider the iterated covering of the free
braid word monoids that mimicks the covering~$\BBB\nn$ of Section~\ref{S:Flip}.

\begin{nota}
We denote by~$\WW\nn$ the free monoid consisting of all positive $\nn$-strand braid words, and by~$\WWW\nn$
the atomic iterated covering of~$\WW\nn$ based on the sequence~$\GGen_\nn$---the same as in the case
of~$\BBB\nn$.
\end{nota}

We shall now use the $\WWW\nn$-decomposition of a word in~$\WW\nn$.
As in Section~\ref{S:Flip}, it is convenient to take advantage of the recursive definition of
the covering~$\WWW\nn$, and to introduce the counterpart of the $\nn$-splitting. 

\begin{defi}
For $\nn \ge 3$ and $\ww$ in~$\WW\nn$, the \emph{$\nn$-splitting} of~$\ww$ is defined to be the unique
sequence $(\ww_\pp, \Ldots, \ww_1)$ of words in~$\WW\nno$ such that $(\ff\nn^\ppo(\ww_\pp), \Ldots,
\ff\nn(\ww_2), \ww_1)$ is the $(\ff\nn(\WW\nno), \WW\nno)$-decomposition of~$\ww$.
\end{defi}

As being a right divisor in a free monoid is equivalent to being a suffix, Proposition~\ref{P:AltDec} implies that
$(\ww_\pp, \Ldots, \ww_1)$ is the $\nn$-splitting of~$\ww$ if and only if, for each~$\rr$, the word~$\ww_\rr$ is
the longest suffix of~$\ff\nn^{\pp-\rr}(\ww_\pp) \cdot ... \cdot \ww_\pp$ that lies in~$\WW\nno$. 

\begin{exam}
\label{X:WordSplitting}
Let $\ww$ be the $4$-strand braid word $\ss3\ss2\s_1^2\ss2\ss3
\ss2\s_1^2\ss2 \s_1^2$. The longest suffix of~$\ww$ that lies in~$\WW3$ is $\ss2\s_1^2\ss2 \s_1^2$, and the 
remaining prefix is $\ss3\ss2\s_1^2\ss2\ss3$, \ie,
$\ff4(\ww^{(1)})$ with $\ww^{(1)} =
\ss1\ss2\s_3^2\ss2\ss1$. The longest suffix of~$\ww^{(1)}$
that does not contain~$\ss3$ is $\ss2\ss1$, with remaining
prefix~$\ss1\ss2\s_3^2$, \ie, $\ff4(\ww^{(2)})$ with
$\ww^{(2)} = \ss3\ss2\s_1^2$. The longest suffix
of~$\ww^{(2)}$ that does not contain~$\ss3$ is
$\ss2\s_1^2$, with remaining prefix~$\ss3$. So, by definition, the $4$-splitting of the word~$\ww$ is the sequence
of $3$-strand braid words $(\ss1,\, \ss2\s_1^2,\, \ss2\ss1,\, \ss2\s_1^2\ss2 \s_1^2)$. 
\end{exam}

Imitating for braid words the notation used for braids in Section~\ref{S:Flip}, we put:

\begin{nota}
For $\ww$ in~$\WW\nn$, we denote the $\WWW\nn$-decomposition of~$\ww$ by~$\TTT\nn\ww$, and its
exponent sequence by~$\TTTe\nn\ww$.
\end{nota}

By construction, the iterated sequence~$\TTT\nn\ww$ is a certain bracketing of~$\ww$. Before giving an example,
we note the following connection between the $\WWW\nn$-decomposition and the splitting.

\begin{lemm}
\label{L:DecompWord}
For $\nn \ge 3$ and  $\ww$ in~$\WW\nn$, we have
\begin{equation}
\label{E:DecompWord}
\TTT\nn\ww=(\ff\nn^\ppo(\TTT\nno{\ww_\pp}), \Ldots,
\ff\nn(\TTT\nno{\ww_2}), \TTT\nno{\ww_1}).
\end{equation}
where $(\ww_\pp, \Ldots, \ww_1)$ is the $\nn$-splitting of~$\ww$.
\end{lemm}

The proof is exactly similar to that of Lemma~\ref{L:DecompBraid}.

\begin{exam}
\label{X:Bracketing}
Let again $\ww$ be the $4$-strand braid word $\ss3\ss2\s_1^2\ss2\ss3
\ss2\s_1^2\ss2 \s_1^2$. We saw in Example~\ref{X:WordSplitting} that the $4$-splitting of~$\ww$ is
$(\ss1, \ss2\s_1^2, \ss2\ss1, \ss2\s_1^2\ss2 \s_1^2)$. Then, we can easily see that the
$3$-splitting of the word $\ss2\s_1^2\ss2 \s_1^2$ is $(\ss2, \s_1^2, \ss2, \s_1^2)$, etc.
Using~\eqref{E:DecompWord}, we conclude that the $\WWW4$-decomposition of~$\ww$ is the
$2$-sequence
\begin{equation}
\label{E:ExampleWord}
\TTT4\ww = ( (\ss3),  (\ss2,\s_1^2),  (\ss2, \ss3),  (\ss2, \s_1^2, \ss2, \s_1^2)  ).
\end{equation}
\end{exam}

The braid word~$\ww$ considered in Example~\ref{X:Bracketing} is the $\Phi$-normal form of~$\D_4^2$.   By
comparing~\eqref{E:XFlip} and~\eqref{E:ExampleWord}, we see that, up to identifying the word~$\s_\ii^\ee$ with
the braid it represents, the $\WWW4$-decomposition of the word~$\ww$ is the
$\BBB4$-decomposition of~$\D_4^2$. This phenomenon is general. 

\begin{lemm}
\label{L:BracketNormal}
If $\ww$ is a $\Phi$-normal $\nn$-strand braid word, we have
$\TTT\nn\ww =
\DDf\nn{\CL\ww}$.
\end{lemm}

\begin{proof}
We use induction on~$\nn$. For $\nn=2$, the result
is obvious. Otherwise, let $(\xx_\pp, \Ldots,
\xx_1)$ be the $\nn$-splitting of~$\CL\ww$, and, for
each~$\rr$, let $\ww_\rr$ be the $\Phi$-normal form
of~$\xx_\rr$. By construction, each word~$\ww_\rr$ with
$\rr \ge 2$ finishes with~$\ss1$, so $(\ww_\pp, \Ldots, \ww_1)$
is the $\nn$-splitting of~$\ww$. The induction hypothesis
implies 
$\TTT\nno{{\ww_\rr}} = \DDf\nno{\xx_\rr}$ for
each~$\rr$. Applying~\eqref{E:DecompWord}, we deduce
\begin{align*}
\TTT\nn\ww 
&=  (\ff\nn^\ppo(\TTT\nno{{\ww_\pp}}), \Ldots, 
\ff\nn(\TTT\nno{{\ww_2}}), 
\TTT\nno{{\ww_1}})\\
&=(\ff\nn^\ppo(\DDf\nno{\xx_\pp}), \Ldots, 
\ff\nn(\DDf\nno{\xx_2}),  
\DDf\nno{\xx_1}).
\end{align*}
By~\eqref{E:BiFlip21}, the latter sequence
is~$\DDf\nn{\CL\ww}$.
\end{proof}

At this point, we can easily establish the connection between our current notion of $\WWW\nn$-decomposition and
Burckel's notion of ``the tree of a braid word''.

\begin{lemm}
\label{L:Bracketing}
Assume $\nn\ge3$ and  $\ww \in \WW\nn$ with $\TTT\nn\ww = (\www_\pp, \Ldots,
\www_1)$. Then, for $1 \le \ii \le \nno$, and assuming $\www_\pp =\nobreak \TTT\nno{\ww_\pp}$, we have 
$$\TTT\nn{{\ss\ii}\ww} = 
\begin{cases}
((...(\ss\ii )...),  \www_\pp, \Ldots, \www_1)
&\text{for $\pp$ even and $\ii = 1$},\\
&\text{\quad and for $\pp$ odd and $\ii=\nno$},\\
{(\TTT\nno{\ff\nn^\ppo(\ss\ii) \ww_\pp}, 
\www_\ppo, \Ldots, \www_1)}
&\text{otherwise}.
\end{cases}
$$
\end{lemm}

\begin{proof}
Let $(\ww_\pp, \Ldots, \ww_1)$ be the $\nn$-splitting of~$\ww$. Then the $\nn$-splitting of~$\ss\ii \ww$ is
$(\ss1, \ww_\pp, \Ldots, \ww_1)$ for $\pp$ even and $\ii = 1$, and for $\pp$ odd and $\ii = \nno$,
and it is $(\ff\nn^\ppo(\ss\ii) \, \ww_\pp,  \ww_\ppo, \Ldots,
\ww_1)$ otherwise. Indeed, the point is whether the additional letter~$\ss\ii$ can be
incorporated in the same entry as~$\ww_\pp$. Taking the
flips into account, this depends on whether $\ff\nn^\ppo(\ss\ii)$
is~$\ss\nno$ or not. The value of~$\TTT\nn{{\ss\ii}\ww}$ directly follows. 
\end{proof}

As the rule of Lemma~\ref{L:Bracketing} directly mimicks the inductive construction of the tree associated
with~$\ww$ in the sense of~\cite{Bur}, we deduce:

\begin{prop}
\label{P:CoincidenceNF}
For each positive $\nn$-strand braid word~$\ww$, the tree
associated with~$\TTT\nn\ww$ coincides with the tree of~$\ww$ as defined in~\cite{Bur}. 
\end{prop}

Before going to Burckel's results, let us observe that the braid ordering~$\lf$ of  of Definition~\ref{D:FlipOrder} 
admits a simple characterization in terms of $\Phi$-normal words. 

\begin{prop}
\label{P:Connection1}
For all~$\xx,\yy$ in~$\BB\nn$, we have
\begin{equation}
\label{E:Connection1}
\xx\lf\yy \quad\Longleftrightarrow\quad
\TTTe\nn\uu \lSL\TTTe\nn\vv,
\end{equation}
where $\uu$ and~$\vv$ are the $\Phi$-normal representatives
of~$\xx$ and~$\yy$.
\end{prop}

\begin{proof}
By Lemma~\ref{L:BracketNormal}, we have
$\DDf\nn\xx=\TTT\nn\uu$ and 
$\DDf\nn\yy=\TTT\nn\vv$, so $\xx\lf\yy$, which is equivalent to
$\DDfe\nn\xx\lSL\DDfe\nn\yy$ by Lemma~\ref{L:OrderCo1}, is also equivalent to
$\TTTe\nn\uu \lSL\TTTe\nn\vv$.
\end{proof}

\begin{rema}
For $\ww$ in~$\WW\nn$, define a $\WWW\nn$-bracketing of~$\ww$ to be any $(\nnt)$-sequence~$\www$ such
that the unbracketing of~$\www$ is~$\ww$ and, for each address~$\h$ of length~$\nn\!-\!2$, the
entry~$\www_\h$ belongs to~$\WW{\cl\h}$ when it exists. So a $\WWW\nn$-bracketing of~$\ww$ is any way of
adding brackets in~$\ww$ so that the resulting iterated sequence has its entries correctly dispatched with respect to
the skeleton of the iterated covering~$\WWW\nn$. By construction, $\TTT\nn\ww$ is always a
$\WWW\nn$-bracketing of~$\ww$, but it is not the only one. For instance, both $(\s_2^2, \ss1)$ and $(\ss2, \e,
\ss2, \e, \e, \ss1)$ are
$\WWW3$-bracketings of the word~$\s_2^2\ss1$. Then it is easy to check that $\TTT\nn\ww$ is, among all
$\WWW\nn$-bracketings of~$\ww$, the one that has the $\lSL$-smallest exponent sequence. 
\end{rema}

\subsection{The Burckel normal form}

We now appeal to Burckel's results in~\cite{Bur} to state a connection between the braid ordering~$\ls$ and the
$\lSL$-ordering of the exponent sequences.  

\begin{defi}
A positive $\nn$-strand braid word~$\ww$ is said to
be \emph{Burckel normal} if the exponent
sequence~$\TTTe\nn\ww$ is
$\lSL$-minimal among all~$\TTTe\nn{\ww'}$
with~$\ww'\equiv\ww$.
\end{defi}

\begin{exam}
Let us consider the two positive $3$-strand braid words that represent~$\D_3$, namely $\ss1\ss2\ss1$ and
$\ss2\ss1\ss2$. Then we find $\TTT3{\ss1\ss2\ss1} = (\ss1, \ss2, \ss1)$, and $\TTT3{\ss2\ss1\ss1} = (\ss2, \ss1,
\ss2, \e)$---here we use the empty word $\e$ to emphasize that we consider words. So we have $\TTTe3{\ss1 \ss2
\ss1} = (1, 1, 1)$, and $\TTTe3{\ss2 \ss1 \ss2} = (1, 1, 1, 0)$. As $(1, 1, 1)$ is shorter, hence $\lSL$-smaller, than
$(1, 1, 1, 0)$, we conclude that $\ss1 \ss2 \ss1$ is Burckel normal, while $\ss2\ss1\ss2$ is not.
\end{exam}

Burckel normal words are called \emph{irreducible}
in~\cite{Bur}. As the $\ShortLex$-ordering of $\nn$-sequences on~$\Nat$ is a
well-ordering, each nonempty set of $\nn$-sequences in~$\Nat$
contains a $\lSL$-least element. Therefore, each
positive braid admits a unique Burckel normal representative. 

\begin{thrm} [Burckel, \cite{Bur}]
\label{T:Burckel}
For~$\xx,\yy$ in~$\BB\nn$, we have
\begin{equation}
\label{E:Connection2}
\xx\ls\yy \quad\Longleftrightarrow\quad
 \TTTe\nn\uu \lSL \TTTe\nn\vv,
\end{equation}
where $\uu$ and~$\vv$ are the Burckel normal
representatives of~$\xx$ and~$\yy$.
\end{thrm}

Burckel's proof of Theorem~\ref{T:Burckel} is quite subtle for $\nn\ge4$ and requires
a transfinite induction. The point is to define a combinatorial operation called reduction so
that, if a braid word~$\ww$ is not Burckel normal, then its
reduct~$\ww'$ is equivalent to~$\ww$ and satisfies
$\TTTe\nn{\ww'} \lSL \TTTe\nn\ww$. 

In the sequel, we shall only use the following consequence of
Theorem~\ref{T:Burckel}.

\begin{coro}
\label{C:Burckel}
If $\uu$ and $\vv$ are the Burckel normal representatives
of~$\xx$ and~$\xx\ss\ii$, then $\TTTe\nn\uu \lSL
\TTTe\nn\vv$ holds.
\end{coro}

\begin{proof}
By definition, we have $\xx \ls \xx\ss\ii$, as the quotient $\xx\inv \xx\ss\ii$ has an expression,
namely~$\ss\ii$, in which the generator with highest index appears positively only.
\end{proof}

\subsection{Connecting the normal forms}
\label{S:Connection}

At this point, two distinguished word representatives have
been introduced for each positive braid, namely its
$\Phi$-normal form, and its Burckel normal form.  We shall now
prove that these \emph{a priori} unrelated normal representatives actually
coincide.

\begin{prop}
\label{P:Connection}
The Burckel normal form coincides with the $\Phi$-normal
form.
\end{prop}

\begin{figure}[htb]
\begin{picture}(106,25)
\put(0,3){\includegraphics{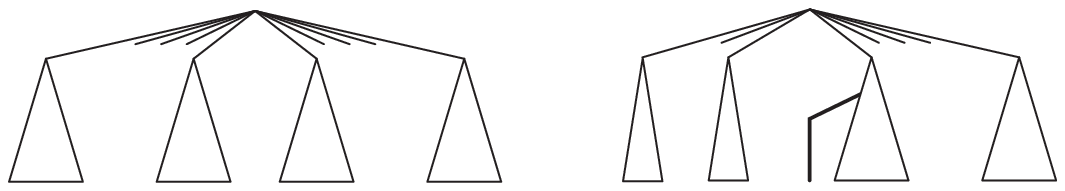}}
\put(20,23){$\TTT\nn\ww$}
\put(75,23){$\TTT\nn{\ww'}$}
\put(2,0){$\ww_\pp$}
\put(17,0){$\ww_\rr$}
\put(27,0){$\ww_\rro$}
\put(44,0){$\ww_1$}
\put(44,0){$\ww_1$}
\put(62,0){$\ww'_{\pp'}$}
\put(71,0){$\ww'_\rr$}
\put(80,0){$\ss\ii$}
\put(84,0){$\ww_\rro$}
\put(101,0){$\ww_1$}
\end{picture}
\caption{\sf Proof of Proposition~\ref{P:Connection}: if
$\ww$ is not $\Phi$-normal because some $\ss\ii$ with $\ii \ge
2$ right divides the braid associated with the $\pp\!-\nobreak\!\rr\!+\nobreak\!1$
left factors, then that $\ss\ii$ can be removed from
the left part and incorporated in the next factor; Corollary~\ref{C:Burckel} guarantees that the new left part is
smaller than the old one, so the new word~$\ww'$ is
equivalent to~$\ww$ but its exponent sequence is
smaller than that of~$\ww$.}
\label{F:Proof}
\end{figure}

\begin{proof}
As each braid admits a unique Burckel normal
representative and a unique $\Phi$-normal representative,
proving one implication is sufficient. Here we prove using
induction on~$\nn\ge2$ that an $\nn$-strand braid word that is not
$\Phi$-normal is not either Burckel normal. For $\nn=2$, every
word, namely every power of~$\ss1$, is normal in both
senses. Assume $\nn\ge3$, and assume that $\ww$ is a word in~$\WW\nn$ that is not $\Phi$-normal. We aim at
proving that $\ww$ is not Burckel normal. Owing to the definition of a Burckel normal word, it is enough to exhibit
a word~$\ww'$ that represents the same braid as~$\ww$ and is such that $\TTTe\nn{\ww'}$ is
$\ShortLex$-smaller than~$\TTTe\nn\ww$.

Let $(\ww_\pp, \Ldots, \ww_1)$ be the $\nn$-splitting of~$\ww$. By
Lemma~\ref{L:DecompWord}, the value of~$\TTTe\nn\ww$ is
\begin{equation}
\label{E:FirstSeq}
(\TTTe\nno{{\ww_\pp}}, \Ldots,
\TTTe\nno{{\ww_2}},
\TTTe\nno{{\ww_1}})
\end{equation} 
---as we consider exponent sequences, we can forget about flips. The hypothesis that $\ww$ is not $\Phi$-normal may
have two causes, namely that one of the words~$\ww_\rr$ is not $\Phi$-normal, or that all words~$\ww_\rr$ are $\Phi$-normal but $(\CL{\ww_\pp}, \Ldots, \CL{\ww_1})$ is not the $\nn$-splitting of the braid~$\CL\ww$.

Assume first that some word~$\ww_\rr$ is not $\Phi$-normal. By induction hypothesis, $\ww_\rr$ is not Burckel
normal either. Hence there exists a word~$\ww'_\rr$ equivalent to~$\ww_\rr$ satisfying 
$$\TTTe\nno{\ww'_\rr} \lSL \TTTe\nno{\ww_\rr}.$$
Let $\ww'$ be the word obtained from~$\ww$ by replacing the subword~$\ff\nn^\rro(\ww_\rr)$
with $\ff\nn^\rro(\ww'_\rr)$. Then
$\ww'$ is equivalent to~$\ww$, and, by construction, one has
$$\TTTe\nn{\ww'} \lSL \TTTe\nn{\ww},$$
hence $\ww$ cannot be Burckel normal.

Assume now that each word~$\ww_\rr$ is $\Phi$-normal and $(\CL{\ww_\pp}, \Ldots, \CL{\ww_1})$ is not the
$\nn$-splitting of~$\CL\ww$. Then there exists~$\rr$ such that the
braid represented by 
$$\vv = \ff\nn^{\pp-\rr}(\ww_\pp) \cdot ... \cdot \ff\nn(\ww_{\rr+1}) \cdot \ww_\rr$$
 is right divisible by some~$\ss\ii$ with $\ii \ge 2$. We shall show that the factor~$\ss\ii$ can be
removed from~$\ww_\rr$ and incorporated in the next factor~$\ww_\rro$,
so as to give rise to a new word~$\ww'$ equivalent to~$\ww$ and satisfying $\TTTe\nn{\ww'} \lSL
\TTTe\nn\ww$---see Figure~\ref{F:Proof}.

Indeed, let $\vv'$ be the Burckel normal form of~$\CL{\vv} \OVER{\s_\ii}$, and let $\ww'$ be the word
$\ff\nn^\rro(\vv') \cdot \ff\nn^\rrt(\ss{\nn-\ii} \ww_\rro) \cdot ... \cdot \ff\nn(\ww_2) \cdot \ww_1$. By
construction, $\ww'$ is equivalent to~$\ww$. The $\nn$-splitting of~$\vv$ is $(\ww_\pp, \Ldots, \ww_\rr)$.
Let $(\ww'_{\pp'}, \Ldots, \ww'_\rr)$ be that of~$\vv'$. Then the $\nn$-splitting of~$\ww'$ is
$(\ww'_{\pp'}, \Ldots, \ww'_\rr, \ss{\nn-\ii}\ww_\rro, \ww_\rrt, \Ldots, \ww_1)$, 
and so, by Lemma~\ref{L:DecompWord}, the value of~$\TTTe\nn{\ww'}$ is
\begin{multline}
\label{E:SecondSeq}
(\TTTe\nno{\ww'_{\pp'}}, \Ldots, \TTTe\nno{\ww'_\rr}, 
\TTTe\nno{\ss{\nn-\ii}\ww_\rro},  \TTTe\nno{\ww_\rrt}, \Ldots, \TTTe\nno{\ww_1}).
\end{multline}
Now---this is the point---Corollary~\ref{C:Burckel} implies $\TTTe\nn{\vv'}
\lSL \TTTe\nn\vv$, \ie, always by  Lemma~\ref{L:DecompWord}, 
$$(\TTTe\nno{\ww'_{\pp'}}, \Ldots, \TTTe\nno{\ww'_\rr}) 
\lSL (\TTTe\nno{\ww_\pp}, \Ldots, \TTTe\nno{\ww_\rr})$$
---hence in particular $\pp' \le \pp$. Adding~$\rr-1$ entries on the right of the above sequences does not change
their order, and we deduce that the sequence of~\eqref{E:SecondSeq} is $\ShortLex$-smaller than that
of~\eqref{E:FirstSeq}, \ie, $\TTTe\nn{\ww'}$ is $\ShortLex$-smaller than~$\TTTe\nn\ww$. This shows that $\ww$
is not Burckel normal.
\end{proof}

\begin{rema*}
It is natural to wonder whether  Proposition~\ref{P:Connection} extends to every dense atomic
covering~$\MMM$ of a locally right monoid~$\MM$, \ie, whether the $\MMM$-normal form of an element~$\xx$
of~$\MM$ is always the representative whose $\underline{\MMM}$-decomposition---defined in the
obvious way from the considered atoms---has the $\lSL$-minimal exponent sequence. This is
\emph{not} the case. Indeed, as in Remark~\ref{R:Shortest}, consider the $2$-covering~$\MMM$ of~$\BB4$
based on $((\ss3,\ss2), (\ss3,\ss1))$. Then, the $\MMM$-normal form of
$\ss3\ss2\ss1\ss2\s_3^2\ss2$ turns out to be the word $\ss1\s_2^2\ss1\ss3\ss2\ss1$. Now the
$\underline{\MMM}$-decompositions of the words
$\ss3\ss2\ss1\ss2\s_3^2\ss2$ and $\ss1\s_2^2\ss1\ss3\ss2\ss1$ respectively are
$$(  (\ss3,\ss2), (\ss1), (\ss2, \s_3^2, \ss2), (\e)
 )
\text{\quad and \quad}
( (\ss1), (\s_2^2), (\ss1), (\ss3, \ss2),g (\ss1)  ),$$ with lengths $4$ and $5$---this is the same
example as in Remark~\ref{R:Shortest}. The latter has a
$\lSL$-larger exponent sequence, so the $\MMM$-normal
form does not correspond to the smallest exponent sequence.  Technically, the point is that the counterpart to
Corollary~\ref{C:Burckel} fails: the breadth may decrease under right
multiplication. For instance, for $\yy =
\ss1\s_2^2\ss1\ss3\ss2\s_1^2$, the
$\MMM$-decomposition of~$\yy$ is $(\ss1 , \s_2^2 ,
\ss1\ss3 , \ss2 , \s_1^2)$, which has length~$5$, while that
of~$\yy\ss2$ is $(\ss3\s_2^2 , \ss1 , \ss2\s_3^2\ss2 , \ss1)$,
which has length~$4$. This shows that the covering~$\BBB\nn$ is quite specific.
\end{rema*}

\subsection{Applications}

Once we know that the $\Phi$-normal form and the Burckel
normal form coincide, each one inherits the properties
of the other, and we easily deduce several consequences, in
particular in terms of braid orderings.

\begin{prop}
\label{P:CompOrder}
For~$\xx,\yy\in\BB\infty$, the relations $\xx\ls\yy$
and $\xx\lf\yy$ are equivalent.
\end{prop}

\begin{proof}
Let $\uu$ and $\vv$ be the $\Phi$-normal representatives
of~$\xx$ and~$\yy$. By Proposition~\ref{P:Connection},
$\uu$ and~$\vv$ also are the Burckel normal
representatives of~$\xx$ and~$\yy$. The equivalences
$$\xx \ls \yy 
\quad\Longleftrightarrow\quad 
\TTTe\nn\uu \lSL \TTTe\nn\vv
\quad\Longleftrightarrow\quad 
\xx \lf \yy$$
then follow from Proposition~\ref{P:Connection1}
and Theorem~\ref{T:Burckel}.
\end{proof}

We deduce that the standard braid ordering~$\ls$ inherits the recursive definition of the
ordering~$\lf$, which is Theorem~A$(ii)$ in the introduction:

\begin{coro}
\label{C:RecOrder}
Let $\xx, \yy$ be positive $\nn$-strand braids. Let $(\xx_\pp, \Ldots, \xx_1)$ and $(\yy_\qq,\Ldots,\yy_1)$
be the $\nn$-splittings of~$\xx$ and~$\yy$. Then $\xx<\yy$ holds in~$\BB\nn$ if and only if we
have either $\pp<\qq$, or $\pp=\qq$ and there exists~$\rr \le \pp$ such that we have $\xx_{\rr'}=\yy_{\rr'}$ for
$\pp\ge\rr'>\rr$ and $\xx_\rr<\yy_\rr$ in~$\BB\nno$.
\end{coro}

In the other direction, we deduce that the ordering~$\lf$ satisfies the known properties of the ordering~$\ls$:

\begin{coro}
\label{C:Compat2}
The order~$\lf$ is compatible with multiplication on the left, and $\xx\lf\xx\ss\ii$ always holds.
\end{coro}

Futher consequences involve the algorithmic complexity. The following result deals with the braid order~$\ls$,
and it is Theorem~B in the introduction.

\begin{coro}
\label{C:Complexity} 
For each~$\nn$, the braid order~$\ls$ on~$B_\nn$ can be
decided in quadratic time: if $\ww$ is a (non necessarily
positive) $\nn$-strand braid word of length~$\ell$, then
whether $\CL\ww\gs1$ holds can be decided in
time~$O(\ell^2\nn^3\log\nn)$.
\end{coro}

\begin{proof}
We first observe
that, if $\uu,\vv$ are positive $\nn$-strand braid words of
length at most~$\ell$, then $\CL\uu\ls\CL\vv$ can be
decided in time~$O(\ell^2\nn\log\nn)$. Indeed, by
Proposition~\ref{P:BraidNF}$(ii)$, we can compute the
decompositions~$\DDf\nn{\CL\uu}$
and~$\DDf\nn{\CL\vv}$ within the indicated amount of time;
the extra cost of subsequently comparing the corresponding
exponent sequences with respect to the
$\ShortLex$-ordering is linear in~$\ell\nn$. 

If $\ww$ is an arbitrary $\nn$ strand braid word of
length~$\ell$, according to~\cite[Chapter~9]{Eps}, we can
find two positive braid words~$\uu, \vv$ of length
in~$O(\ell\nn^2)$ such that $\ww$ is equivalent
to~$\uu\inv\vv$ in time~$O(\ell^2\nn\log\nn)$. Then
$\CL\ww\gs1$ is equivalent to $\CL\uu\ls\CL\vv$, which,
by the above observation, can be decided in
time~$O(\ell^2\nn^5\log\nn)$. Actually, we can lower the
exponent of~$\nn$ to~$3$ because an upper bound for the
computation of the $\Phi$-normal form is~$O(\ell\ell_c\nn\log\nn)$, where
$\ell_c$ is the canonical length, \ie, the
number of divisors of~$\D_\nn$ involved in the right
greedy normal form. When we go from~$\ww$
to~$\uu\inv\vv$, the canonical lengths of~$\uu$
and~$\vv$ are bounded above by that of~$\ww$, leading
to~$O(\ell\ell_c\nn^3\log\nn)$ for the whole comparison.
\end{proof}

Finally, another application is that, for each~$\nn$, the Burckel normal form of a positive
$\nn$-strand braid word can be computed in quadratic time w.r.t.\,the
length of the initial word, which is clear from Proposition~\ref{P:BraidNF} and the fact that the Burckel normal form
coincides the $\Phi$-normal form. In the approach of~\cite{Bur}, the Burckel normal form comes as the final result of an
iterated reduction process whose convergence is guaranteed by the fact that an ordinal decreases, and no complexity
analysis has been published so far.

\section{Open questions and further work}

\subsection{The $\Phi$-normal form}

We have seen in Proposition~\ref{P:Constraints} that an arbitrary sequence of braids in~$\BB\nno$ need not be the
$\nn$-splitting of a braid in~$\BB\nn$. An obvious question is whether the constraints of 
Proposition~\ref{P:Constraints} are sufficient conditions.

\begin{ques}
Assume that $\xx_\pp, \Ldots, \xx_1$ are braids of~$\BB\nno$ that satisfy
\begin{equation}
\xx_\pp \gse \ss1, \quad
\xx_\rr \gse \ddd\nno \ss1
\text{\  for $\pp > \rr \ge 3$}, \quad
\xx_2\gse \ddd\nno
\text{\  if $\pp \ge 3$ holds}.
\end{equation}
Does there exist a braid in~$\BB\nn$ whose $\nn$-splitting is $(\xx_\pp, \Ldots, \xx_1)$?
\end{ques}

The only case where a (positive) answer is known is $\nn = 3$.

\begin{prop}
A sequence $(\s_1^{\ee_\pp}, \Ldots, \s_1^{\ee_1})$ is the $3$-splitting of a braid of~$\BB3$ if and only if the
numbers~$\ee_\rr$ satisfy the inequalities: 
\begin{equation}
\ee_\pp \ge 1, \quad
\ee_\rr \ge 2 \text{\ for $\pp > \rr \ge 3$}, \quad
\ee_2\ge 1 \text{\  if $\pp \ge 3$ holds}.
\end{equation}
\end{prop}

\begin{proof}
What remains to be shown is that, if at least one of the above conditions fails, then $(\s_1^{\ee_\pp}, \Ldots,
\s_1^{\ee_1})$ is not a $3$-splitting. Now, by Lemma~\ref{L:Dense}, no gap may exist in a $3$-splitting, so
$\ee_\rr = 0$ is impossible for $\pp > \rr \ge 2$.

On the other hand, assume $\ee_\rr = 1$ with $\pp > \rr \ge 3$. As we have $\s_1^{\ee_{\pp+1}} \ss2
\s_1^{\ee_\ppo} = \s_1^{\ee_{\pp+1}-1} \s_2^{\ee_\ppo} \ss1 \ss2$, the braid $\ff3^{\ppo-\rr}(\s_1^{\ee_\rr})
\cdot ... \cdot \s_1^{\ee_\ppo}$ is right divisible by~$\ss2$, contradicting the characteristic property of a
$3$-splitting.
\end{proof}

The result can be restated as

\begin{coro}
\label{C:Normal3}
Set $\emin_1 = 0$, $\emin_2 = 1$, and $\emin_\rr = 2$ for $\rr \ge 3$. Then a positive $3$-strand braid word
$\s_{\cl\pp}^{\ee_\pp} \, ... \, \s_2^{\ee_2} \, \s_1^{\ee_1}$ with $\ee_\pp \ge 1$ is $\Phi$-normal if and only if the
inequality $\ee_\rr \ge \emin_\rr$ is satisfied for all indices~$\rr$ except possibly~$\pp$.
\end{coro}

\begin{rema}
A priori, the $\Phi$-normal form of a positive braid is completely different from its right greedy normal form
of~\cite[Chapter~9]{Eps}. However, it was observed by J.~Mairesse (private communication) that, in the case of
$3$-strands, there is a rather simple connection: starting from the right greedy normal form of a positive $3$-strand
braid, we can obtain its $\Phi$-normal form by replacing the final factor~$\D_3^\ee$ with its $\Phi$-normal
form, and, depending on the parity of~$\ee$ and on the final letter in the next factor, possibly push some
factors~$\s_\ii^\dd$ through~$\D_3^\ee$---see~\cite{Dhq} for details.
\end{rema}

\subsection{The braid ordering}

The proof of Proposition~\ref{P:Connection} heavily
relies on Burckel's Theorem~\ref{T:Burckel}, a highly
non trivial combinatorial result in the general case.

\begin{ques}
\label{Q:Final}
Is there a direct proof for the following results?

$(i)$ The orders~$\lf$ and~$\ls$ coincide. 

$(ii)$ The order~$\lf$ is compatible with multiplication on
the left.

$(iii)$ The relation $\xx\lf\xx\ss\ii$ always holds.
\end{ques}

So far we have no general answer. We mention
below partial results toward a positive answer to
Question~\ref{Q:Final}$(i)$, namely proving that, for
all braids~$\xx,\yy$, the relation $\xx\lf\yy$ implies
$\xx\ls\yy$---as we are dealing with linear orders, one
implication is enough. Here we consider special values
for~$\yy$. By Proposition~\ref{P:ShortLex}$(ii)$, we already know that
$\xx\lf\ss\nno$ is equivalent to~$\xx\ls\ss\nno$, as
both are equivalent to
$\xx\in\BB\nno$. Here is another result of this kind.

\begin{prop}
\label{P:Compar}
For every~$\xx$ in~$\BB\nn$, the relation
$\xx\lf\DDhat\nn\dd$ implies $\xx<\DDhat\nn\dd$.
\end{prop}

\begin{proof}
Assume $\xx\lf\DDhat\nn\dd$. By
Proposition~\ref{P:UpperBound}, the $\nn$-breadth of~$\xx$ is
at most~$\dd+1$, and we can write $\xx =
\ff\nn^\dd(\xx_{\dd+1}) \cdot ... \cdot
\ff\nn(\xx_2) \cdot \xx_1$ for some $\xx_{\dd+1}, \Ldots, \xx_1$
in~$\BB\nno$. An easy computation using~\eqref{E:Delta} and the equalities $\ff\nn(\xx_\rr\inv)=\D_\nn
\xx_\rr\inv \D_\nn\inv$ gives
\begin{equation}
\label{E:Quotient}
\xx\inv\,\DDhat\nn\dd = 
\xx_1\inv \cdot \D_\nn\xx_2\inv \cdot \D_\nn\xx_3\inv ...
\cdot \D_\nn \xx_{\dd+1}\inv \cdot \D_\nno^{-\dd}.
\end{equation}
This leads to an expression of the quotient
$\xx\inv \DDhat\nn\dd$ in which the letter~$\ss\nno$
occurs $\dd$~times, while neither
$\sss\nno$ nor any letter~$\s_\jj^{\pm1}$ with
$\jj\ge\nn$ does. Indeed, each factor~$\D_\nn$ admits a
positive expression in which $\ss\nno$ occurs once, namely
the one arising from the decomposition $\D_\nn=\DDhat\nn1
\D_\nno$, while the negative factors~$\xx_\rr\inv$ and
$\D_\nno^{-\dd}$ belong to~$B_\nno$ and therefore can
be expressed using neither~$\ss{n-1}$ nor~$\sss\nno$.
Therefore $\xx\ls\DDhat\nn\dd$ holds.
\end{proof}

It is not hard to deduce that, for every~$\xx$ in~$\BB\nn$, the relation
$\xx\lf\D_\nn^\dd$ implies $\xx<\D_\nn^\dd$, as well as various similar compatibility results between~$\lf$
and~$\ls$. But, so far, we have no complete answer to Question~\ref{Q:Final}$(i)$ in the
general case. 

It is however easy to provide such an answer in the case $\nn = 3$. Indeed, in this special case, the exact form of $\Phi$-normal words is known, and a direct computation similar to that of Proposition~\ref{P:Compar} shows that, for $\xx,
\yy$ in~$\BB3$, the relation~$\xx \lf \yy$ implies that the braid
$\xx\inv \yy$ are an expression where $\ss2$ occurs but $\sss2$ does not, or an expression where $\ss1$ occurs
but none of~$\sss1, \ss2, \sss2$ does, \ie, that $\xx \ls \yy$ holds.

By~\cite{Lve} and~\cite{Bur}, we know that $(\BB3, \ls)$ is a well-ordering of order type~$\om^\om$. Hence
the position of every braid of~$\BB3$ is unambiguously specified by an ordinal number, called the \emph{rank}
of~$\xx$, namely the order type of the initial segment of~$(\BB3, \ls)$ determined by~$\xx$. Using the
formula for the $\Phi$-normal form given in Corollary~\ref{C:Normal3}, we deduce the following explicit
value for the rank of a $3$-strand braid.

\begin{prop}
\label{P:Rank}
The rank of the braid with $\Phi$-normal form $\s_{\cl\pp}^{\ee_\pp} \, ... \, \s_2^{\ee_2} \, \s_1^{\ee_1}$ in the
well-ordering $(\BB3, \ls)$ is the ordinal number
\begin{equation}
\label{E:Rank}
\om^\ppo \cdot {\ee_\pp} + \sum_{\rr = \ppo}^{\rr = 1} \om^\rro \cdot (\ee_\rr - \emin_\rr),
\end{equation}
where the (absolute) numbers~$\emin_\rr$ are those of Corollary~\ref{C:Normal3}.
\end{prop}

\begin{proof}
The point is to determine which $\Phi$-normal words correspond to braids smaller than the considered one. By
Corollary~\ref{C:Normal3}, $\Phi$-normal words are characterized by the inequalities $\ee_\rr \ge \emin_\rr$ for
$\rr < \pp$, and \eqref{E:Rank} follows.
\end{proof}

For instance, we saw in Lemma~\ref{L:dddd} that the $3$-splitting of~$\D_3^\dd$ is the
length~$\dd+\nobreak 2$ sequence $(\ss1, \s_1^2, \Ldots, \s_1^2, \ss1, \s_1^\dd)$. Proposition~\ref{P:Rank}
shows that,  for each~$\dd$, the rank of~$\D_3^\dd$ in~$(\BB3, \ls\nobreak)$ is the ordinal~$\om^{\dd+1} +
\dd$: only the initial~$1$ and the final~$\dd$ contribute here, as all intermediate exponents have the minimal legal
value~$\emin_\rr$.

\begin{ques}
\label{Q:Rank}
Does there exist a similar explicit formula for the rank of an arbitrary positive braid in the
well-ordering~$(\BB\infty, \ls)$?
\end{ques}

We refer to~\cite{But} for partial results about Question~\ref{Q:Rank}, and to~\cite{Dhq} for further applications,
consisting of unprovability statements involving braids.

\subsection{Artin--Tits monoids and other Garside monoids}

We proved in Section~\ref{S:Iter} that
$\MMM$-decompositions exist in every locally right Garside
monoid~$\MM$ in which enough closed submonoids exist.
This is in particular the case for every Artin--Tits monoid
with respect to the standard set of generators~$\SS$, as every
subset of~$\SS$ generates a closed submonoid that is closed.
Thus, dense atomic coverings exist for every Artin--Tits monoid~$\MM$, and each of them leads to
$\MMM$-decompositions similar to those of Section~\ref{S:Flip}. Then, we can
adapt Section~\ref{S:FlipOrder} and define a linear
ordering~$<_{\MMM}$ of~$\MM$ using the
$\ShortLex$-ordering on $\MMM$-decompositions.

\begin{ques}
\label{Q:Artin}
Let $\MM$ be an Artin--Tits monoid. Is any of the 
linear orders~$<_{\MMM}$ invariant under left multiplication?
\end{ques}

In type~$A_\nn$, \ie, if $\MM$ is a braid monoid,
Corollary~\ref{C:Compat2} provides a positive  answer. But the
proof depends on the connection between the
orders~$\lf$ and~$\ls$ and it is quite specific. More general positive results would presumably entail a direct
proof in the case of braids,
\ie, an answer to Question~\ref{Q:Final}$(ii)$.

Another possible extension of the current approach consists in
addressing braids again, but in connection with
other monoids. Laver's proof of Theorem~\ref{T:Order}$(ii)$
implies that the restriction of~$\ls$ to any finitely generated
submonoid of~$B_\infty$ generated by conjugates of
the~$\ss\ii$'s is a well-ordering. In particular, the
restriction of~$\ls$ to the dual braid monoids
of~\cite{BKL} is a well-ordering. The latter are Garside monoids,
and they are directly relevant for our approach. Natural analogs to the $\Phi$-normal forms exist, and
investigating their connection with the braid ordering is an
obvious task, recently achieved by J.\,Fromentin in~\cite{Fro}. It turns out that the dual framework is more
suitable than the standard one, in that a positive answer to the counterpart of Question~\ref{Q:Final} can be given,
with a direct proof that requires no transfinite induction.

\subsection{Geometric and dynamic properties}

Not much is known about the $\Phi$-normal form
of braids. As every braid admits a canonical
decomposition as a fraction~$\xx\yy\inv$ with $\xx,\yy$ in~$\BB\infty$ with no common right divisor, we can
extend the $\Phi$-normal form of~$\BB\infty$ into a unique normal form
on~$B_\infty$. Experiments suggest that the behaviour of this
normal form is rather different from that of the greedy
normal form, and many questions arise about the geometry
it induces on the Cayley graph of~$B_\nn$. In particular, it is
natural to ask for a possible associated automatic structure.
The answer seems to be negative.

\begin{prop}
$(i)$ For each~$\nn$, the set of all (positive) $\Phi$-normal
$\nn$-strand braid words is rational, \ie, recognized by a finite state
automaton.

$(ii)$ For $\nn\ge3$, $\Phi$-normal words do not satisfy the Fellow Traveler Property \cite{Eps}
with respect to multiplication on the right.
\end{prop}

\begin{proof}[Proof (sketch)]
$(i)$ By Proposition~\ref{P:NormalBraid}, a positive $\nn$-strand braid word~$\ww$ is $\Phi$-normal if and only if
each letter occurring in~$\ww$ is the smallest~$\ss\ii$ that right divides the braid represented by the
prefix finishing at that letter, with respect to an ordering of~$\{\ss1, \Ldots, \ss\nno\}$ that depends on the suffix
starting at that letter (actually at the next one). It is easy to construct an automaton that, when reading a braid word,
returns the set of all~$\ss\ii$ that right divide the braid represented by that word. Similarly, it is easy to
construct a reversed automaton that, reading a braid word from the right, returns the local ordering 
of~$\{\ss1, \Ldots,
\ss\nno\}$ that is involved in the above construction. Standard techniques from the theory of automata enable one
to mix both constructions, and to build an automaton that recognizes the family of all $\Phi$-normal $\nn$-strand
braid words.

$(ii)$ For odd (\resp even) $\dd\ge0$, the $\Phi$-normal
form of~$\D_3^\dd$ is~$\uu_\dd = \DDhat3\dd
\s_1^\dd$, while that of~$\D_3^\dd \ss2$ is~$\vv_\dd =
\ss1 \uu_\dd$ (\resp $\ss2 \uu_\dd$)---as $\DDhat\nn\dd$ is a braid that admits a unique positive word
representative, there is no danger here in using the same notation for the word and the braid. For
$\ell = 1,
\Ldots, 3\dd+1$, the successive distances between the length~$\ell$
prefixes of~$\uu_\dd$ and~$\vv_\dd$ turn out to be

$0, 2, 4, 4, 6, 6, \Ldots, 2(\dd-1), 2(\dd-1), 2\dd, 2\dd, 2\dd,
2(\dd-1), \Ldots, 6, 4, 2, 1$.\\ There is no uniform upper bound for
the above distances, hence the $\kk$-Fellow Traveler Property
fails for every~$\kk$.
\end{proof}

Investigating the dynamical properties of the $\Phi$-normal
form along the lines addressed in~\cite{BNV, NeV, MaM, Mai,
MaM2} is also a natural task. The generic problem is to study
growth and stabilization in random walks through~$B_\nn$
or~$\BB\nn$:  one compares the successive normal forms,
typically looking at whether the first factors become
eventually constant. Each new normal form induces a new
problem. Let
$\bb(\xx)$ denote the $\nn$-breadth of~$\xx$, and
$\cc_\rr(\xx)$ denote the $\rr$th entry, starting from the
right, in the $\nn$-splitting of~$\xx$. 

\begin{ques}
Let $\XX$ be the random walk
through~$\BB\nn$ defined by $\XX_{\kk+1}=\nobreak \ss\ii\, \XX_\kk$
with $\ii$ equidistributed in $\{1,\Ldots, \nn-1\}$. What are the
distributions of $\frac1\kk\bb(\XX_\kk)$ and
$\frac1\kk\lg{\cc_\rr(\XX_\kk)}$ for each fixed~$\rr$?
\end{ques}

Experiments suggest that the length of~$\cc_0(\XX_\kk)$
might grow like~$\kk/(\nn+2)$, while
$\cc_\rr(\XX_\kk)$ with $\rr\ge1$ tends to stabilize
to~$\ddd\nno\ss1$, of constant length,
and $\bb(\XX_\kk)$ might be connected with~$\sqrt{\kk}$.

\end{document}